\documentclass[reqno]{amsart}

\def\E{\ifmmode{\mathbb E}\else{$\mathbb E$}\fi} 
\def\N{\ifmmode{\mathbb N}\else{$\mathbb N$}\fi} 
\def\R{\ifmmode{\mathbb R}\else{$\mathbb R$}\fi} 
\def\Q{\ifmmode{\mathbb Q}\else{$\mathbb Q$}\fi} 
\def\C{\ifmmode{\mathbb C}\else{$\mathbb C$}\fi} 
\def\H{\ifmmode{\mathbb H}\else{$\mathbb H$}\fi} 
\def\Z{\ifmmode{\mathbb Z}\else{$\mathbb Z$}\fi} 
\def\P{\ifmmode{\mathbb P}\else{$\mathbb P$}\fi} 
\def\T{\ifmmode{\mathbb T}\else{$\mathbb T$}\fi} 
\def\SS{\ifmmode{\mathbb S}\else{$\mathbb S$}\fi} 
\def\DD{\ifmmode{\mathbb D}\else{$\mathbb D$}\fi} 

\renewcommand{\a}{\alpha}
\renewcommand{\b}{\beta}

\newcommand{\e}{\varepsilon}
\newcommand{\g}{\gamma}
\newcommand{\G}{\Gamma}

\newcommand{\s}{\sigma}

\renewcommand{\t}{\tau}

\renewcommand{\O}{\Omega}
\renewcommand{\o}{\omega}

\newcommand{\del}{\partial}

\newcommand{\ben}{\begin{enumerate}}
\newcommand{\een}{\end{enumerate}}
\newcommand{\be}{\begin{equation}}
\newcommand{\ee}{\end{equation}}
\newcommand{\bea}{\begin{eqnarray}}
\newcommand{\eea}{\end{eqnarray}}
\newcommand{\bc}{\begin{center}}
\newcommand{\ec}{\end{center}}

\theoremstyle{theorem}
\newtheorem{thm}{Theorem}[section]
\newtheorem{cor}[thm]{Corollary}
\newtheorem{lem}[thm]{Lemma}
\newtheorem{prop}[thm]{Proposition}

\theoremstyle{definition}
\newtheorem{defn}{Definition}[section]
\newtheorem{rem}[defn]{Remark}

\newtheorem{ques}[defn]{Question}

\newtheorem{exm}[defn]{Example}

\numberwithin{equation}{section}

\def\R{{\mathbb R}}

\def\E{{\mathbb E}}
\def\Z{{\mathbb Z}}
\def\C{{\mathbb C}}
\def\R{{\mathbb R}}

\def\N{{\mathbb N}}

\def\HH{{\mathcal H}}

\def\SS{{\mathcal S}}
\def\LL{{\mathcal L}}

\def\DD{{\mathcal D}}

\def\JJ{{\mathcal J}}
\def\AA{{\mathcal A}}
\def\PP{{\mathcal P}}

\def\FF{{\mathcal F}}

\def\11{{\mathbb I}}

\def\C{\mathbb{C}}
\def\Z{\mathbb{Z}}

\def\T{\mathbb{T}}
\def\L{\mathbb{L}}

\def\Q{\mathbb{Q}}

\def\E{\ifmmode{\mathbb E}\else{$\mathbb E$}\fi} 
\def\N{\ifmmode{\mathbb N}\else{$\mathbb N$}\fi} 
\def\R{\ifmmode{\mathbb R}\else{$\mathbb R$}\fi} 
\def\Q{\ifmmode{\mathbb Q}\else{$\mathbb Q$}\fi} 
\def\C{\ifmmode{\mathbb C}\else{$\mathbb C$}\fi} 
\def\H{\ifmmode{\mathbb H}\else{$\mathbb H$}\fi} 
\def\Z{\ifmmode{\mathbb Z}\else{$\mathbb Z$}\fi} 
\def\P{\ifmmode{\mathbb P}\else{$\mathbb P$}\fi} 
\def\SS{\ifmmode{\mathbb S}\else{$\mathbb S$}\fi} 
\def\DD{\ifmmode{\mathbb D}\else{$\mathbb D$}\fi} 

\def\R{{\mathbb R}}

\def\E{{\mathbb E}}
\def\Z{{\mathbb Z}}
\def\C{{\mathbb C}}
\def\R{{\mathbb R}}

\def\N{{\mathbb N}}
\def\LL{{\mathcal L}}

\def\JJ{{\mathcal J}}
\def\FF{{\mathcal F}}

\def\a{\alpha}
\def\b{\beta}

\def\e{\varepsilon} 
  
\def\g{\gamma}  \def\G{\Gamma}

\def\o{\omega}  \def\O{\Omega}
\def\p{\psi}  \def\P{\Psi}
\def\s{\sigma}  
  
\def\t{\tau}

\def\CC{{\mathcal C}}

\def\CF{{\mathcal F}}

\def\CQ{{\mathcal Q}}

\def\rd{\partial}

\def\darr#1{\raise1.5ex\hbox{$\leftrightarrow$}
\mkern-16.5mu #1}

\def\Fr#1#2{{#1\over#2}}

\def\roughly#1{\raise.3ex\hbox{$#1$\kern-.75em
\lower1ex\hbox{$\sim$}}}

\def\opname#1{\mathop{\kern0pt{\rm #1}}\nolimits}

\def\dim{\opname{dim}}

\def\pr{\prime}

\begin{document}
\quad \vskip1.375truein

\def\mq{\mathfrak{q}}
\def\mp{\mathfrak{p}}
\def\mH{\mathfrak{H}}
\def\mh{\mathfrak{h}}
\def\ma{\mathfrak{a}}
\def\ms{\mathfrak{s}}
\def\mm{\mathfrak{m}}
\def\mn{\mathfrak{n}}

\def\Hoch{{\tt Hoch}}
\def\mt{\mathfrak{t}}
\def\ml{\mathfrak{l}}
\def\mT{\mathfrak{T}}
\def\mL{\mathfrak{L}}
\def\mg{\mathfrak{g}}
\def\md{\mathfrak{d}}

\title[Deformations of Coisotropic Submanifolds]
{Deformations of Coisotropic Submanifolds \\
and Strong Homotopy Lie Algebroids}

\author{Yong-Geun Oh
\and Jae-Suk Park}
\thanks{The first named author is partially supported
by the NSF grant \#DMS 0203593,
a grant of the 2000 Korean Young Scientist Prize, and the Vilas
Research Award of the University of Wisconsin}

\address{
Department of Mathematics, University of Wisconsin, Madison, WI
53706 \& Korea Institute for Advanced Study, 207-43
Cheongryangri-dong Dongdaemun-gu, Seoul 130-012, KOREA,
oh@math.wisc.edu; \\
Korea Institute for Advanced Study,  207-43
Cheongryangri-dong Dongdaemun-gu, Seoul 130-012, KOREA,
jaesuk@kias.re.kr}

\begin{abstract}
In this paper, we  study deformations of coisotropic submanifolds
in a symplectic manifold. First we derive the equation that
governs $C^\infty$ deformations of coisotropic submanifolds and
define the corresponding $C^\infty$-moduli space of coisotropic
submanifolds modulo the Hamiltonian isotopies. This is a
non-commutative and non-linear generalization of the well-known
description of the local deformation space of Lagrangian
submanifolds as the set of graphs of {\it closed} one forms in the
Darboux-Weinstein chart of a given Lagrangian submanifold. We then
introduce the notion of  {\it strong homotopy Lie algebroid} (or
{\it $L_\infty$-algebroid}) and
associate a canonical isomorphism class of strong homotopy Lie algebroids
to each pre-symplectic manifold $(Y,\omega)$ and identify the formal
deformation space of coisotropic embeddings into a symplectic
manifold in terms of this strong homotopy Lie algebroid. The
formal moduli space then is provided by the gauge equivalence
classes of solutions of a version of the {\it Maurer-Cartan
equation} (or the {\it master equation}) of the strong homotopy
Lie algebroid, and plays  the role of the classical part of the
moduli space of quantum deformation space of coisotropic
$A$-branes. We provide a criterion for the
unobstructedness of the deformation problem and analyze a family
of examples that illustrates that this deformation problem is
obstructed in general and heavily depends on the
geometry and dynamics of the null foliation.
\end{abstract}

\maketitle

\bigskip

\tableofcontents
\section{Introduction}\label{intro}

The well-known Darboux-Weinstein theorem [We1] states that a
neighborhood of any Lagrangian submanifold $L$ of any symplectic
manifold $(X,\omega_X)$ (with real dimensions $2n$)
is diffeomorphic to a neighborhood
of the zero section of the cotangent bundle $T^*L$ with the
standard symplectic form
\begin{equation}\label{eq:omega0}
\omega = -d\theta, \quad \theta = \sum_{i =1}^n p_idq^i
\end{equation}
where $\theta$ is the canonical one form defined by
$$
\theta_p(\xi) = p(T\pi(\xi))
$$
for $p \in T^*L$, $\xi \in T_p(T^*L)$ and $\pi:T^*L \to L$ is the
canonical projection. Furthermore it is also well-known that for
any section $\alpha: L \to T^*L$ of $T^*L$, i.e.,
for any one form $\alpha$ on $L$,
we have the identity
\begin{equation}\label{eq:alpha}
\alpha^*\theta = \alpha.
\end{equation}
{}From this it follows that any Lagrangian submanifold $C^1$-close
to the zero section is the graph of a {\it closed} one form. And
two such Lagrangian graphs of $\alpha$ and $\alpha'$ are
Hamiltonian isotopic if and only if $\alpha - \alpha' = d\beta$,
i.e, exact. Therefore the (local) moduli space near the given
Lagrangian submanifold $L \subset (X,\omega_X)$ is diffeomorphic
to a neighborhood of $0 \in H^1(L;\R)$.
In particular the
local moduli problem of the Lagrangian submanifold up to
Hamiltonian isotopy is {\it linear} and {\it commutative}. It
depends only on the manifold $L$ but is independent of where the
abstract manifold $L$ is embedded into as a Lagrangian
submanifold.

The main purpose of the present paper is to  describe the moduli
space of {\it coisotropic} submanifolds modulo the Hamiltonian
isotopy in a symplectic manifold, and its formal counterpart.
Recall that a submanifold $i: Y
\hookrightarrow (X,\omega_X)$ is called {\it coisotropic} if the
symplectic orthogonal $(TY)^\omega$ satisfies
$$
(TY)^\omega \subset TY
$$
and has constant rank. Then the pull-back $\omega=i^*\omega_X$ is a
closed two form with constant rank $2k$ on $Y$. Such a pair
$(Y,\omega)$ in general is called a {\it pre-symplectic} manifold.
Gotay [Go] proved that any given pre-symplectic manifold
$(Y,\omega)$ can be embedded into a symplectic manifold
$(X,\omega_X)$ as a coisotropic submanifold so that $i^*\omega_X =
\omega$. Furthermore the neighborhoods of any two such embeddings
are diffeomorphic regardless of $(X,\omega_X)$. In fact, a
neighborhood of any coisotropic submanifold $Y \subset
(X,\omega_X)$ is locally diffeomorphic to a neighborhood of the zero section
$o_{E^*} \cong Y$ of the bundle
$$
E^* \to Y; \quad E = (TY)^\omega \subset TY
$$
with a symplectic form $\omega_{E^*}$ defined in the neighborhood.
One drawback of Gotay's theorem [Go] is that it does not provide
the symplectic form $\omega_{E^*}$ explicitly. This hindered our
study of the current deformation problem in the beginning.  We
would like to recall that in the Lagrangian case we have the
explicit canonical symplectic form $\omega_X = -d\theta$ on the
model space $T^*L$ which enables one to study the moduli problem
explicitly. See section 3 for more comments on this point.

We first note that $E=T\FF$ is the tangent bundle of the {\it null
foliation} of $(Y,\omega)$ and so $E^*$ is the {\it cotangent
bundle} $T^*\FF$ of the {\it foliation} $\FF$. Then we introduce an
explicitly given one form $\theta_G$ on $T^*\FF$ which is the analogue to the
canonical one form $\theta$ on the cotangent bundle of the
{\it manifold}. This one form will, however, depend on the choice of the splitting
$$
TY = G \oplus E.
$$
This splitting can be regarded as a ``connection'' of the
``$E$-bundle'' $TY \to Y/\sim$ where $Y/\sim$ is the space of leaves
of the null foliation of $Y$. Using this  one form $\theta_G$ we
will write down the explicit symplectic form $\omega_{E^*}$ in
terms of $(Y,\omega)$ and $\theta_G$. And then we will write
down the defining equation for a section $s$ of $\pi: E^* \to Y$
whose graph is to be coisotropic in $(E^*, \omega_{E^*})$ (in a
neighborhood of the zero section). The moduli problem for
this general coisotropic case is {\it non-commutative} and {\it
fully non-linear}. The moduli problem becomes trivial when the
foliation is one dimensional (i.e., the case of hypersurfaces) and
becomes {\it quadratic} when the null-foliation also allows a
transverse foliation. The moduli problem up to Hamiltonian isotropy
is obstructed in general.

In fact, it turns out that the equation for the formal moduli
problem involves the structure of a {\it strong homotopy
Lie-algebra}
$$
(\Omega^\bullet(\FF), \frak m),
$$
where $\frak m = \{\frak m_i\}_{i=1}^\ell$, $\Omega^\bullet(\FF):=
\oplus_{\ell = 0}^{n-k} \Gamma(\Lambda^\ell(E^*))$: In general, we
call a Lie algebroid $E \to Y$ a {\it strong homotopy Lie
algebroid} if its associated graded group $\Omega^\bullet(E): =
\oplus \Gamma(\Lambda^\ell(E^*))$ has the structure of strong
homotopy Lie algebra (or $L_\infty$-algebra) where
$\Lambda^\ell(E^*)$ is the set of $\ell$-wedge product of $E^*$ so
that $\frak m_1$ is the $E$-differential induced by the Lie
algebroid structure on $E \to Y$.  In our case, $\frak m_1$ is the
leafwise differential associated to the null foliation and $\frak
m_2= \{\cdot, \cdot\}$ a (graded) bilinear map, which is an analog to
a Poisson bracket but not
satisfying the Jacobi identity whose failure is then measured by
$\frak m_3$ and so on. We refer to later sections, especially
section \ref{sec:algebroid} for precise details.

The main purpose of the present paper is to unravel the geometric
and algebraic structures that govern the deformation problem of
coisotropic submanifolds. It turns out that even setting up the
proper framework for the study of this deformation problem
requires definitions of many new geometric notions which have not
appeared in the literature before, as far as we know. We have
discovered them first via tensor calculations using coordinates
and then associated the relevant geometric structures to them.
Partly because of this and also because the coordinate
calculations are better suited for the super-calculus in the
Batalin-Vilkovisky formalism as outlined in the appendix, we
prefer to carry out the coordinate calculations first and then
provide the corresponding invariant descriptions. However the
latter can be done only after one develops an appropriate
invariant calculus in the properly formulated geometric framework.
It appears that a full systematic study of those should be a theme
of separate study which is not our main interest in the present
paper. In the present paper, we refrain  from developing
the full invariant calculus, but carry out only those essential
for a self-contained description of our deformation problem of
coisotropic submanifolds. A more thorough study of the invariant
calculus in the context of general foliation theory will be
made elsewhere.

The present work is a mathematical spin-off of our paper [OP] in
preparation in which we provide an off-shell description of
topological open/closed $A$-strings and $A$-branes on symplectic
manifold (``bulk''). In [OP], we have derived in the natural
framework of the Batalin-Vilkovisky formalism that the set of
$A$-branes of the topological open $\sigma$-model on the
symplectic manifold is the set of coisotropic submanifolds (at
least the bosonic part thereof), which we also outline in the
appendix of the present paper. A possible significance of
coisotropic branes in homological mirror symmetry was first
observed by Kapustin and Orlov [KaOr]. After this paper was
originally submitted and circulated, a paper [CF2] by Catteneo and
Felder has appeared in
which a similar discussion on the role of coisotropic submanifolds
in the context of Poisson manifolds is carried out in a same
spirit to ours given in Appendix in relation to the open string
$\sigma$-model. There has also appeared a paper [R] by W.-D. Ruan
in which Ruan studies a deformation problem of a restricted
class of coisotropic submanifolds which he calls
{\it integral coisotropic submanifolds} :  a coisotropic submanifold
$Y \subset (X,\omega)$ is said to be integral if its associated null
foliation is defined by a smooth fibration $\pi: Y \to S$ over a
smooth (Hausdorff) manifold $S$. He proved that deformations inside the class
of integral coisotropic submanifolds are unobstructed and the moduli
space thereof modulo Hamiltonian diffeomorphisms is smooth and finite
dimensional [R].

This research was initiated while both authors were visiting the Korea
Institute for Advanced Study in the winter of 2002. We thank KIAS for its financial
support and excellent research environment. The first named author
thanks M. Zambon for explaining the example mentioned in section
\ref{sec:moduli} in an IPAM conference in April, 2003 and for asking some
questions on the preliminary version of this paper. He himself
independently studied neighborhoods of coisotropic submanifolds
[Za] and found an example of a coisotropic submanifold which
illustrates the fact that
$C^1$-close coisotropic submanifolds of the same nullity do not
form a Fr\'echet manifold in general.

The second named author thanks the mathematics department of
POSTECH, especially to his host B.\ Kim, for the financial support
and excellent research environment, and to the mathematics
department of University of  the Wisconsin-Madison for its hospitality
during his visit. We also thank J. Stasheff for some historical
comments on the $L_\infty$-structure and the Gerstenhaber bracket,
and I. Vaisman for attracting our attention to his paper [V]
after the paper was submitted. The first named author also
thanks N. Kieserman for pointing out and correcting
some inaccuracy in our coordinate calculations.
Last, but not least, we are very much
indebted to the referees for making many valuable suggestions and
corrections which have led to much improvement and clarification of
the contents and presentation of the paper.

\section{Geometry of coisotropic Grassmanians}
\label{Coiso}

In this section, we will summarize some linear algebraic facts on
the coisotropic subspace $C$ (with real dimensions $n+k$ where $0\leq k\leq n$)
in $\C^n$ with respect to the standard
symplectic form $\omega = \omega_{0,n}$.  We denote by
$C^{\omega}$ the $\omega$-orthogonal complement of $C$ in
$\R^{2n}$ and by $\Gamma_k$ the set of coisotropic subspaces of
$(\R^{2n}, \omega)$. In other words,
\begin{equation}
\label{gammak} \Gamma_k = \Gamma_k(\R^{2n}, \omega) = :\{ C \in
Gr_{n+k}(\R^{2n}) \mid C^\omega \subset C\}.
\end{equation}
{}From the definition, we have the
canonical flag,
$$
0 \subset C^\omega \subset C \subset \R^{2n}
$$
for any coisotropic subspace.
We call $(C,C^\omega)$ a {\it coisotropic pair}. Combining this
with the standard complex structure on $\R^{2n} \cong \C^n$, we
have the splitting
\begin{equation}\label{eq:csplit}
C = H_C \oplus C^\omega
\end{equation}
where $H_C$ is the complex subspace of $C$.

\begin{prop}
Let $ 0 \leq k \leq n$ be fixed. The unitary group $U(n)$ acts
transitively on $\Gamma_k$. The corresponding homogeneous space is
given by
\begin{equation}
\Gamma_k \cong U(n)/U(k) \times O(n-k)
\end{equation}
where $U(k) \times O(n-k) \subset U(n)$ is the isotropy group of
the coisotropic subspace $\C^k \oplus \R^{n-k} \subset \C^n$. In
particular we have
$$
\dim \Gamma_k(\R^{2n}, \omega) = \frac{(n+3k +1)(n-k)}{2}.
$$
\end{prop}
\begin{proof}
Let $C \subset \C^n$ be a coisotropic subspace with rank $2k$ and
$C^{\omega}$ be its null space. Since $C^{\omega} \subset C$,
it follows
$$
g(C, iC^{\omega}) = 0,
$$
where $g$ is the Euclidean inner product. Now if we write
$$
H_C = \{x \in C \mid g( x, C^\omega) = 0 \},
$$
then $H_C$ is the Hermitian orthogonal complement of $C^\omega$ which is
a complex subspace. Similarly $C^\omega\oplus iC^\omega$ is also a
complex subspace. Therefore we have obtained the Hermitian
orthogonal decomposition
$$
\C^n = H_C \oplus (C^\omega \oplus iC^\omega).
$$
It then follows that there is a unitary matrix $A \in U(n)$ such
that $C = A \cdot (\C^k \oplus \R^{n-k})$. This proves that $U(n)$
acts transitively on $\Gamma_k$. Now it is enough to show that the
isotropy group of $\C^{k}\oplus \R^{n-k}$ is $U(k) \times O(n-k)
\subset U(n)$ which is obvious. This proves the proposition.
\end{proof}

Next we give a parametrization of all the coisotropic subspaces
near given $C \in \Gamma_k$. Up to the unitary change of
coordinates we may assume that $C$ is the canonical model
$$
C = \C^k \oplus \R^{n-k}.
$$
We denote the (Euclidean) orthogonal complement of $C$ by $C^\perp
= i \R^{n-k}$ which is canonically isomorphic to $(C^\omega)^*$
via the isomorphism $\widetilde\omega: \C^n \to (\C^n)^*$.
Then any nearby subspace of dimension $\dim C$ that
is transverse to $C^\perp$ can be written as the graph of the
linear map
$$
A: C \to C^\perp \cong (C^\omega)^*
$$
i.e., has the form
\begin{equation}\label{eq:ca}
C_A:=\{ (x, Ax) \in C \oplus C^\perp = \R^{2n} \mid x \in C \}.
\end{equation}

Denote $A = A_H \oplus A_I$ where
\begin{align}
A_H & : H = \C^k \to C^\perp \cong (C^\omega)^*,\nonumber\\
A_I & : C^\omega = \R^{n-k}  \to C^\perp\cong (C^\omega)^*.\nonumber
\end{align}
Note that the symplectic form $\omega$ induce the canonical
isomorphism
\begin{eqnarray*}
\widetilde\omega^H & : & \C^k \to (\C^k)^* \\
\widetilde\omega^I & : & \R^{n-k}=C^\omega
\to (C^\omega)^*\cong C^\perp = i
\R^{n-k} \end{eqnarray*}
With this identification, the symplectic form $\omega$ has the
form
\begin{equation}\label{eq:omega0}
\omega = \pi^*\omega_{0,k} + \sum_{i=1}^{n-k} dx_i\wedge dy^i,
\end{equation}
where $\pi: \C^n \to \C^k$ is the projection and $(x_1,\cdots,
x_{n-k})$ the standard coordinates of $\R^{n-k}$ and $(y^1,
\cdots, y^{n-k})$ its dual coordinates of $(\R^{n-k})^*$. We also
denote by $\pi_H: (\C^k)^* \to \C^k$ the inverse of the above
mentioned canonical isomorphism $\widetilde \omega^H$. Then we
have the following

\begin{prop} The subspace $C_A$ is coisotropic if and only
if $A_H$ and $A_I$ satisfies
\begin{equation}\label{eq:grapca}
A_I - (A_I)^* + A_H \pi_H (A_H)^* = 0.
\end{equation}
\end{prop}
\begin{proof} We need to study under what conditions on $A_H$ and
$A_I$, the relation $(C_A)^{\omega} \subset C_A$ holds and vice versa.
Let $(\xi_H, \xi_I, \xi_I^*) \in (C_A)^{\omega} \subset \C^n
=\C^k \oplus \R^{n-k}\oplus i\R^{n-k}$, i.e, let
\begin{equation}\label{eq:vhvi}
\omega((v_H, v_I, A(v_H, v_I)), (\xi_H, \xi_I, \xi_I^*)) = 0
\end{equation}
for all $(v_H, v_I) \in \C^k\oplus \R^{n-k}$. It follows from
(\ref{eq:vhvi}) and from the above identifications, we have
\begin{eqnarray}
0 & = & \omega(v_H, \xi_H)  +  A(v_H, v_I)(\xi_I) -
\xi_I^*(v_I) \nonumber \\
& = & \omega(v_H, \xi_H) + (A_H(v_H)(\xi_I) + A_I(v_I)(\xi_I)) -
\xi_I^*(v_I)
\label{eq:vhvi0}
\end{eqnarray}
for all $v_H, \, v_I$. Substituting $v_I = 0$ we get
\begin{equation}\label{eq:vi0j}
\omega(v_H,\xi_H) + A_H(v_H)(\xi_I) = 0
\end{equation}
for all $v_H \in \C^k$. With the above identification, we derive
$$
-{\widetilde \omega}_0^H(\xi_H)(v_H) + A_H^*(\xi_I)(v_H) = 0
$$
for all $v_H \in \C^k$. Therefore we derive
\begin{equation}\label{eq:vi0xi}
\xi_H= \pi_H A_H^*(\xi_I).
\end{equation}
And substituting $v_H =0$, we derive
$$
A_I(v_I)(\xi_I) - \xi_I^*(v_I) = 0
$$
for all $v_I \in \R^{n-k}$ and hence
\begin{equation}\label{eq:vh0}
A_I^*(\xi_I) = \xi_I^*.
\end{equation}
Therefore it follows from (\ref{eq:vi0xi}) and (\ref{eq:vh0}) that $C_A$ is
coisotropic if and only if
\begin{equation}\label{eq:ahj}
A_I^*(\xi_I)= A(\pi_H A_H^*(\xi_I), \xi_I) =
A_H\oplus A_I (\pi_H A_H^*(\xi_I), \xi_I)
\end{equation}
for all $\xi_I$. The latter becomes
$$
A_I - (A_I)^* + A_H\pi_H(A_H)^* = 0,
$$
which finishes the proof.
\end{proof}

\begin{rem} Note that when $k=0$, this reduces to the standard
parametrization of Lagrangian subspaces by the set of symmetric
matrices.
\end{rem}

\section{Canonical symplectic neighborhoods}
\label{sec:neighborhoods}

We first recall some basic properties of  coisotropic
submanifolds and the coisotropic neighborhood theorem [Go]. We
will mostly adopt the notations used in [Go]. First $(TY)^\omega:
= E$ defines a distribution, the so called {\it characteristic
distribution} on $Y$ which is integrable since $\omega$ is closed.
We call the corresponding foliation the {\it null foliation} on
$Y$ and denote it by $\CF$.
The null foliation carries a natural {\it transverse symplectic form}
but the space of leaves may {\it not} be a Hausdorff space in general.
This space of leaves provides a symplectic invariant of coisotropic
submanifolds up to the Hamiltonian isotopy, or equivalently an
invariant of the pre-symplectic manifold $(Y,\omega)$.
We refer to the next section for the detailed description of the geometry of
null foliation.

We now consider the dual bundle $\pi: E^* \to Y$ of $E$.
The bundle
$TE^*|_Y$ where $Y \subset E^*$ is the zero section of $E^*$
carries the canonical decomposition
$$
TE^*|_Y = TY \oplus E^*.
$$
It is easy to check that the canonical isomorphism
$$
\widetilde \omega: TX \to T^*X
$$
maps $TY^\omega$ to the conormal $N^*Y \subset T^*X$, and induces
an isomorphism between $NY = TX/TY$ and $E^*$. In the standard
notation in the foliation theory, $E$ and $E^*$ are denoted by
$T\FF$ and $T^*\FF$ and called the tangent bundle (respectively
cotangent bundle) of the foliation $\FF$.

Following Gotay [Go],  we choose a splitting
\begin{equation}\label{eq:splitting}
TY = G \oplus E, \quad E = (TY)^\omega.
\end{equation}
Using this splitting, we can write a symplectic form on a
neighborhood of the zero section $Y \hookrightarrow E^*$ in the
following way (see [V]). We denote by
$$
p_G: TY \to E
$$
the projection to $E$ along $G$ in the splitting
(\ref{eq:splitting}). We have the bundle map
$$
TE^* \stackrel{T\pi}\longrightarrow TY \stackrel{p_G}
\longrightarrow E.
$$
Let $\alpha \in E^*$ and $\xi \in T_\alpha E^*$. We define the one
form $\theta_G$ on $E^*$ by its value
\begin{equation}\label{eq:thetag}
\theta_{G,\alpha}(\xi): = \alpha(p_G\circ T\pi(\xi))
\end{equation}
at each $\alpha \in E^*$.
Then we define the closed (indeed exact) two form on $E^*$ by
$$
- d\theta_G.
$$
It is easy to see that the closed two form
\begin{equation}\label{eq:omega*}
\omega_{E^*}:= \pi^*\omega - d\theta_G
\end{equation}
is non-degenerate in a neighborhood $U \subset E^* $ of the zero
section (See the coordinate expression (\ref{eq:dthetaG}) of $d
\theta_G$ and $\omega_U$). We denote the restriction of
$\omega_{E^*}$ by $\omega_U$. Then the pair $(U, \omega_U)$
provides an explicit normal form of the symplectic neighborhood of
the pair $(Y,\omega)$ which depends only on $(Y,\omega)$ and the
splitting (\ref{eq:splitting}). By Weinstein's uniqueness theorem
[We1], this normal form is unique up to diffeomorphism. We call
the pair $(U,\omega_U)$ a {\it (canonical) symplectic thickening}
of the pre-symplectic manifold $(Y,\omega_Y)$. We refer to Corollary
\ref{diffeo} for the precise statement on the uniqueness.

Here is the coisotropic analog to (1.2) whose proof we omit.

\begin{lem}\label{coiso(alpha)} For any section of $s: Y \to U \subset E^*$,
we have the identity
$$
s^*\theta_G = p_G^* s
$$
where $p_G^*: E^* \to T^*Y$ is the adjoint to the projection
$p_G: TY \to E$.
\end{lem}

We next introduce morphisms between pre-symplectic manifolds
and automorphisms of  $(Y,\omega)$.

\begin{defn}\label{morphism} Let $(Y,\omega)$ and
$(Y^\prime, \omega^\prime)$ be two pre-symplectic manifolds. A
diffeomorphism $\phi: Y\to Y'$ is called {\it pre-symplectic} if
$\phi^*\omega^\prime = \omega$.
\end{defn}

Lemma \ref{coiso(alpha)} immediately implies the following

\begin{cor}\label{closed}
Suppose that the graph $i: \mbox{Graph } s \hookrightarrow U$ of  a section
$s: Y \to U$ is a coisotropic submanifold i.e., $i^*\omega_U$
induces a pre-symplectic structure. Then the map $s: (Y,\omega)
\to (\mbox{Graph} s, i^*\omega_U)$ is a pre-symplectic
diffeomorphism if and only if $p_G^*s$ defines a closed one form
on $Y$.
\end{cor}
\begin{proof}
We first note that
\begin{equation}\label{eq:s*i*omegaU}
s^*(i^*\omega_U) = s^*(\omega_U) =
s^*\pi^*\omega - s^*(d\theta_G) = \omega - ds^*\theta_G
= \omega - d(p_G^*s).
\end{equation}
The corollary immediately follows from this.
\end{proof}

We would like to emphasize that unlike Lagrangian submanifolds for
which {\it there is no intrinsic structure}, coisotropic
submanifolds carry an intrinsic geometric structure, the
pre-symplectic form $\omega$.  Therefore the hypothesis in this
corollary is not an automatic proposition even when $p_G^*s$ is a
closed one form on $Y$. In the next several sections, we will
provide a description of the condition under which  the graph of
$s: Y \to U$ becomes coisotropic.

However, (\ref{eq:s*i*omegaU}) shows that up to a diffeomorphism
$s\circ i: Y \to \mbox{Graph}(s)$, the pre-symplectic form of any
nearby coisotropic submanifold is cohomologous to the given
$\omega$ on $Y$. Therefore all the pre-symplectic structures that
occur in the study of local deformations of a coisotropic
submanifold is special in that they are cohomologous to one
another up to the pull-back by diffeomorphisms.

The following question seems to be an interesting nontrivial
question to ask in general.

\begin{ques}\label{converse}
Let $(Y,\omega), \, (Y,\omega')$ be pre-symplectic structures of
the same rank such that
$$
[\omega] = [\omega'] \quad \mbox{in } H^2(Y, \R),
$$
i.e., $\omega' - \omega = d\theta$ for some one form on $Y$.
Suppose that $\theta$ is sufficiently $C^\infty$ small. Are they
diffeomorphic or can they be connected by an isotopy of
pre-symplectic forms of constant rank? This question is
closely related to the question whether the set of
coisotropic submanifolds of constant rank is locally
path-connected (see [Za]).
\end{ques}

\begin{defn}\label{vectorfield}
Let $(Y,\omega)$ be a pre-symplectic manifold. A vector field $\xi$
on $Y$ is called or {\it locally
pre-Hamiltonian} if $\LL_\xi(\omega) = d(\xi
\rfloor \omega) = 0$. We call $\xi$ (globally) {\it pre-Hamiltonian} if the
one form $\xi \rfloor \omega$ is exact. We call the
diffeomorphisms generated by $\xi$ {\it locally pre-Hamiltonian}
(respectively {\it pre-Hamiltonian}) diffeomorphisms.
We call {\it pre-symplectic} any diffeomorphism $\phi$ that
satisfy $\phi^*\omega = \omega$.
\end{defn}

We denote by $\PP\HH am(Y,\omega)$ the set of pre-Hamiltonian
diffeomorphisms and by $\PP Symp(Y,\omega)$ the set of pre-symplectic
diffeomorphisms of $(Y,\omega)$.  According to our definitions, the set of
locally pre-Hamiltonian diffeomorphisms is the identity
component of $\PP Symp(Y,\omega)$. Therefore we will denote
the latter by $\PP Symp_0(Y,\omega)$.

The following will be important in our formulation of the moduli
problem of coisotropic submanifolds later. We will give its proof
in section \ref{sec:Ham}.

\begin{thm}\label{extension}
Any locally pre-Hamiltonian (respectively, pre-Hamiltonian) vector
field $\xi$ on a pre-symplectic manifold $(Y,\omega)$ can be
extended to a locally Hamiltonian (respectively, Hamiltonian)
vector field on the thickening $(U,\omega_U)$.
\end{thm}

\section{Leaf space connection and curvature}
\label{sec:splittings}

However the one-form $\theta_G$ defined in the previous section
depends on the splitting (\ref{eq:splitting}).  We now
describe this dependence more systematically. In
this section and the next, we will study the {\it intrinsic}
geometry of pre-symplectic manifold $(Y,\omega)$ and
 the {\it extrinsic} geometry of its symplectic thickening
in section \ref{sec:connections}.
We first need to develop some invariant calculus
``over the leaf space'' which will play an important role in our
study of the deformation problem of coisotropic submanifolds
later. This section applies to {\it any foliation}, not just to
our null foliation. As far as we know, this calculus has not been
introduced in the literature yet.

Let $\FF$ be an arbitrary foliation on a smooth manifold $Y$.
Following the standard notations in the
foliation theory, we define the normal bundle $N\FF$ and conormal
bundle $N^*\FF$ of the foliation $\FF$ by
$$
N_y\FF:= T_yY/E_y, \quad N^*_y\FF: = (T_y/E_y)^* \cong E_y^\circ
\subset T_y^*Y.
$$
In this vein, we will denote $E=T\FF$ and $E^*=T^*\FF$
respectively, whenever it makes our discussion more transparent.
We have the natural exact sequences
\begin{eqnarray}\label{eq:exactseq}
0 & \to &  T\FF \to TY \to N\FF \to 0, \\
0 &\leftarrow & T^*\FF \leftarrow T^*Y \leftarrow N^*\FF \leftarrow 0.
\end{eqnarray}
The choice of splitting $TY = G\oplus T\FF$ may be regarded as a
``connection'' of the ``$E$-bundle'' $TY \to Y /\sim$ where $Y/\sim$
is the space of leaves of the foliation on $Y$. Note that
$Y/\sim$ is {\it not} Hausdorff in general. We will indeed call
a choice of
splitting {\it a leaf space connection of $\FF$} in general.

We can also describe the splitting in a more invariant way as
follows: Consider bundle maps $\Pi: TY \to TY$ that satisfy
$$
\Pi^2_x = \Pi_x, \, \operatorname{im }\Pi_x = T_x\FF
$$
at every point of $Y$, and denote the set of such projections by
$$
\AA_E(TY) \subset \Gamma(Hom(TY,TY)) = \Omega^1_1(Y).
$$
There is a one-one correspondence between the choice of splittings
(\ref{eq:splitting}) and the set $\AA_E(TY)$ provided by the correspondence
$$
\Pi \leftrightarrow G:= \ker \Pi.
$$
If necessary, we will denote by $\Pi_G$ the element with $\ker \Pi
=G$ and by $G_\Pi$ the complement to $E$ determined by $\Pi$. We
will use either of the two descriptions, whichever is more
convenient depending on the circumstances.

The following is easy to see by using the isomorphism
$\pi_{G}: G \to N\FF$ where $\pi_{G}$ is the restriction
to $G$ of the natural projection $\pi_{\Pi}: TY \to N\FF$.
We omit its proof.

\begin{lem}\label{splittings}
The space of splittings (\ref{eq:splitting}) is an infinite
dimensional (Frechet) manifold modelled by
$$
\Gamma(Hom(G_0,T\FF)) \cong \Gamma(N^*\FF\otimes T\FF):
$$
for any reference choice $\Pi_0$, and for other $\Pi$, we have
\begin{eqnarray*}
G_\Pi & = & \{\eta  \oplus B_{\Pi_0\Pi}\circ \pi_{\Pi_0}(\eta) \in
TY
\mid \eta \in G_0, B_{\Pi_0\Pi} \in \Gamma(Hom(N\FF, T\FF)) \} \\
& = & \{\pi_{\Pi}^{-1}(y) \oplus B_{\Pi_0\Pi}(y) \in TY \mid y
\in N\FF, B_{\Pi_0\Pi} \in \Gamma(Hom(N\FF, T\FF)) \}
\end{eqnarray*}
Moreover, it is weakly contractible.
\end{lem}

Next we introduce the analogue of ``curvature'' of the above
``connection''. To define this, we recall some basic facts about
the {\it foliation coordinates}. We can choose coordinates on $Y$
adapted to the foliation in the following way. Since the
distribution $E$ is integrable, the Frobenius theorem provides
coordinates
$$
(y^1, \cdots, y^{\ell}, y^{\ell+1}, \cdots, y^{m})
$$
on an open subset $V \subset Y$, such that the plaques of $V$ are
given by the equation
\begin{equation}\label{eq:leaves}
y^1 = c^1, \cdots, y^{\ell} = c^{\ell}, \quad \mbox{$c^i$'s constant}.
\end{equation}
In particular, we have
\begin{equation}\label{eq:Esubx}
E_x = T_x\FF = \mbox{span }\Big \{\frac{\del}{\del y^{\ell+1}}, \cdots,
\frac{\del}{\del y^{m}}\Big \}.
\end{equation}
We denote
$$
q^\alpha = y^{\ell +\alpha}, \quad 1 \leq \alpha \leq m- \ell.
$$
For the given splitting $TY = G \oplus E$, we can write
\begin{equation}\label{eq:Gsubx}
G_x = \mbox{span }\Big\{\frac{\del}{\del y^i} + \sum_{\alpha
=1}^{m-\ell} R_i^\alpha \frac{\del}{\del q^\alpha}\Big\}_{1 \leq i
\leq \ell}
\end{equation}
for some $R_i^\alpha$'s, which are uniquely determined by the
splitting and the given coordinates. Here $R_i^\alpha$'s can be
regarded as the ``Christoffel symbols'' for the ``connection''
$\Pi$.

{}From now on, we will use the summation convention for repeated
indices, whenever there is no danger of confusion.
Roman indices run over $1, \cdots, \ell$, and
Greek ones over $1, \cdots, m-\ell$.

Now we are ready to provide the definition of the ``curvature'' of
the $\Pi$-connection.

\begin{defn}\label{curvature}
Let $\Pi \in \AA_E(TY)$ and denote by $\Pi: TY = G_\Pi \oplus
T\FF$ the corresponding splitting. The {\it transverse
$\Pi$-curvature} of the foliation $\FF$ is a $T\FF$-valued
two form defined on $N\FF$ as follows: Let $\pi: TY \to N\FF$ be
the canonical projection and
$$
\pi_\Pi: G_\Pi \to N\FF
$$
be the induced isomorphism. Then we define
$$
F_\Pi: \Gamma(N\FF) \otimes \Gamma(N\FF) \to \Gamma(T\FF)
$$
by
\begin{equation}\label{eq:Fdef}
F_\Pi(\eta_1, \eta_2): = \Pi([X,Y])
\end{equation}
where $X = \pi_\Pi^{-1}(\eta_1)$ and $Y= \pi_\Pi^{-1}(\eta_2)$
and  $[X,Y]$ is the Lie bracket on $Y$.
\end{defn}

The following proposition, which is straightforward to check,
shows that $F_\Pi$ is a tensorial object, justifying the name {\it
transverse $\Pi$-curvature}. This  tensor will play a crucial role
in our description of the strong homotopy Lie algebroid associated
to the pre-symplectic manifold $(Y,\omega_Y)$ (and so of
coisotropic submanifolds) and its Maurer-Cartan equation.

\begin{prop}\label{welldefined} Let $F_\Pi$ be as above.
For any smooth functions $f, \, g$ on $Y$
and sections $\eta_1, \, \eta_2$ of $N\FF$, we have the identity
$$
F_\Pi(f\eta_1,g\eta_2) = fgF_\Pi(\eta_1,\eta_2)
$$
i.e., the map $F_\Pi$ defines a well-defined section as an element in
$\Gamma(\Lambda^2(N^*\FF)\otimes T\FF)$.
\end{prop}
\begin{proof}
Let $X, \, Y$ be the unique lifts of $\eta_1, \, \eta_2$ in $\Gamma(G)
\subset \Gamma(TY)$ as in the definition. Then it
follows that $fX, \, gY$ are the lifts of $f\eta_1, \, g \eta_2$.  We have
$$
[fX, gY] = fg[X,Y] + fX[g] Y - gY[f] X.
$$
Since $X, \, Y$ are tangent to $G$, we derive
$$
\Pi ([fX,gY]) = fg \Pi([X,Y])
$$
which finishes the proof.
\end{proof}

In the foliation coordinates $(y^1,\cdots,y^{\ell},q^1, \cdots,
q^{m-\ell})$, $F_\Pi$ has the expression
\begin{equation}
F_\Pi = F^\beta_{ij} \frac{\del}{\del q^\beta} \otimes dy^i \wedge
dy^j \in \Gamma(\Lambda^2(N^*\CF) \otimes T\FF),
\end{equation}
where
\begin{equation}\label{eq:curvature}
F^\beta_{ij}= \frac{\del R^\beta_j}{\del y^i} - \frac{\del
R^\beta_i}{\del y^j} + R^\gamma_i \frac{\del R^\beta_j}{\del
q^\gamma} -R^\gamma_j \frac{\del R^\beta_i}{\del q^\gamma}.
\end{equation}

We next derive the relationship between $F_{\Pi_0}$ and $F_{\Pi}$.
Note that with respect to the given splitting
$$
\Pi_0:\, TY = G_0\oplus T\FF \cong N\FF \oplus T\FF
$$
any other projection $\Pi: TY \to TY$ can be written as the
following block matrix
$$
\Pi = \Big(\begin{matrix}
0 & 0 \\
B & Id
\end{matrix}\Big)
$$
where $B = B_{\Pi_0\Pi}\circ \pi_{G_0}: G_0 \to T\FF$ is the
bundle map introduced in  Lemma \ref{splittings} which is uniquely
determined by $\Pi_0$ and $\Pi$ and vice versa. The following
lemma shows their relationship in coordinates.

\begin{lem}\label{relation}
Let $F_{\Pi}$ and $F_{\Pi_0}$ be the transverse $\Pi$-curvatures
with respect to $\Pi$ and $\Pi_0$ respectively, and let
$B=B_{\Pi_0\Pi}$ be the bundle map mentioned above.
In terms of the foliation coordinates, we have
\begin{eqnarray}
F^\beta_{ij} & = F^\beta_{0,ij} + \Big(\frac{\del B_j^\beta}{\del y^i}
-\frac{\del B^\beta_i}{\del y^j} +
R^\alpha_i\frac{\del B^\beta_j} {\del q^\alpha} -
R^\alpha_j\frac{\del B^\beta_i} {\del q^\alpha} +
B^\alpha_i\frac{\del R^\beta_j} {\del q^\alpha} -
B^\alpha_j\frac{\del R^\beta_i} {\del q^\alpha}
\Big) \nonumber \\
& \quad + \Big(B_i^\alpha \frac{\del B^\beta_j} {\del q^\alpha}
-  B^\alpha_j\frac{\del B^\beta_i} {\del q^\alpha}\Big)
\label{eq:coord-relation}
\end{eqnarray}
\end{lem}
\begin{proof} Let $(y^1, \cdots, y^{\ell},q^1,\cdots, q^{m-\ell})$ be a
foliation coordinates of $\FF$ and let $R_j^\beta$ be the
``Christoffel symbols'' for $\Pi_0$. Let $B_{\Pi_0\Pi}: N\FF \to
T\FF$ be the above bundle map associated to the pair $\Pi_0, \,
\Pi$ of splittings. Using the above block decomposition of $\Pi$,
it follows that the corresponding ``Christoffel symbols'' for
$\Pi$ are given by
$$
R^\beta_j + B^\beta_j
$$
where $(B^\beta_j)$ is the matrix of $B_{\Pi_0\Pi}$ in terms of
the bases of $N\FF$ and $T\FF$ associated to the foliation
coordinates. Then the formula (\ref{eq:coord-relation}) follows
immediately from (\ref{eq:curvature}) by substituting $R^\beta_j$
by $R^\beta_j + B^\beta_j$.
\end{proof}

Now we provide an invariant description of the above formula
(\ref{eq:coord-relation}). Consider the sheaf
$\Lambda^\bullet(N^*\FF)\otimes T\FF$ and denote by
$$
\Omega^\bullet(N^*\FF; T\FF): =
\Gamma(\Lambda^\bullet(N^*\FF)\otimes T\FF)
$$
the group of (local) sections thereof. For an invariant
interpretation of the above basis of $G_x$ and the transformation
law (\ref{eq:coord-relation}), we need to use the notion of {\it
basic vector fields} (or {\it projectable vector fields})
which is standard in the foliation
theory (see e.g., [MM]) :   Consider the Lie subalgebra
$$
L(Y,\FF) = \{ \xi \in \Gamma(TY) \mid ad_\xi(\Gamma(T\FF)) \subset
\Gamma(T\FF) \}
$$
and its quotient Lie algebra
$$
\ell(Y,\FF) = L(Y,\FF) / \Gamma(T\FF).
$$
An element from $\ell(Y,\FF)$ is called a {\it transverse vector
field} of $\FF$. In general,
there may not be a global basic lifting $Y$ of a given transverse vector
field. But the following lemma shows that  this is always possible
locally.

\begin{lem}\label{PiYx=v} Let $x_0 \in Y$ and
$v \in N_{x_0}\FF$. Then there exists a local basic vector
field $\xi$ in a neighborhood of $x_0$ such that it is tangent to $G$
$$
\pi(\xi(x_0)) = v
$$
where $\pi: TY \to N\FF$ is the canonical projection.
\end{lem}
\begin{proof} Using the foliation coordinates, it is well-known and easy to
check that $N_{x_0}\FF$ is spanned by (local) basic vector
fields. Let $\xi'$ be any such transverse vector field with
$\pi(Y'(x_0) )=v$. Then we just take $\xi = \xi' - \Pi(\xi')$ which is
obviously transverse because $Y'$ is transverse and $\Pi(\xi')$ is
tangent to $\FF$.
\end{proof}

\begin{defn} Let $\FF$ be a foliation on $Y$. Let $\Pi \in
\AA_E(TY)$ and $\Pi : TY = G_\Pi \oplus T\FF$ be the
$\Pi$-splitting. We call a basic vector field $\xi$ tangent to $G_\Pi$
a {\it $\Pi$-basic vector field} or a {\it $G$-basic vector field}.
\end{defn}

In this point of view, the vector field
$$
Y_i: = \frac{\del}{\del y^i} + \sum_{\alpha
=1}^{n-k} R_i^\alpha \frac{\del}{\del q^\alpha}
$$
is the unique $G$-basic vector field that satisfies
$$
Y_j \equiv \frac{\del}{\del y^i} \mod T\FF,
$$
i.e., defines the same transverse vector field as $\frac{\del}{\del y^i}$.
\begin{defn}\label{PiLie}
Let $X$ be any (local) basic vector field of $\FF$ tangent to
$G_\Pi$. We define the {\it $\Pi$-Lie derivative} of $B$
with respect to $X$ by the formula
\begin{equation}\label{eq:PiLie}
L_X^\Pi B = \sum_{i_1 < \cdots < i_\ell} L_X(
B_{i_1i_2\cdots i_\ell}) dy^{i_1}\wedge \cdots \wedge
dy^{i_\ell}
\end{equation}
where $B_{i_1i_2\cdots i_\ell}$ is a local section of $T\FF$
given by the local representation of $B$
$$
B = \sum_{i_1< \cdots < i_\ell} B_{i_1 \cdots
i_\ell} dy^{i_1}\wedge \cdots \wedge dy^{i_\ell}
$$
in any given foliation coordinates. Here  $B_{i_1 \cdots
i_\ell}$ is the (locally defined) leafwise tangent
vector field given by
$$
B_{i_1 \cdots
i_\ell} = B_{i_1 \cdots i_\ell}^\beta \frac{\del}{\del q^\beta}.
$$
\end{defn}
From now on without mentioning further, we will always
assume that $B$ is locally defined, unless otherwise stated.

\begin{defn}\label{dPi}
For any element $B \in \Gamma(\Lambda^\ell(N^*\FF);TF)$, we define
$$
d^\Pi B \in \Gamma(\Lambda^{\ell+1}(N^*\FF);TF)
$$
by the formula
\begin{equation}\label{eq:dPi}
d^\Pi B = \sum_{j = 1}^{2k} dy^j \wedge L^\Pi_{Y_j} B
\end{equation}
where we call the operator $d^\Pi$ the {\it $\Pi$-differential}.
\end{defn}

For given splitting $\Pi$ and a vector field $\xi$, we denote by
$\xi^\Pi$ the projection of $\xi$ to $G=G_\Pi$, i.e.,
$$
\xi^\Pi = \xi - \Pi(\xi).
$$
Then the definition of $d^\Pi$ can be also given  by the same kind of
formula as that of the usual exterior derivative $d$: For given $B
\in \Omega^k(N^*\FF;T\FF)$ and local sections $\eta_1, \cdots,
\eta_{k+1} \in N_x\FF$, we define
\begin{eqnarray}
d^\Pi B(v_1, & \cdots& , v_k, v_{k+1}) \nonumber \\
& = & \sum_i (-1)^{i-1}X_i(B(\eta_1, \cdots, \widehat{\eta_i},
\cdots, \eta_{k+1}))
\nonumber\\
\, & + & \sum_{i < j}(-1)^{i+j -1}B(\pi([X_i,X_j]), \eta_1, \cdots,
\widehat{\eta_i}, \cdots, \widehat{\eta_j}, \cdots, \eta_{k+1}) :
\label{eq:inv-dPi}
\end{eqnarray}
Here $X_i$ is a $\Pi$-basic vector field with $\pi(X_i(x)) =
\eta_i(x)$ for each given point $x \in Y$.

It is straightforward
to check that this definition coincides with (\ref{eq:dPi}).

Next we introduce the analog of the ``bracket''
$$
[ \cdot, \cdot]_\Pi  : \Omega^{\ell_1}(N^*\FF; T\FF)
\otimes \Omega^{\ell_2}(N^*\FF;
T\FF) \to \Omega^{\ell_1 + \ell_2}(N^*\FF;T\FF).
$$
\begin{defn} Let $B \in \Omega^{\ell_1}(N^*\FF; T\FF), \, C \in
\Omega^{\ell_2}(N^*\FF; T\FF)$. We define their bracket
$$
[B,C]_\Pi \in \Omega^{\ell_1 + \ell_2}(N\FF;T\FF)
$$
by the formula
\begin{eqnarray}
& & [B,C]_\Pi(v_1,\cdots, v_{\ell_1},
v_{\ell_1+1},\cdots,v_{\ell_1
+ \ell_2 })  \nonumber \\
& = & \sum_{\sigma \in S_n}
\frac{\mbox{sign}(\sigma)}{(\ell_1 +\ell_2)!}
 [B(X_{\sigma(1)}, \cdots, X_{\sigma(\ell_1)}),
C(X_{\sigma(\ell_1 +1)}, \cdots, X_{\sigma(\ell_1 + \ell_2)}]
\label{eq:wedge1}\\
& = & \sum_{\tau \in Shuff(n)} \frac{\mbox{sign}(\tau)}{\ell_1!
\ell_2!} [B(X_{\tau(1)}, \cdots, X_{\tau(\ell_1)}),\nonumber \\
&\quad & \hskip2.0in C(X_{\tau(\ell_1 +1)}, \cdots, X_{\tau(\ell_1
+ \ell_2)}]\label{eq:wedge2}
\end{eqnarray}
for each $x\in Y$ and $v_i \in N_x\FF$, and $X_i$'s are (local)
$\Pi$-basic vector fields such that $\pi(X_i(x)) =
v_i$ as before. Here $S_n$ is the symmetric group with size $n$
and $Shuff(n) \subset S_n$ is the subgroup of all ``shuffles''. $[\cdot, \cdot]$
is the usual Lie bracket of leafwise vector fields.
\end{defn}

For the case $\ell_1 = \ell_2 = 1$, we derive the coordinate
formula
\begin{equation}\label{eq:[B,C]incoord}
[B, C]_\Pi = \Big(B_i^\alpha \frac{\del C^\beta_j} {\del q^\alpha} -
C^\alpha_j\frac{\del B^\beta_i} {\del q^\alpha}\Big) \frac{\del}{\del q^\beta}
\otimes  dy^i \wedge dy^j .
\end{equation}

With these definitions, we have the following ``Bianchi identity''
in our context.
\begin{prop}
Let $\Pi: TY = G \oplus T\FF$ and $d^\Pi$ be the associated
$\Pi$-differential.  Then we have
\begin{eqnarray*}
d^\Pi F_\Pi & = & 0 \label{eq:Bianchi}\\
(d^\Pi)^2 B & = & [F_\Pi, B]_\Pi.
\label{eq:dPi2} \\
\end{eqnarray*}
\end{prop}
\begin{proof} As before, we denote
$Y_i = \frac{\del}{\del y^i} + R^\alpha_i\frac{\del}{\del
q^\alpha}$.  We compute
$$
d^\Pi F_\Pi  =  \frac{1}{2!} dy^i \wedge L_{Y_i}(F_{jk})dy^j\wedge
dy^k, \quad F_{jk} = F^\alpha_{jk}\frac{\del}{\del q^\alpha}.
$$
By definition of $F_\Pi$, we have
$$
L_{Y_i}(F_{jk}) = L_{Y_i}(\Pi[Y_j,Y_k])
$$
And we have
$$
L_{Y_i}(\Pi[Y_j,Y_k]) = L_{Y_i}\Pi([Y_j,Y_k]) +
\Pi([Y_i,[Y_j,Y_k]]).
$$
Here after taking the cyclic sum over $i, \, j, \, k$, the second term
vanishes by the Jacobi identity of the Lie bracket. On the other
hand, differentiating $\Pi^2 = \Pi$, we derive
$$
(L_X\Pi)([Y,Z])  = \Pi (L_X\Pi)([Y,Z])+(L_X\Pi\cdot \Pi)([Y,Z])
$$
in general. Furthermore, {\it if we restrict to the
basic vector fields} $X, \, Y, \, Z$ and so
$[Y,Z]$ is also basic, then it immediately follows
from $\Pi^2 = \Pi$ that the
second term vanishes. The first term becomes
\begin{eqnarray*}
\Pi (L_X\Pi)([Y,Z]) & = & \Pi (L_X(\Pi[Y,Z] ) - \Pi \cdot \Pi
L_X([Y,Z])
\\ & = & \Pi ([X, \Pi[Y,Z]] - [X, [Y,Z]]) = -\Pi ([X, [Y,Z] -
\Pi[Y,Z]).
\end{eqnarray*}
{\it When this identity is applied to}
$$
X=Y_i, \, Y= Y_j, \, Z= Y_k,
$$
this term also vanishes because we have
$$
[Y_j, Y_k] = \Big[\frac{\del}{\del y^j} +
R^\alpha_j\frac{\del}{\del q^\alpha},\frac{\del}{\del y^k} +
R^\alpha_k\frac{\del}{\del q^\alpha}\Big] = F^\alpha_{jk}\frac{\del}
{\del q^\alpha}
$$
which is tangent to $\FF$. This finishes proof of the Bianchi
identity, $d^\Pi F_\Pi = 0$.

For the proof of (\ref{eq:dPi2}), we consider two $\Pi$-basic vector
fields
$$
Y_i = \frac{\del}{\del y^j}+ R^\alpha_i\frac{\del}{\del q^\alpha},
\, Y_j = \frac{\del}{\del y^j}+ R^\alpha_j\frac{\del}{\del
q^\alpha}.
$$
We recall
$$
d^\Pi B = dy^j \wedge L_{Y_j}B = \frac{1}{\ell !}L_{Y_j}(
B_{i_1\cdots i_\ell}) dy^j \wedge dy^{i_1}\wedge\cdots \wedge
dy^{i_\ell}.
$$
Therefore we have
\begin{eqnarray*}
(d^\Pi)^2 B & = &\frac{1}{\ell !}
L_{Y_i}L_{Y_j}(B_{i_1\cdots i_\ell})dy^i \wedge dy^j \wedge
dy^{i_1}\wedge\cdots \wedge dy^{i_\ell} \\
& = & \sum_{i< j}\frac{1}{\ell !}(L_{Y_i}L_{Y_j} -
L_{Y_j}L_{Y_i})(B_{i_1\cdots i_\ell})dy^i \wedge dy^j
\wedge dy^{i_1}\wedge\cdots \wedge
dy^{i_\ell}\\
& = &\sum_{i< j}\frac{1}{\ell !} L_{[Y_i,Y_j]}(B_{i_1\cdots
i_\ell})dy^i \wedge dy^j \wedge dy^{i_1}\wedge\cdots \wedge
dy^{i_\ell}.
\end{eqnarray*}
Here we note that  $[Y_i,Y_j]$ is tangent to $\FF$ and hence
$$
[Y_i,Y_j](x) = \Pi[Y_i,Y_j](x) = F_\Pi(u_i, u_j)
$$
where $u_i = \pi(Y_i)(x)$. Therefore we have
$$
(d^\Pi)^2 B = [F_\Pi , B]_\Pi
$$
on $\Omega^{\ell_1 +\ell_2}(N^*\FF;T\FF)$
which finishes the proof of (\ref{eq:dPi2}).
\end{proof}

Combining the above discussion, the transformation law
(\ref{eq:coord-relation}) in coordinates is translated into the
following invariant form.

\begin{prop}\label{inv-relation}
Let $\Pi, \, \Pi_0$ be two splittings as in Lemma \ref{relation}
and $B_{\Pi_0\Pi} \in \Gamma(N^*\FF\otimes T\FF)$ be the
associated section. Then we have
\begin{equation}\label{eq:inv-relation}
F_\Pi = F_{\Pi_0} + d^{\Pi_0}B_{\Pi_0\Pi} + [B_{\Pi_0\Pi},
B_{\Pi_0\Pi}]_{\Pi_0}.
\end{equation}
\end{prop}

\begin{rem}
We would like to emphasize that the bracket $[\cdot, \cdot]_\Pi$
we defined is not bilinear over $C^\infty(Y)$, but linear over the
so-called  subalgebra of {\it basic functions}:
a smooth function $f: Y \to \R$ is called basic, if  it is constant
along the leaves(see [MM], [To] for the definition). It seems that a
good formulation of invariant objects ``over the leaf space''
should be  in terms of  the Haefliger-type cocycles and germs of
$\Omega^\bullet(T;T\FF)$ over transverse sections $T$ of the
foliation $\FF$ as in [Ha] for his definition of
$\Omega^\bullet(Tr\FF)$. We postpone elsewhere a full disclosure
of geometric structures that arise in the study of the deformation
problem of foliations $\FF$ and the role of the transverse
curvature $F_\Pi$.
\end{rem}

\section{Geometry of the null foliation}
\label{sec:nullfoliation}

In this section, we will apply the leaf space calculus developed
in the previous section to the null foliation. Let $(Y,\omega)$ be
a pre-symplectic manifold and denote by $\FF$ the associated null
foliation. About the range of indices of $i$ and $\alpha$ from section
 \ref{sec:splittings}, we use
$$
m = n+k, \, \ell = 2k, \quad 1 \leq i \leq 2k, \, 1 \leq \alpha \leq n-k
$$
for the null foliation for $(Y,\omega)$.

Using Lemma \ref{splittings}, we now state a uniqueness statement
in the symplectic thickening of $(Y,\omega)$. This precise form of
the neighborhood theorem will be a crucial ingredient for our
proof in section \ref{sec:gauge} of the gauge equivalence of
strong homotopy Lie algebroids that we associate to the splittings
$\Pi$. We denote by $\omega_\Pi$ the symplectic form given in
(\ref{eq:omega*}) associated to the splitting $\Pi$.

\begin{prop}\label{diffeo} For given two splittings $\Pi_0, \, \Pi$, there exist
neighborhoods $U, \, U'$ of the zero section $Y \subset E^*$ and a
diffeomorphism $\phi: U \to U'$ such that
\begin{enumerate}
\item $\phi^*\omega_\Pi = \omega_{\Pi_0}$, \item $\phi|_Y \equiv
id$, and $T\phi|_{T_YE^*} \equiv id$ where $T_YE^*$ is the
restriction of $TE^*$ to $Y$.
\end{enumerate}
\end{prop}
\begin{proof}
Since $\AA_E(TY)$ is contractible, we can choose a smooth family
$$
\{\Pi_t\}_{0 \leq t\leq 1}, \quad \Pi_0=\Pi_0, \, \Pi_1 = \Pi.
$$
Denoting $\omega_t :=\omega_{\Pi_t}$, we have
$$
\omega_t - \omega_0 = d(\theta_0 - \theta_t)
$$
$\theta_t$ is the one-form $\theta_G$ associated to $\Pi_t$ as
defined in (\ref{eq:thetaG}). From the definition, it follows
$\theta_t|_{T_YE^*} \equiv 0$ and hence
$$
(\theta_0- \theta_t)|_{T_YE^*} \equiv 0.
$$
for all $0 \leq t \leq 1$. With these, the proof follows from the
standard Moser's homotopy method [We1], since $\omega_t$ are all
nondegenerate and homologous to each other in a sufficiently small
neighborhood of the zero section $Y \subset E^*$.
\end{proof}

For the study of the deformation problem of pre-symplectic
structures it is crucial to understand the transverse geometry of
the null foliation. First we note that the pre-symplectic form
$\omega$ carries a natural {\it transverse symplectic form}. This
defines the symplectic analog to the much-studied Riemannian
foliation (see [Mo], [To] for example).

\begin{prop}\label{tr-sy-form}  Let $\FF$ be the null foliation of
the pre-symplectic manifold $(Y,\omega)$. Then the pre-symplectic form
$\omega$ defines a transverse symplectic form on $\FF$ in the following sense:
\enumerate
\item $\ker (\omega_x) = T_x\FF$ for any $x \in M$, and
\item $L_X\omega = 0$ for any vector field on $M$ tangent to $\FF$
\endenumerate
\end{prop}
\begin{proof} The first statement is trivial by definition of the
null foliation and the second is an immediate consequence of the
Cartan identity
$$
L_X\omega = d(X \rfloor \omega) + X \rfloor d\omega.
$$
The second term on the right hand side vanishes since $\omega$
is closed, and the first also
vanishes if $X$ is tangent to the null foliation $\FF$.
\end{proof}

One immediate consequence of the presence of the transverse symplectic
form, together with the fact that the pre-symplectic form $\omega$
is closed, is that any {\it transverse section}  $T$ of the
foliation $\FF$ carries a natural symplectic form: in any
foliation coordinates, it follows from $E = \ker \omega =
\operatorname{span}\{\frac{\del}{\del q^\alpha}\}_{1 \leq \alpha
\leq n-k}$ that we have
\begin{equation}\label{eq:piomegaY}
\pi^*\omega = \frac{1}{2}\omega_{ij} dy^i \wedge dy^j,
\end{equation}
where $\omega_{ij} = \omega(\frac{\del}{\del y^i},\frac{\del}{\del
y^j})$ is skew-symmetric and invertible. And closedness of
$\omega$ implies that $\omega_{ij}$ is independent of
$q^\alpha$'s. Note that this expression is {\it independent} of
the choice of splitting as long as $y^1, \cdots, y^{2k}$ are those
coordinates that characterize the leaves of $E$ by (\ref{eq:leaves}).

The proof of the following proposition is straightforward by
definition of the holonomy map and is omitted (see [Proposition 2.5,
MM] for a proof of its Riemannian analog).

\begin{prop}\label{holonomy}
Let $L$ be a leaf of the null foliation $\FF$ on $(Y,\omega)$,
$\lambda$ a path in $L$, and let $T$ and $S$  be transverse
sections of $\FF$ with $\lambda(0) \in T$ and $\lambda(1) \in S$.
Then the holonomy map
$$
hol^{S,T}(\lambda): (T,\lambda(0)) \to (S,\lambda(1))
$$
defines the germ of a symplectic diffeomorphism.
\end{prop}

Now for the case of the null foliation $\FF$ which carries a
transverse symplectic structure, we can go one step further with
the curvature $F_\Pi$ from the case of general foliations: we can
do the operation of ``raising indices'' using symplectic form like
a metric. We recall that $\omega$ defines a non-degenerate
bilinear form on $N\FF$. We denote by $\omega^{-1}$ the natural
bilinear form induced on $N^*\FF$. We refer to section
\ref{sec:algebroid} for more on this.

Now we introduce the following notion of {\it
symplectic mean transverse $\Pi$-curvature} of the null foliation
$\FF$.
This generalizes the {\it Reeb vector field} of a contact form
on the contact manifold $(Y,\xi)$ to arbitrary pre-symplectic
manifolds.

\begin{defn}\label{rhoFF}
Define the section $\rho_\Pi \in \Gamma(T\FF)$ by
\begin{equation}\label{eq:rhoFF}
\rho_\Pi = \frac{1}{2k}\operatorname{trace}_\omega F_\Pi
: = \frac{1}{2k}\langle F_\Pi,\omega^{-1} \rangle
\end{equation}
where $\omega^{-1}$ is the inverse of $\omega$ on $N\FF$
and $\langle\cdot, \cdot \rangle$ is the natural pairing
between $\Lambda^2G^*$ and $\Lambda^2(N\FF)$.
In a foliation coordinates, it is given by
\begin{equation}\label{eq:coordrhoFF}
\rho_\Pi = \frac{1}{2k} F^\beta_{ij}\omega^{ij}\frac{\del}{\del q^\beta}
\end{equation}
where $(\omega^{ij})$ is the inverse of $(\omega_{ij})$ with
$\omega = \frac{1}{2}\omega_{ij}dy^i\wedge dy^j$ (see
 section \ref{sec:connections}).
\end{defn}

\begin{thm}\label{contact}
Let $(Y,\xi)$ be a contact manifold of dimension $2n-1$
and with $\xi$ its contact
distribution. Choose a contact one-form
$\theta$ and consider the pre-symplectic form
$\omega = -d\theta$.
Denote by $X_\theta$ the associated Reeb vector field.
We set $E = \ker d\theta$ and $\FF_\theta$ be the associated
foliation (or the line field).
Consider the projection $\Pi: TY \to TY$ such that
$$
G_\Pi = \xi.
$$
Then for any $X,\, Y \in \xi \cong N\FF_\theta$, we have
\begin{equation}F_\Pi(X,Y)  =  \theta([X,Y])X_\theta
= -d\theta(X, Y)\cdot X_\theta \label{eq:contactF}
\end{equation}
and
\begin{equation}
\rho_\Pi = X_\theta. \label{eq:contactrho}
\end{equation}
\end{thm}
\begin{proof}
Recall that the contact vector field $X_\theta$ is the unique
vector field
that satisfies
$$
X_\theta \rfloor \theta = 1, \quad X_\theta \rfloor d\theta =0.
$$
Therefore it immediately implies
$$
F_\Pi(X,Y) = \Pi([X,Y]) = \theta([X,Y])X_\theta
$$
for any $X, \, Y \in \xi$
since we have $\theta([X,Y] - \theta([X,Y])X_\theta) = 0$, i.e,
$$
[X,Y] - \theta([X,Y])X_\theta \in \xi = G_\Pi.
$$
The second identity of (\ref{eq:contactF}) follows from
the definition of the exterior derivative $d\theta$ and the
defining equation of contact form
$\xi = \ker \theta$.

The identity (\ref{eq:contactrho}) immediately follows from
the choice of pre-symplectic form $\omega = -d\theta$ and
(\ref{eq:contactF}).
\end{proof}

The following is an interesting consequence of this theorem.

\begin{cor}\label{hypercont}
Let $(X,\omega_X)$ be a given symplectic manifold
and $J$ be a compatible almost complex structure.
Denote by $\HH yper(X,\omega_X, J)$ the space of hypersurfaces
with the induced pre-symplectic form $\omega$ and with the
Riemannian metric
induced from $g_J=\omega_X(\cdot, J\cdot)$.
We choose the orthogonal splitting $\Pi: TY = G \oplus T\FF$
for each hypersurface $Y \in \HH yper(X,\omega_X, J)$
and denote by $\rho_{(J,Y)}$ the corresponding symplectic
mean transverse curvature.
Let $\CC ont(X,\omega_X, J)$ the subset of contact hypersurface
with the contact form $\theta_{(J,Y)} : = \omega_X(N_J, \cdot)$
where $N_J$ is the unit (positive) normal vector field
and denote by $X_{(J,Y)}$ the Reeb vector field of $\theta_{(J,Y)}$.
Then the assignment
$$
(J,Y) \in \JJ(X,\omega) \times \HH yper(X,\omega_X,J)
\mapsto \rho_{(J,Y)}
$$
is continuous with respect to the $C^\infty$-topology
which extends the assignment of the Reeb vector fields
$$
(J,Y) \in \JJ(X,\omega) \times \CC ont(X,\omega_X, J)
\mapsto X_{(J,Y)}.
$$
\end{cor}

We believe that this corollary will play some role in the study of
Hamiltonian dynamics on the hypersurfaces. This will be a subject
of the future study. We also refer to [Oh] for a study of
coisotropic submanifolds in K\"ahler manifolds with respect to the
canonical orthogonal splitting $TY = N_J\FF \oplus T\FF$ where
$N_J\FF$ is the normal bundle of $\FF$ with respect to the
K\"ahler metric $g_J = \omega(\cdot, J \cdot)$.

Another important geometric structure related to the bundle
$$
T\FF= E = \ker \omega \to Y
$$
is the structure of a {\it Lie algebroid} (see e.g. [CW], [MM], [NT]
for its definition) : the {\it anchor map} $E \to TY$ is just the
inclusion map and the Lie bracket on the sections of $E$ is
just the usual bracket of vector fields on $Y$. This induces the
leafwise differential $d_\FF$ and defines the leafwise de Rham
complex $\Omega^*(\FF)$, which will play an important role in the
deformation problem of pre-symplectic structures and coisotropic
submanifolds later.

\section{Geometry of the symplectic thickening}
\label{sec:connections}

In this section and the next two, we will unravel the geometric
and algebraic structure that governs the deformation problem of
coisotropic submanifolds up to Hamiltonian isotopy. As we
showed in section 2, it is enough to study this in the model space
$(U,\omega_U)$ constructed in section 2. Precise formulation of
the problem is in order.

Again we start with a splitting
$$
TY = G \oplus E,
$$
the associated bundle projection $\Pi: TY \to TY$,  the associated
canonical one form $\theta_G$, and the symplectic form
$$
\omega_{E^*} = \pi^*\omega - d\theta_G
$$
on $U \subset E^*$. Note the projection map $\pi: U \to Y$ induces
a foliation $\pi^{-1}(\FF)$ on $U$ in a canonical way. The leaves
of $\pi^{-1}(\FF)$ are the preimages of the leaves of $\FF$, which
are symplectic submanifolds of $U$. When we choose  foliation
coordinates $(y^1,\cdots,y^{2k}, q^1, \cdots, q^{n-k})$ on $Y$, we
can extend these coordinates to  foliation coordinates of
$\pi^{-1}(\FF)$
$$
(y^1, \cdots, y^{2k},q^1, \cdots, q^{n-k}, p_1, \cdots, p_{n-k})
$$
so that
$$
y^1 = c^1, \cdots, y^{2k} = c^{2k}, \quad \mbox{$c^i$'s constant}
$$
defines the leaves of the foliation $\pi^{-1}(\FF)$. We will construct these
coordinates explicitly below. We have the following commutative
diagram of exact sequences
\begin{equation}\label{eq:CDES}
\begin{matrix}
0  \quad \to  &  T\pi^{-1}(\FF) \quad \to
& TU  \to  & N(\pi^{-1}(\FF)) \quad  \to  & 0 \\
\qquad  & \downarrow T\pi  &  \downarrow T\pi &\downarrow \cong  \,&\, \\
0  \quad \to &  T\FF \quad   \to  & TY \quad  \to &   N\FF \quad \to & 0
\end{matrix}
\end{equation}

Note that for a given splitting $\Pi: TY = G\oplus T\FF$, due to
the presence of the symplectic form $\omega_U$, there exists the
unique splitting of $TU$
\begin{equation}\label{eq:TUtildeG}
TU = G^\sharp \oplus T\pi^{-1}(\FF)
\end{equation}
that satisfies
\begin{equation}\label{eq:tildeG}
G^\sharp = (T_\alpha\pi^{-1}(\FF))^{\omega_U}
\end{equation}
for any $\alpha \in U$, which is invariant under the action of
symplectic diffeomorphisms on $(U, \omega_U)$ that preserve the
leaves of $\pi^{-1}(\FF)$.

\begin{defn}\label{transymp}
We call the above unique splitting the {\it leafwise symplectic
connection of $U \to Y$ compatible to the splitting $\Pi: TY = G
\oplus T\FF$} or simply a {\it leafwise symplectic
$\Pi$-connection} of $U \to Y$.
\end{defn}
We would like to emphasize that this connection is not
a vector bundle connection of $E^*$ although $U$ is a subset of
$E^*$, which reflects {\it nonlinearity} of this connection.

Note that the splitting $\Pi$ naturally induces the splitting
$$
\Pi_*: T^*Y = (T\FF)^\circ \oplus G^\circ
$$
where $(T\FF)^\circ$ and $G^\circ$ are the annihilators of $T\FF$
and $G$ respectively. We denote by $\Pi^\sharp: TU \to TU$ the
projection to $T(\pi^{-1}\FF)$ associated to (\ref{eq:TUtildeG}).
We now derive the coordinate expression of the leafwise symplectic
$\Pi$-connection. Let
$$
(y^1, \cdots, y^{2k}, y^{2k+1}, \cdots, y^{n+k})
$$
be coordinates on $Y$ adapted to the null foliation on an open
subset $V \subset Y$ as before. By choosing the frame
$$
\{f_1^*, \cdots, f_{n-k}^*\}
$$
of $E^*$ that is dual to the frame $\{\frac{\del}{\del q^1}, \cdots,
\frac{\del}{\del q^{n-k}}\}$ of $E$,
we introduce the {\it canonical coordinates} on $E^*$ by writing
an element $\alpha \in E^*$ as a linear combination of $\{f_1^*,
\cdots, f_{n-k}^*\}$
$$
\alpha = p_\beta f_\beta^*,
$$
and taking
$$
(y^1, \cdots, y^{2k}, q^1, \cdots, q^{n-k}, p_1, \cdots, p_{n-k})
$$
as the associated coordinates. To derive the coordinate expression
of $\theta_G$, we compute
\begin{eqnarray*}
\theta_G\Big(\frac{\del}{\del y^i}\Big) & = & \alpha\Big(p_G\circ
T\pi(\frac{\del}{\del y^i})\Big)
= \alpha\Big(p_G(\frac{\del}{\del y^i})\Big) \\
& = & p_\beta f_\beta^*\Big(-R_i^\alpha \frac{\del}{\del q^\alpha}\Big) = -
p_\alpha
R_i^\alpha, \\
\theta_G\Big(\frac{\del}{\del q^\beta}\Big)  & = & p_\beta, \qquad
\theta_G\Big(\frac{\del}{\del p_\beta}\Big) = 0.
\end{eqnarray*}
Hence we derive
\begin{equation}\label{eq:thetaG}
\theta_G = p_\beta(dq^\beta - R_i^\beta dy^i).
\end{equation}
Here we note that
$$
(dq^\beta - R^\beta_idy^i)|_{G_x} \equiv 0.
$$
This shows that if we identify $E^* = T^*\FF$ with
$$
G^\circ:= \operatorname{im}(\Pi)^*|_{E^*} \subset T^*Y
$$
via the embedding $(\Pi)^*|_{E^*}: E^* \hookrightarrow T^*Y$
induced by the splitting $\Pi: TY = G\oplus T\FF$, then we may
write the dual frame on $T^*\FF$ as
\begin{equation}\label{eq:fbeta*}
f_\beta^* = dq^\beta - R^\beta_idy^i.
\end{equation}
Motivated by this, we write
\begin{eqnarray}
d\theta_G & = & dp_\beta \wedge (dq^\beta - R_i^\beta dy^i)
- p_\beta dR_i^\beta \wedge dy^i
\label{eq:dthetaG} \\
& = & dp_\beta \wedge (dq^\beta - R_i^\beta dy^i)
 - p_\beta\frac{\del R_i^\beta}{\del y^j} dy^j \wedge dy^i \nonumber\\
&\quad &  - p_\beta\frac{\del R_i^\beta}{\del q^\gamma} dq^\gamma
\wedge dy^i.
\end{eqnarray}

Combining (\ref{eq:dthetaG}) and (\ref{eq:piomegaY}), we have
\begin{eqnarray}
\omega_U & = & \pi^*\omega - d\theta_G 
\nonumber\\
& = & \Big(\frac{1}{2}\omega_{ij} - p_\beta \frac{\del
R_j^\beta}{\del y^i}\Big) dy^i \wedge dy^j  -dp_\delta \wedge
(dq^\delta - R_i^\delta dy^i) \nonumber \\
& \quad & + p_\nu \frac{\del R_i^\nu}{\del q^\gamma} dq^\gamma \wedge dy^i
\nonumber\\
& = & \Big(\frac{1}{2}\omega_{ij} - p_\beta \frac{\del
R_j^\beta}{\del y^i} - p_\beta \frac{\del R_i^\beta}{\del
q^\gamma} R_j^\gamma \Big) dy^i
\wedge dy^j \nonumber \\
&\quad &  - (dp_\delta + p_\beta \frac{\del R_i^\beta}{\del
q^\delta}
dy^i) \wedge (dq^\delta - R_j^\delta dy^j) \nonumber\\
& = & \Big(\frac{1}{2}\omega_{ij} - p_\beta (\frac{\del
R_j^\beta}{\del y^i} - \frac{\del R_i^\beta}{\del q^\gamma}
R_j^\gamma) \Big) dy^i
\wedge dy^j \nonumber \\
& \quad &  - (dp_\beta + p_\beta \frac{\del R_i^\beta}{\del
q^\gamma} dy^i) \wedge (dq^\delta - R_j^\delta dy^j)\nonumber\\
& = & \frac{1}{2}\Big(\omega_{ij} - p_\beta F^\beta_{ij} \Big)
dy^i \wedge dy^j \nonumber \\
& \quad &  - (dp_\delta+ p_\beta \frac{\del R_i^\beta}{\del
q^\delta} dy^i) \wedge (dq^\delta- R_j^\delta dy^j)
\label{eq:omegaU}
\end{eqnarray}
where $F^\beta_{ij}$ are the components of the transverse
$\Pi$-curvature of the null-foliation. Here the last identity
comes from the skew-symmetry of $\omega_{ij}$ and by
anti-symmetrizing the first term of the identity right before.

Note that we have
$$
T\pi^{-1}(\FF) = \operatorname{span} \Big\{\frac{\del}{\del q^1},
\cdots, \frac{\del}{\del q^{n-k}},  \frac{\del}{\del p_1}, \cdots,
\frac{\del}{\del p_{n-k}} \Big\}
$$
which is independent of the choice of the above induced foliation
coordinates of $TU$.

Now we compute $(T\pi^{-1}(\FF))^{\omega_U}$ in $TU$ in terms of
these induced foliation coordinates. We will determine when the
expression
$$
a^j(\frac{\del}{\del y^j} + R_j^\alpha\frac{\del}{\del q^\alpha})
+ b^\beta \frac{\del}{\del q^\beta} + c_\gamma \frac{\del}{\del
p_\gamma}
$$
satisfies
$$
\omega_U\Big( a^j(\frac{\del}{\del y^j} +
R_j^\alpha\frac{\del}{\del q^\alpha}) + b^\beta \frac{\del}{\del
q^\beta}+ c_\gamma \frac{\del}{\del p_\gamma},
T\pi^{-1}(\FF)\Big)= 0.
$$
It is immediate to see by
pairing with $\frac{\del}{\del p_\mu}$
\begin{equation}\label{eq:bbeta}
b^\beta = 0, \quad \beta = 1, \cdots, n-k.
\end{equation}
Next we study the equation
$$
0 = \omega_U\Big (a^j(\frac{\del}{\del y^j} +
R_j^\alpha\frac{\del}{\del q^\alpha}) + c_\gamma \frac{\del}{\del
p_\gamma}, \frac{\del}{\del q^\nu}\Big)
$$
for all $\nu = 1, \cdots, n-k$. A straightforward check provides
\begin{equation}\label{eq:cnu}
- a^jp_\beta\frac{\del R_j^\beta}{\del q^\nu} - c_\nu = 0
\end{equation}
for all $\nu$ and $j$. Combining (\ref{eq:bbeta}) and
(\ref{eq:cnu}), we have obtained
\begin{equation}\label{eq:Tpi-1FFomegaU}
(T\pi^{-1}(\FF))^{\omega_U} = \operatorname{span} \Big\{
\frac{\del}{\del y^j} + R_j^\alpha \frac{\del}{\del q^\alpha} -
p_\beta \frac{\del R_j^\beta}{\del q^\nu}\frac{\del}{\del p_\nu}
\Big\} _{1 \leq j \leq 2k}.
\end{equation}

\begin{rem} Just as we have been considering $\Pi: TY = G \oplus T\FF$ as a
``connection'' over the leaf space, we may consider the splitting
$\Pi^\sharp: TU = G^\sharp \oplus T(\pi^{-1}\FF)$ as the leaf
space connection {\it canonically induced from $\Pi$ under the
fiber-preserving map}
$$
\pi: U \to Y
$$
over the same leaf space $Y/\sim$: Note that the space of leaves
of $\FF$ and $\pi^{-1}\FF$ are canonically homeomorphic.
\end{rem}

\section{Master equation in coordinates}
\label{sec:master}

We will now study the condition that the graph of a section $s: Y
\to E^*\cong NY$ is coisotropic with respect to $\omega_U$. We
call the corresponding equation the {\it classical part} of the
master equation. We study the full (local) moduli problem of
coisotropic submanifolds by analyzing the condition that the graph
of  a section $s: Y \to U$ in the symplectic thickening $U$ is to
be coisotropic with respect to $\omega_U$.

Recall that an Ehresmann connection of $U \to Y$ with a structure
group $H$ is a splitting of the exact sequence
$$
0 \to VTU \longrightarrow TU \stackrel{T\pi} \longrightarrow TY
\to 0
$$
that is invariant under the action of the group $H$. Here $H$ is
not necessarily a finite dimensional Lie group. In other words, an
Ehresmann connection is a choice of decomposition
$$
TU = HTU \oplus VTU
$$
that is invariant under the fiberwise action of $H$. Recalling
that there is a canonical identification $V_\alpha TU \cong
V_\alpha TE^* \cong E^*_{\pi(\alpha)}$, a connection can be
described as a {\it horizontal} lifting $HT_\alpha U$ of $TY$ to
$TU$ at each point $y \in Y$ and $\alpha \in U \subset E^*$ with
$\pi(\alpha) = y$.  We denote by $F^\# \subset HTU$ the horizontal
lifting of a subbundle $F \subset TY$ in general.

Let $(y^1,\cdots, y^{2k},q^1,\cdots, q^{n-k})$ be a foliation
coordinates of $\FF$ on $Y$ and
$$
(y^1,\cdots, y^{2k},q^1,\cdots, q^{n-k}, p_1, \cdots, p_{n-k})
$$
be the induced foliation coordinates of $\pi^{-1}(\FF)$ on $U$.
Then $G^\# = (T\pi^{-1}(\FF))^{\omega_U}$ has the natural basis
given by
\begin{eqnarray*}\label{eq:ej}
e_j & = & \frac{\del}{\del y^j} + R_j^\alpha \frac{\del}{\del
q^\alpha} - p_\beta \frac{\del R_j^\beta}{\del
q^\nu}\frac{\del}{\del p_\nu}
\end{eqnarray*}
which are basic vector fields of $T(\pi^{-1}\FF)$. We also
denote
$$
f_\alpha = \frac{\del}{\del q^\alpha}
$$
We define a local lifting of $E$
\begin{equation}\label{eq:Esharp}
E^\sharp = \operatorname{span} \Big\{f_1, \cdots, f_{n-k}\Big\}.
\end{equation}

\begin{rem}\label{Esharp}
\begin{enumerate}
\item When a splitting $\Pi: TY = G\oplus T\FF$ is given,
$G^\sharp$ is canonically defined which is independent of the
foliation coordinates.  Unlike $G^\sharp$, $E^\sharp$ depends on
the choice of foliation coordinates. However for another choice of
foliation coordinates $(\underline{y}^1, \cdots,
\underline{y}^{2k}, \underline{q}^1, \cdots,
\underline{q}^{n-k})$, $E^\sharp$ will have a basis in the form
$$
\frac{\del}{\del \underline{q}^\alpha} + A^\alpha_\beta
\frac{\del}{\del \underline{p}_\beta}
$$
for $\alpha = 1, \cdots, n-k$, where the matrix $(A^\alpha_\beta)$
is symmetric. The latter follows from the fact that $E^\sharp$ is
isotropic.
\item The choice of splitting
$$
T(\pi^{-1}\FF) = E^\sharp \oplus VTU
$$
may be considered as the analog of {\it $\FF$-partial connection}
associated to the principal $G$-bundle $P \to B$ for a Lie
groupoid whose associated Lie algebroid is the one associated to
the null foliation $\FF$, {\it except} that in our case we do not
have associated the global Lie groupoid in sight. We refer to
[section 5, MM] for the precise definition of principal
$G$-bundles of a Lie groupoid $G$ and the $\FF$-partial
connection. The relevant groupoid-like object in our case is the
symplectic thickening $\pi: U \to Y$ which can be interpreted as
the ``integrated'' object of our strong homotopy Lie algebroid
that we introduce in the next section. We hope to have more
detailed discussion elsewhere.
\end{enumerate}
\end{rem}
\medskip

(\ref{eq:Esharp}) provides a local splitting
$$
TU = (G^\sharp \oplus E^\sharp) \oplus VTU \to U
$$
and defines a locally defined Ehresmann connection
where $VTU$ is the vertical tangent bundle of $TU$. From the
expression (\ref{eq:omegaU}) of $\omega_U$, it follows that
$G^\sharp \oplus E^\sharp$ is a coisotropic lifting of $TY$ to $T
U$. We denote by $\Pi^v: TU \to VTU$ the vertical projection with
respect to this splitting.

With this preparation, we are finally ready to derive the master
equation. Let $s: Y \to U \subset E^*$ be a section and denote
\begin{equation}\label{eq:nablas}
\nabla s: = \Pi^v\circ ds
\end{equation}
its locally defined ``covariant derivative''. In coordinates
$(y^1, \cdots, y^{2k},q^1, \cdots, q^{n-k})$, we have
\begin{eqnarray*}
ds\Big(\frac{\del}{\del y^j}\Big) & = &\frac{\del}{\del y^j} +
\frac{\del s_\alpha}{\del y^j} \frac{\del}{\del p_\alpha} \\
& = & e_j - R^\alpha_j\frac{\del}{\del q^\alpha} +
s_\beta\frac{\del R_j^\beta}{\del q^\nu}\frac{\del}{\del p_\nu} +
\frac{\del
s_\alpha}{\del y^j}\frac{\del}{\del p_\alpha} \\
& = & e_j - R^\alpha_j f_\alpha + \Big(
s_\beta\frac{\del R_j^\beta}{\del
q^\nu}\frac{\del}{\del p_\nu} + \frac{\del s_\alpha}{\del
y^j}\frac{\del}{\del p_\alpha}\Big).
\end{eqnarray*}
Therefore we have derived
\begin{equation}\label{eq:nablas}
\nabla s\Big(\frac{\del}{\del y^j}\Big) = \Big(\frac{\del s_\nu}{\del
y^j} +  s_\beta\frac{\del
R_j^\beta}{\del q^\nu}\Big)\frac{\del}{\del p_\nu}.
\end{equation}
Similarly we compute
\begin{eqnarray*}
ds\Big(\frac{\del}{\del q^\nu}\Big) & = &\frac{\del}{\del q^\nu} +
\frac{\del s_\alpha}{\del q^\nu}\frac{\del}{\del p_\alpha} \\
& = & \frac{\del}{\del q^\nu} +
\frac{\del s_\alpha}{\del q^\nu} \frac{\del}{\del p_\alpha},
\end{eqnarray*}
and so
\begin{equation}\label{eq:nablanus}
\nabla s\Big(\frac{\del}{\del q^\nu}\Big) =
\frac{\del s_\alpha}{\del q^\nu} \frac{\del}{\del p_\alpha}.
\end{equation}

Recalling that $T_\alpha U = (E_\alpha^\# \oplus VT_\alpha
U)^{\omega_U} \oplus E_\alpha ^\# \oplus VT_\alpha U$, we conclude
that the graph of $ds$ with respect to the frame
$$
\Big\{e_1, \cdots, e_{2k}, f_1, \cdots, f_{n-k}, \frac{\del}{\del
p_1}, \cdots, \frac{\del}{\del p_{n-k}} \Big\}
$$
can be expressed by the linear map
\begin{eqnarray*}
A_H: (E^\#\oplus VTU)^{\omega_U}
& \to & VTU \cong E^*;\quad (A_H)_\alpha^i = \nabla_i s_\alpha, \\
A_I: E^\# & \to &VTU \cong E^*;  \quad (A_I)_\alpha^\beta =
\nabla_\beta s_\alpha,
\end{eqnarray*}
where
\begin{eqnarray}\label{eq:nablaincoord}
\nabla s\Big(\frac{\del}{\del y^i}\Big) & = & (\nabla_i s_\alpha)
\frac{\del}{\del q^\alpha}, \quad  \nabla_i s_\alpha
: =  \frac{\del s_\alpha}{\del y^j}
+ s_\beta\frac{\del R_j^\beta}{\del q^\alpha}, \nonumber\\
\nabla s\Big(\frac{\del}{\del q^\beta}\Big)  & = & (\nabla_\beta s_\alpha)
\frac{\del}{\del q^\alpha}, \quad
\nabla_\beta s_\alpha  : =  \frac{\del s_\alpha}{\del q^\beta}, \nonumber
\end{eqnarray}

Finally we note that
$$
\omega_U(s)(e_i,e_j) = w_{ij} - s_\beta F^\beta_{ij} : =
\widetilde \omega_{ij}
$$
and denote its inverse by $(\widetilde \omega^{ij})$. Note that
$(\widetilde \omega_{ij})$ is invertible if $s_\beta$ is
sufficiently small, i.e., if the section $s$ is $C^0$-close to the
zero section, or its image stays inside of $U$.
Now Proposition 2.2 immediately implies

\begin{thm} Let $\nabla s$ be the vertical projection of $ds$
as in (\ref{eq:nablas}). Then  the graph of the section $s: Y \to
U$ is coisotropic with respect to $\omega_U$ if and only if $s$
satisfies
\begin{equation}
\nabla_i s_\alpha \widetilde \omega^{ij} \nabla_j s_\beta
= \nabla_\beta s_\alpha - \nabla_\alpha s_\beta
\end{equation}
for all $\alpha > \beta$
or
\begin{equation}\label{eq:masterincoord}
\frac{1}{2}
(\nabla_i s_\alpha \widetilde \omega^{ij} \nabla_j s_\beta)
f_\alpha^* \wedge f_\beta^* =
(\nabla_\beta s_\alpha)  f_\alpha^* \wedge f_\beta^*
\end{equation}
where $f_\alpha^*$ is the dual frame of
$\{\frac{\del}{\del q^1}, \cdots, \frac{\del}{\del q^{n-k}}\}$
defined by (\ref{eq:fbeta*}).
We call any of the two the master equation for the pre-symplectic
manifold $(Y,\omega)$.
\end{thm}

\begin{rem} It is straightforward to check that
both sides of (\ref{eq:masterincoord}) are independent of the local
lifting $E^\sharp$ but depends only on the splitting $\Pi: TY = G
\oplus T\FF$. In fact, the right hand side does not depend on this
splitting either but depends only on the null foliation $\FF$,
which is nothing but $-d_\FF(s)$ where $d_\FF$ is the exterior
derivative of $s \in \Omega^1(\FF) = \Gamma(T^*\FF)$ along the
foliation (See e.g., [Ha] or section 8 later for the definition of
$d_\FF$). The latter property follows from Remark \ref{Esharp}.
On the other hand, the invariant meaning of the left hand side
is not completely clear. One may interpret it as a kind of
twisted bracket which varies along the value of the section
itself.
\end{rem}

Note that (\ref{eq:masterincoord}) involves terms of all order of
$s_\beta$ because the matrix $(\widetilde \omega^{ij})$
is the inverse of the matrix
$$
\widetilde \omega_{ij} = \omega_{ij} - s_\beta F^\beta_{ij}.
$$
There is a special case where the curvature vanishes
i.e., satisfies
\begin{equation}\label{eq:special}
F_G = F^\beta_{ij}\frac{\del}{\del q^\beta}\otimes dy^i \wedge
dy^j= 0
\end{equation}
in addition to (\ref{eq:leaves}). In this case, $\widetilde
\omega_{ij} = \omega_{ij}$ which depends only on $y^i$'s and so
does $\omega^{ij}$. Therefore (\ref{eq:masterincoord}) is reduced
to the quadratic equation
\begin{equation}\label{eq:quad-master}
\frac{1}{2}(\nabla_i s_\alpha \omega^{ij} \nabla_j s_\beta)
f_\alpha^* \wedge f_\beta^* = (\nabla_\beta s_\alpha) f_\alpha^*
\wedge f_\beta^*.
\end{equation}

\section{Hamiltonian equivalence}
\label{sec:Ham}

In this section, we give the proof of Theorem \ref{extension} and
clarify the relation between the {\it intrinsic} equivalence
between the pre-symplectic structures and the {\it extrinsic}
equivalence between coisotropic embeddings in $U$.  The intrinsic
equivalence is provided by the pre-Hamiltonian diffeomorphisms on
the pre-symplectic manifold $(Y,\omega)$ and the extrinsic ones by
Hamiltonian deformations of its coisotropic embedding into
$(U,\omega_U)$ (and so into any symplectic $(X,\omega_X)$ that
allows a coisotropic embedding of $(Y,\omega))$, as far as the
Hamiltonian deformations are small enough).

We first recall the theorem stated in section
\ref{sec:neighborhoods} and provide its proof.

\begin{thm}
Any locally pre-Hamiltonian (respectively, pre-Hamiltonian) vector
field $\xi$ on a pre-symplectic manifold $(Y,\omega)$ can be
extended to a locally Hamiltonian (respectively, Hamiltonian)
vector field on the thickening $(U,\omega_U)$.
\end{thm}
\begin{proof}
Let $\FF$ be the null foliation of $(Y,\omega)$ and $\pi^{-1}(\FF)$ be
the induced foliation on the canonical thickening $\pi: U \to Y$.
Choose any splitting $TY = G\oplus E$ and write
$$
\xi =  \xi_G + \xi_E.
$$
First consider the case when $\xi$ is locally pre-Hamiltonian,
i.e., $d(\xi \rfloor \omega) = 0$. The pre-Hamiltonian case will follow
immediately from the proof of this.
Since $\omega|_E \equiv 0$, $\xi\rfloor \omega
= \xi_G \rfloor \omega$ and so we
have $d(\xi_G \rfloor \omega) = 0$. We denote by
$$
G^\sharp: = (T(\pi^{-1}(\FF)))^{\omega_U}
$$
and decompose $TU = G^\sharp \oplus T(\pi^{-1}(\FF))$. We will
find a vector field $\Xi$ on $U$ of the form
$$
\Xi = \Xi_G + \Xi_E
$$
where $\Xi_G(\alpha) \in G^\sharp_\alpha$ and $\Xi_E(\alpha) \in
T_\alpha (\pi^{-1}(\FF))$.

First we determine $\Xi_E$. We define a function $f_\xi: U \to \R$
by
$$
f_\xi(\alpha) = \langle \alpha, \xi_E \rangle
$$
for $\alpha \in U \subset E^*$ and by $X_{f_\xi}$ its Hamiltonian
vector field on $U$ with respect to the symplectic form
$\omega_U$. Note that $X_{f_\xi}|_Y = \xi_E$. Motivated by this,
we just set
\begin{equation}\label{eq:XiE}
\Xi_E = X_{f_\xi}.
\end{equation}

Determination of $\Xi_G$ is now in order. With the above
choice of $\Xi_E$, we have
\begin{equation}\label{eq:dwidetildexi}
d(\Xi \rfloor \omega_U) = d(\Xi_G \rfloor \omega_U).
\end{equation}
We would like to solve the equation
\begin{equation}\label{eq:XiG}
\begin{cases}
& d(\Xi_G \rfloor \omega_U) = 0 \\
& \Xi_G|_Y = \xi_G
\end{cases}
\end{equation}
Note the initial condition and the condition $\Xi_G \in G^\sharp$
imply
$$
(\Xi_G \rfloor \omega_U)|_Y = \xi \rfloor \omega.
$$
Noting that $\omega_U$ is nondegenerate, we will therefore set
$$
\Theta: = \Xi_G \rfloor \omega_U
$$
and solve the following extension problem of one-forms instead.
\begin{equation}\label{eq:Theta}
d\Theta = 0, \, \quad \Theta|_Y = \xi \rfloor \omega.
\end{equation}
This can be solved in a neighborhood of $Y \subset E^*$ by
the standard homotopy method [We1], when the initial one-form
$\xi \rfloor \omega$ on $Y \subset U$ is closed
which is precisely the condition
for $\xi$ to be locally pre-Hamiltonian. Furthermore the homotopy
method also leads to an exact extension when the initial form
$\xi \rfloor \omega$ is exact, i.e., when $\xi$ is pre-Hamiltonian.
This finishes the proof.
\end{proof}

Next, we consider the coordinate expression of (\ref{eq:XiG}). In
the canonical coordinates
$$
(y^1,\cdots, y^{2k}, q^1, \cdots,
q^{n-k}, p_1, \cdots, p_{n-k}),
$$
we write
$$
\Xi_G = \Xi_G^j e_j = \Xi_G^j \Big( \frac{\del}{\del y^j} + R_j^\alpha
\frac{\del}{\del q^\alpha} - p_\beta \frac{\del R_j^\beta}{\del
q^\nu}\frac{\del}{\del p_\nu}\Big)
$$
where $\Xi_G^j$'s are the coordinate functions of $\Xi_G$ with respect
to the basis $\{e_j\}_{1 \leq j \leq 2k}$ of $G^\sharp$. Then a
straightforward calculation using the formula (\ref{eq:omegaU})
shows that (\ref{eq:XiG}) becomes
\begin{equation}\label{eq:XiGincoord}
\begin{cases}
& d(\Xi^i_G(\omega_{ij} - p_\beta F^\beta_{ij}) dy^j) = 0 \\
& \Xi^i_G(y^1, \cdots, y^{2k}, q^1, \cdots, q^{n-k}, 0, \cdots, 0) =
\xi^i_G(y^1, \cdots, y^{2k}, q^1,\cdots, q^{n-k})
\end{cases}
\end{equation}
where $\xi^i_G$ is the coordinates of $\xi_G$
$$
\xi_G = \xi^i_G
\Big(\frac{\del}{\del y^i} + R_i^\alpha\frac{\del}{\del q^\alpha}\Big).
$$
In particular we have the following immediate corollary.

\begin{cor}\label{F=0case}
Suppose $F_\Pi \equiv 0$. Then $\Xi = \xi_G^\sharp + X_{f_\xi}$
is the locally Hamiltonian (respectively Hamiltonian)
vector field on $(U,\omega_U)$ which extends
the locally pre-Hamiltonian (respectively pre-Hamiltonian)
vector field $\xi$ on $(Y,\omega)$.
\end{cor}
\begin{proof}
When $F_\Pi \equiv 0$, (\ref{eq:XiGincoord}) becomes
$$
d(\Xi^i_G\omega_{ij}dy^j) = 0.
$$
Since $\omega_{ij}$ depends only on $y^i$'s and
$d(\xi_G^i\omega_{ij}dy^j) = 0$, the $p_\alpha$-independent
function
$$
\Xi^i_G = \xi_G^i
$$
provides a solution for (\ref{eq:XiGincoord}), which precisely
corresponds to the extension
$$
\Xi = \xi_G^\sharp + X_{f_\xi}.
$$
This finishes the proof.
\end{proof}

\section{Strong homotopy Lie algebroid}
\label{sec:algebroid}

In this section, we unravel the algebraic structure that provides
an invariant description of the master equation
(\ref{eq:masterincoord}) in the formal level.

\subsection{Lie algebroid and its cohomology}
\label{subsec:liealgebroid}

We start with recalling the definition of Lie algebroid and its
associated $E$-de Rham complex and $E$-cohomology. The leafwise de
Rham complex $\Omega^\bullet(\FF)$ is a special case of the {\it
$E$-de Rham complex} associated to the general Lie algebroid $E$,
We quote the following definitions from [NT].

\begin{defn}\label{liealgebroid} Let $M$ be a smooth manifold. A
{\it Lie algebroid} on $M$ is a triple $(E,\rho, [\, ,\, ])$,
where $E$ is a vector bundle on $M$, $[\, ,\, ]$ is a Lie algebra
structure on the sheaf of sections of $E$, and $\rho$ is a bundle
map, called the {\it anchor map},
$$
\rho: E \to TM
$$
such that the induced map
$$
\Gamma(\rho): \Gamma(M;E) \to \Gamma(TM)
$$
is a Lie algebra homomorphism and, for any sections $\sigma$ and
$\tau$ of $E$ and a smooth function $f$ on $M$, the identity
$$
[\sigma, f\tau] = \rho(\sigma)[f]\cdot \tau + f\cdot [\sigma,
\tau].
$$
\end{defn}

\begin{defn}\label{E-deRham} Let
$(E,\rho,[\, ,\, ])$ be a Lie algebroid on $M$. The {\it $E$-de Rham
complex} $(^E\Omega^\bullet(M), ^Ed)$ is defined by
\begin{eqnarray*}
^E\Omega(\Lambda^\bullet(E^*)) & = &\Gamma(\Lambda^\bullet(E^*))\\
^Ed\omega(\sigma_1, \cdots, \sigma_{k+1})
& = &\sum_i(-1)^i \rho(\sigma_i)\omega(\sigma_1, \cdots,
\widehat{\sigma_i}, \cdots, \sigma_{k+1}) \\
&\qquad &+
\sum_{i<j}(-1)^{i+j-1}\omega([\sigma_i,\sigma_j],\sigma_1,
\cdots,\widehat \sigma_i, \cdots, \widehat \sigma_j, \cdots,
\sigma_{k+1}).
\end{eqnarray*}
The cohomology of this complex will be denoted by $^EH^*(M)$ and
called the $E$-de Rham cohomology of $M$.
\end{defn}

In our case, $M = Y$ and $E: = TY^\omega=T\FF$ and the anchor map
$\rho: E \to TY$ is nothing but the inclusion map $i: TY^\omega
\to TY$. The integrability of $TY^\omega$ implies that the
restriction of the Lie bracket on $\Gamma(TY)$ to
$\Gamma(TY^\omega)$ defines the Lie bracket $[\, ,\, ]$ on
$\Gamma(E)$. Therefore the triple
$$
(E=TY^\omega, \rho = i, [\, ,\, ])
$$
defines the structure of Lie algebroid and hence the
$E$-differential. In our case, the
corresponding $E$-differential is nothing but $d_\FF$ the exterior
derivative along the null foliation $\FF$ and its cohomology,
the cohomology $H^*(\FF)$ of the foliation $\FF$.

Now we derive, as a first step towards the study of the master
equation, the linearized version of the master equation which
characterizes the infinitesimal deformation space of coisotropic
submanifolds. For this, we introduce the space
$$
{\mathcal Coiso}_k = {\mathcal Coiso}_k(X,\omega_X)
$$
the set of coisotropic submanifolds with nullity $n-k$ for $0\leq
k \leq n$ and characterize its infinitesimal deformation space at
$Y \subset E^*$, the zero section of $E^*$.
As Zambon [Za] observed (see Example 10.4 later), this space
does not form a smooth Frechet manifold unlike the Lagrangian case.
By the coisotropic neighborhood theorem, the infinitesimal
deformation space, denoted as  $T_{Y}{\mathcal
Coiso}_k(X,\omega_X) = T_{Y}{\mathcal Coiso}_k(U,\omega_U)$ with
some abuse of notion, depends only on $(Y, \omega)$ where $\omega
= i^*\omega_X$, but not on $(X,\omega_X)$. An element in
$T_{Y}{\mathcal Coiso}_k(U,\omega_U)$ is a section of the bundle
$E^*= T^*\FF \to Y$.

The following characterizes the condition for a section $\xi$ to
be an infinitesimal deformation of $Y$.

\begin{thm}\label{master1}
A section $\xi \in \Omega^1(\FF)$ is an infinitesimal deformation of
coisotropic submanifold $Y$ if and only if $\xi$ satisfies \be
\label{eq:master1} d_{\FF}( \xi) =0 \quad \mbox{on } \,
\Omega^2(\FF), \ee where $d_\FF$ is the exterior derivative
along the leaves of the foliation $\FF$.
\end{thm}
\begin{proof}
This is an immediate consequence of higher order nature of the
left hand side in (\ref{eq:masterincoord}) and the definition of
$d_{\FF}$.
\end{proof}

Now we need to mod out the solution space of (\ref{eq:master1}) by
certain gauge equivalence classes. Algebraically, the set of
equivalence classes is the first cohomology of the null foliation $\FF$
$$
H^1(\FF) : = \ker d_{\FF}|_{\Omega^1(\FF)}/
\mbox{im }d_{\FF}|_{\Omega^0(\FF)}
$$
of the graded algebra $(\Omega^\bullet(\FF), d_{\FF})$, which is
the leafwise de Rham complex of the foliation $\FF$. Geometrically
we are considering the set of equivalence classes of
pre-symplectic structures up to the pre-Hamiltonian diffeomorphisms
on $(P,\omega)$.  Or equivalently the latter is (locally)
equivalent to the set of deformations of coisotropic embeddings
into $(U,\omega_U)$ (or any into $(X,\omega_X)$ into which
$(Y,\omega)$ is coisotropically embedded) of $(Y,\omega)$ up to the
ambient Hamiltonian isotopies. See section \ref{sec:Ham}.

Since the leafwise de Rham cohomology $H^\bullet(\FF)$ is determined
by $(Y,\omega)$, we will
denote
$$
H^\bullet(Y,\omega): = H^\bullet(\FF)
$$
and call the {\it cohomology of (the null foliation of)
$(Y,\omega)$}. In fact, it depends only on the foliation $\FF$.
When $k = 0$, this reduces to the standard description of the
deformation problem of the Lagrangian submanifolds and
$H^1(Y,\omega)$ becomes the standard de Rham cohomology
$H^1(Y;\R)$. When $k = n$, it becomes $C^\infty(X,\omega)$ which
carries the Poisson bracket associated to the symplectic form
$\omega$. The cohomology $H^1(Y,\omega)$ is a mixture of these
two. Except the case of $k = 0$, this cohomology is infinite
dimensional in general.

\subsection{Review of strong homotopy Lie algebra}
\label{subsec:Linftylagebra}

In this section we give a brief review of the definition of the
strong homotopy Lie algebra or the $L_\infty$ algebra. We will
follow Fukaya's exposition verbatim on the part of $L_\infty$
algebra from [Fu]. In particular, we will follow the sign
convention used in [Fu] which is also the same as that of [FOOO]
in the context of the $A_\infty$-algebra case. We will extract
only the essentials that are needed to give a self-contained
definition of our strong homotopy Lie algebroid. For more details
of the discussion on the $L_\infty$ algebra, we refer to [Fu] or
[K1] (with different sign convention) also to [AKKS], [P1] for the
exposition in the context of the Batalin-Vilkovisky approach using
the notion of super-manifolds.

Let $C$ be a graded $R$-module where $R$ is the coefficient ring.
In our case, $R$ will be either $\R$ or $\C$. We denote by $C[1]$
its suspension defined by
$$
C[1]^k = C^{k+1}.
$$
We denote by $deg(x)=|x|$ the degree of $x \in C$ before the shift and
$deg'(x)=|x|'$ that after the degree shifting, i.e., $|x|' = |x| - 1$.
Define the {\it bar complex} $B(C[1])$ by
$$
B_k(C[1]) = (C[1])^{k\otimes}, \quad B(C[1]) =
\bigoplus_{k=0}^\infty B_k(C[1]).
$$
Here $B_0(C[1]) = R$ by definition. We provide the degree of
elements of $B(C[1])$ by the rule
\begin{equation}\label{eq:degonBC[1]}
|x_1 \otimes \cdots \otimes x_k|': = \sum_{i=1}^k |x_i|' = \sum_{i =1}^k|x_i|
-k
\end{equation}
where $|\cdot|'$ is the shifted degree. There is a natural
coproduct
$$
\Delta: B(C[1]) \to B(C[1]) \otimes B(C[1])
$$
on $B(C[1])$ defined by
\begin{equation}\label{eq:coproduct}
\Delta(x_1 \otimes \cdots x_k) = \sum_{i = 0}^k (x_1 \otimes
\cdots \otimes x_i) \otimes (x_{i+1} \otimes \cdots \otimes x_k).
\end{equation}

Finally we have the natural projection
$$
\e : B(C[1]) \to (B(C[1])^0\cong R.
$$
Then the triple $(B(C[1]), \Delta, \e)$ defines the structure of
{\it graded coalgebra}.

Now we consider the action of the symmetric group $S_k$ on
$B_k(C[1])$ by the permutation of arguments in the tensor powers
\begin{equation}\label{eq:sigmaaction}
\sigma(x_1 \otimes\cdots  \otimes x_k) = (-1)^{|\sigma|}
x_{\sigma(1)} \otimes \cdots \otimes x_{\sigma(k)}
\end{equation}
where $|\sigma|$ is defined to be
\begin{equation}\label{eq:sigmasign}
|\sigma| = \sum_{i,j; i < j, \, \sigma(i) > \sigma(j)}
|x_i|'|x_j|'.
\end{equation}
We define $E_k(C[1])$ to be the submodule of $B_k(C[1])$
consisting of fixed points of the $S_k$-action defined above, and
$$
E(C[1]) = \bigoplus_{i=0}^\infty E_k(C[1]).
$$
The above coproduct (\ref{eq:coproduct}) naturally induces a
coproduct on $E(C[1])$, which we also denote by $\Delta$. And the
projection $\e$ induces the projection on $E(C[1])$ again denoted
by $\e$. Then we have the following lemma

\begin{lem}{\bf [Lemma 8.3.1, Fu]} The triple $(E(C[1]), \Delta,
\e)$ is a graded-cocommutative coalgebra.
\end{lem}

\begin{defn}{\bf [Definition 8.3.2, Fu]} The structure of
{\it $L_\infty$ algebra} or {\it strong homotopy Lie algebra} is a
sequence of $R$ module homomorphisms
$$
\frak m_k: E_k(C[1]) \to C[1], \quad k = 1, 2, \cdots,
$$
of degree +1 such that the coderivation
$$
\delta = \sum_{k=1}^\infty \widehat{\frak m}_k
$$
satisfies $\delta \delta = 0$. Here we denote by $\widehat{\frak
m}_k: E(C[1]) \to E(C[1])$ the unique extension of $\frak m_k$
as a coderivation on $E(C[1])$.
\end{defn}
One can write the condition $\delta\delta = 0$ more explicitly as
\begin{eqnarray}\label{eq:Linfty}
&&\sum_{k=1}^n\sum_{i =1}^{n-k+1} (-1)^{|x_1|' + \cdots +
|x_{i-1}|'} \frak m_{n-k+1}(x_1 \otimes \cdots \otimes
x_{i-1}\otimes\frak m_k(x_i \otimes
\cdots \otimes x_{i+k-1}) \nonumber \\
&\quad & \hskip1in
\otimes x_{i+k} \otimes \cdots \otimes x_n) = 0 \nonumber
\end{eqnarray}
for any $n = 1, \, 2, \cdots$. In particular, we have $\frak m_1
\frak m_1 = 0$ and so it defines complex $(C,\frak m_1)$. We
define the $\frak m_1$-cohomology by
$$
H(C,\frak m_1) = \mbox{ker }\frak m_1/\mbox{im }\frak m_1.
$$

A {\it weak $L_\infty$-algebra} is defined in the same way, except
that it also includes the $\frak m_0$-term
$$
\frak m_0: R \to E(C[1]).
$$
The first two terms of the $L_\infty$ relation for a weak
$L_\infty$ algebra are given as
\begin{eqnarray}
\frak m_1(\frak m_0(1)) & = & 0 \nonumber \\
\frak m_1\frak m_1 (x) + (-1)^{|x|'}\frak m_2(x, \frak m_0(1)) +
\frak m_2(\frak m_0(1), x) & = & 0. \label{eq:m0m1}
\end{eqnarray}
In particular, for the case of weak $L_\infty$ algebras, $\frak
m_1$ will not satisfy boundary property, i.e., $\frak m_1\frak m_1
\neq 0$ in general. We will explain in the appendix of the present
paper that the vanishing of $\frak m_0$ is closely related to the
coisotropic boundary condition via the Batalin-Vilkovisky
formulation of the open string $A$-model. See also Lemma
\ref{brane}.

\subsection{Definition of strong homotopy Lie algebroid}

Now we are ready to introduce the main definition of our {\it strong
homotopy Lie algebroid}

\begin{defn} Let $E \to Y$ be a Lie algebroid. An
{\it $L_\infty$-structure} over the
Lie algebroid is a structure of strong homotopy Lie algebra
$(\frak l[1], \frak m)$ on the associated $E$-de Rham complex
$\frak l^\bullet = \Omega^\bullet(E)
= \Gamma(\Lambda^\bullet(E^*))$
such that $\frak m_1$ is the $E$-differential $^Ed$ induced by
the Lie algebroid structure on $E$ as described
in section \ref{liealgebroid}.
We call the pair $(E \to Y, \frak m)$ a strong homotopy Lie algebroid.
\end{defn}

With this definition of strong homotopy Lie algebroid, we will
show that for given presymplectic manifold $(Y,\omega)$
each splitting $\Pi: TY = G\oplus T\FF$ induces a
canonical $L_\infty$-structure over the Lie algebroid $T\FF \to Y$.

We recall that $U \subset E^* = T^*\FF$ and
we have chosen the (locally defined) transverse
symplectic connection
$$
TU = (G^\sharp \oplus E^\sharp)\oplus VTU
$$
where $E^\sharp$ is given as in (\ref{eq:Esharp}). This induces a
locally defined leafwise symplectic
connection of $TU \to Y$ which we denote by $\nabla$.

The crucial structures relevant to the invariant description of
this structure will be a linear map
\begin{equation} \label{eq:lin}
\widetilde \omega: \Omega^1(Y;\Lambda^\bullet E^*) \to
\Gamma(\Lambda^{\bullet +1}E^*) = \Omega^{\bullet + 1}(\FF),
\end{equation}
a quadratic map \be\label{eq:quad} \langle \cdot, \cdot
\rangle_{\omega}: \Omega^1(Y;\Lambda^{\ell_1} E^*) \otimes
\Omega^1(Y;\Lambda^{\ell_2} E^*) \to \Omega^{\ell_1 +
\ell_2}(\FF), \ee and the third map that is induced by the
transverse $\Pi$-curvature.

Now we describe those maps. The linear map $\widetilde \omega$ is
defined by
$$
\label{eq:lin} \widetilde\omega(A): = (A|_E)_{skew}:
$$
Here note that an element $A \in \Omega^1(Y; \Lambda^kE^*)$ a
section of $T^*Y\otimes \Lambda^k E^*$. $A|_E$ is the element in
$E^* \otimes \Lambda^k E^*$ obtained by restricting $A$ to $E$ for
the first factor, and then $(A|_E)_{skew}$ is the
skew-symmetrization of $A|_E$.

The quadratic map is defined by
$$
\langle A,B \rangle_{\omega}:= \langle A|\pi |B \rangle - \langle
B|\pi|A \rangle
$$
where $\pi$ is the transverse Poisson bi-vector on $N^*\FF$ associated to the
transverse symplectic form $\omega$ on $N\FF$.

Finally we define two maps involving the transverse $\Pi$-curvature
$F_\Pi$: For any two form $\eta \in\Omega^2(Y)$, we denote
$$
\ker \eta:=\{ v \in TY \mid \eta(v, \cdot) \equiv 0 \},
$$
and the subset $\Omega^2_\omega(Y) \subset \Omega^2(Y)$ by
\begin{equation}\label{eq:rangeF}
\Omega^2_\omega(Y) :=\{ \eta \in \Omega^2(Y)\mid
\ker \eta \supset \ker \omega \}.
\end{equation}
The first map we will use is the ``contraction'' by $F_\Pi$
\begin{equation}\label{eq:F}
\widetilde F: \Omega^\ell(\FF) \to
\Omega^2_\omega(Y)\otimes \Omega^{\ell-1}(\FF); \quad
\widetilde F(\xi) =  F \rfloor \xi,
\end{equation}
where the contraction is taken between $E$ and $E^*$. The second
map is ``raising indices'' of $F$ by $\omega^{-1}=(\omega^{ij})$
on $TY/E \cong G$. We will denote
\begin{equation}\label{eq:Fsharp}
F^\#:= F\omega^{-1} = F^{\alpha j}_i dy^i \otimes
\Big(\frac{\del}{\del y^j} + R^\beta_j \frac{\del}{\del
q^\beta}\Big) \otimes \frac{\del}{\del q^\alpha} \in \Gamma(G^*
\otimes G \otimes E),
\end{equation}
where $F^{\alpha j}_i = F^\alpha_{ik}\omega^{kj}$. Note that we
can identify $\Gamma(G^* \otimes G \otimes E)$ with
$\Gamma(N^*\FF\otimes N\FF \otimes E)$ via the isomorphism $\pi_G:
G \to N\FF$.

For given $\xi \in \Omega^\ell(\FF)$, we denote
\begin{equation}\label{eq:lin}
d_{\FF}(\xi): = (\nabla \xi|_E)_{skew},
\end{equation}
\begin{equation}\label{eq:quad} \{\xi_1, \xi_2
\}_{\Pi}: = \langle \nabla \xi_1, \nabla \xi_2 \rangle_{\omega} =
\sum_{i < j} \omega^{ij}(\nabla_i \xi_1) \wedge (\nabla_j\xi_2 ).
\end{equation}
Here the first map is nothing but the leafwise differential of the
null foliation which is indeed independent of the choice of
splitting $\Pi: TY = G \oplus T\FF$ but depends only on the
foliation. The second is a bracket in the transverse direction
which does not satisfy the Jacobi identity in general because of
the presence of {\it non-zero} transverse curvature of the null
foliation. Dependence on $\Pi$ for the bracket comes from the
covariant derivative $\nabla_i s:= \Pi^v\circ ds(e_i)$. Because of this,
the structure
\begin{equation}\label{eq:dga}
\left(\frak l = \bigoplus_{j=0}^{n-k}\frak l^j; d_{\FF}, \{\cdot,
\cdot \}_{\Pi}\right), \quad \frak l^j= \Omega^j(\FF)
\end{equation}
fails to define a differential graded Lie algebra in
general. However we have

\begin{thm}\label{dga}
Let $\Pi: TY = G\oplus T\FF$ be a splitting such that $F_\Pi
\equiv 0$. Then (\ref{eq:dga}) defines a differential graded
Lie algebra. More precisely, we have the identity
\begin{equation}\label{eq:dgaidentity}
 d_\FF \{\xi_1,
\xi_2\}_\Pi= \{d_\FF(\xi_1), \xi_2\}_\Pi+ (-1)^{|\xi_1|}\{\xi_1,
d_\FF (\xi_2)\}_\Pi.
\end{equation}
\end{thm}
\begin{proof}
Let $\Pi: TY = G\oplus T\FF$ be such a splitting. By definition of
$F_\Pi$, vanishing of $F_\Pi$ means that the distribution $G$ is
integrable. Together with the foliation $\FF$, we can construct a
coordinate system
$$
(y^1, \cdots, y^{2k}, q^1, \cdots, q^{n-k})
$$
such that
$$
G = \mbox{span}\Big\{\frac{\del}{\del y^1}, \cdots,
\frac{\del}{\del y^{2k}}\Big\}, \quad T\FF = \mbox{span} \Big\{
\frac{\del}{\del q^1}, \cdots, \frac{\del}{\del q^{n-k}} \Big\}.
$$
In particular, we have $R^\alpha_i \equiv 0$ and so $[\nabla_i,
\nabla_j] = [\nabla_k, \nabla_\alpha] = [\nabla_\alpha,
\nabla_\beta] = 0$. Once we have these, (\ref{eq:dgaidentity})
immediately follows by computing
$$
d_\FF \{\xi_1, \xi_2\}_\Pi= dq^\alpha \wedge \nabla_\alpha
(\frac{1}{2}\omega^{ij} \nabla_i\xi_1 \wedge \nabla_j \xi_2).
$$
using the fact that $\omega^{ij}$ is independent of $q^\alpha$'s.
This finishes the proof.
\end{proof}

In the general case, it turns out that for each given
splitting $\Pi: TY = G\oplus T\FF$ the pair $(d_{\FF},
\{\cdot, \cdot\}_\Pi)$ can be extended to an infinite family of
graded multilinear maps
\begin{equation}\label{eq:frakm}
\frak m_\ell =
\frak m_\ell^\Pi:(\Omega[1]^\bullet(\FF))^{\otimes\ell} \to
\Omega[1]^{\bullet}(\FF)
\end{equation}
so that the structure
$$
\left(\bigoplus_{j=0}^{n-k} \frak l[1]^j; \{\frak m_\ell\}_{1 \leq
\ell < \infty}\right)
$$
defines a {\it strong homotopy Lie algebroid} on $E=T\FF \to Y$ in
the above sense.  Here  $\Omega[1]^\bullet(\FF)$ is the shifted complex of
$\Omega^\bullet(\FF)$, i.e., $\Omega[1]^k(\FF) = \Omega^{k+1}(\FF)$
and $\frak m_1$ is defined by
$$
\frak m_1(\xi) = (-1)^{|\xi|} d_{\FF}(\xi)
$$
and $\frak m_2$ is given by
$$
\frak m_2(\xi_1, \xi_2) =
(-1)^{|\xi_1|(|\xi_2|+1)}\{\xi_1,\xi_2\}_\Pi.
$$
On the un-shifted group $\frak l$, $d_{\FF}$ defines a
differential of degree 1 and $\{\cdot, \cdot\}_{\omega}$ is a
graded bracket of degree 0 and $\frak m_\ell$ is a map of degree
$2-\ell$.

We now define $\frak m_\ell$ for $\ell \geq 3$. Here enters the
transverse $\Pi$-curvature $F=F_\Pi$ of the splitting $\Pi$ of the
null foliation $\FF$. We define
\begin{equation}\label{eq:mell}
\frak m_\ell(\xi_1, \cdots, \xi_\ell) : =
\sum_{\sigma \in S_\ell}
(-1)^{|\sigma|}\langle \nabla \xi_{\sigma(1)},
(F^\#\rfloor \xi_{\sigma(2)}) \cdots
(F^\#\rfloor \xi_{\sigma(\ell-1)})
\nabla \xi_{\sigma(\ell)}\rangle_\omega
\end{equation}
where $|\sigma|$ is given by the rule written in
(\ref{eq:sigmaaction}). We have obtained our definition of strong
homotopy Lie algebroid associated to the coisotropic submanifolds
first via the language of super-manifolds and the
Batalin-Vilkovisky formalism in the context of coisotropic branes
formulated in [AKSZ] and [P1]. The current definition below is a
literal translation of the one formulated in the language of
formal super manifolds into the language of tensor calculus. For
the reader's convenience, we include this original derivation in
the context of Batalin-Vilkovisky formalism in the appendix. In
fact, the current proof itself is not very different from this in
which we just add more mathematical explanations to it. The proof
is essentially a consequence of the fact that the symplectic form
$\omega_U$ is closed. A purely classical proof using the tensor
calculus should be also possible, but with paying the price of
obscuring the origin, $d\omega_U = 0$, of the
$L_\infty$-structure. Because of this, we do not pursue carrying
out the tensor calculations preferring the proof using the
super-manifold and the relevant super-calculus.

\begin{thm}\label{algebroid}
Let $(Y,\omega)$ be a pre-symplectic manifold and $\Pi: TY =
G\oplus T\FF$ be a splitting. Then $\Pi$ canonically induces a
structure of strong homotopy Lie algebroid on $T\FF$ in that the
graded complex
$$
\left(\bigoplus_\bullet \Omega[1]^\bullet(\FF),
\{\frak m_\ell \}_{1 \leq \ell < \infty}\right)
$$
defines the structure of strong homotopy Lie algebra. We denote by
$\frak l[1]^\infty_{(Y, \omega;\Pi)}$ the corresponding strong
homotopy Lie algebra.
\end{thm}
\begin{proof}
We will prove that the coderivation $\delta = \sum_{\ell=1}^\infty
\widehat{\frak m}_\ell$ satisfies $\delta \delta = 0$.
By restricting to the sections $\xi_i$'s
supported in a coordinate chart,
we will work with coordinate calculations.
We first recall from (\ref{eq:omegaU})
the expression of the symplectic form $\omega_U$ on $U \subset E^*$
\begin{eqnarray}\label{eq:omega}
\omega & = & \frac{1}{2}\Big(\omega_{ij} - p_\beta F^\beta_{ij}
\Big)
dy^i \wedge dy^j \nonumber \\
& \quad &  - (dp_\beta + p_\beta \frac{\del R_i^\beta}{\del
q^\gamma} dy^i) \wedge (dq^\delta - R_j^\delta dy^j)
\end{eqnarray}
In the same coordinates, the corresponding Poisson
bi-vector field $P$ has the form
\begin{equation}\label{eq:Poisson}
P = \frac{1}{2}\widetilde\omega^{ij}e_i \wedge e_j
+ \frac{\del}{\del q^\alpha} \wedge \frac{\del}{\del p_\alpha}
\end{equation}
where the vector field $e_j$ are the ones defined by (\ref{eq:ej})
(see Appendix for the derivation of this formula).
The fact that $d\omega_U =0$
is equivalent to the vanishing of Schouten bracket
\begin{equation}\label{eq:[P,P]=0}
[P, P] = 0.
\end{equation}
The Poisson tensor defines a map
$$
\delta_P(u) = [P,u]
$$
for each multi-vector field $u$. Furthermore, since $[P,P] = 0$,
it satisfies $\delta_P \delta_P = \frac{1}{2}[ [P,P], \cdot] = 0$.
A good way of describing this map $\delta_P$ is to use the
super-language (see Appendix or [Gz] for an elegant description of
this translation). We change the parity of $TU$ along the fiber
and denote by $T[1]U$ the corresponding super tangent bundle of
$U$. One considers a multi-vector field on $U$ as a (fiberwise)
polynomial function on $T^*[1]U$. For example, the Poisson tensor
$P$ defines a quadratic function, which we denote by $H$. This
also coincides with the push-forward of the canonical even
function $H^*: T[1]U \to \R$ induced by the symplectic form
$\omega_U$. On the other hand, the exterior differential $d$
defines an odd vector field on $T[1]U$, which we denote by $Q$.
This vector field is nothing but the pull-back of the Hamiltonian
vector field of $H$ with respect to the {\it canonical (odd)
symplectic form} $\Omega$ on $T^*[1]X$. We warn readers that
$\Omega$ should  {\it not} be confused with the symplectic form
$\omega_U$ itself on $U$. We denote by $\{\cdot, \cdot\}_\Omega$
the (super-)Poisson bracket associated to the odd symplectic form
$\Omega$ on $T[1]X$. Then we have the identity
$$
Q = \{H^*, \cdot\}_\Omega
$$
as a derivation on the set $\mathcal O_{T[1]X}$ of ``functions''
on $T[1]X$: Here $\mathcal O_{T[1]X}$ is the set of differential
forms on $X$ considered as fiberwise polynomial functions on
$T[1]X$. We refer to Appendix or [Gz] for the precise mathematical meaning for
this correspondence.

Next we will be interested in whether one can canonically restrict
the vector field $Q$ to $\L = T\FF[1]$ or equivalently whether the
function $H$ has constant value on $\L$. Here comes the
coisotropic condition naturally.

\begin{lem}\label{brane}
Let $H$ be the even function on $T[1]X$ induced by the symplectic
form $\omega_X$, and $H^*: T^*[1] X \to \R$ be its push-forward by
the isomorphism $\widetilde \omega_X: T[1]X \to T^*[1]X$. When $Y
\subset (X,\omega)$ is a coisotropic submanifold we have
$H^*|_{N^*[1]Y} = 0$. Conversely, any (conic) Lagrangian subspace
$\mathbb{L}^* \subset T^*[1]X $ satisfying $H^*|_{\mathbb{L}^*}=0$
is equivalent to $N^*[1]Y$, for some coisotropic submanifold $Y$
of $(X,\omega)$.
\end{lem}
\begin{proof}
This is essentially a translation of the definitions of
coisotropic submanifolds and the even function $H^*$.
We first note that $T_xY$ being coisotropic in $(T_xX, \omega_X(x))$
is equivalent to $N^*_xY$ being isotropic in $(T_x^*X, \pi_X(x))$
where $\pi_X$ is the Poisson tensor associated to $\omega_X$.
Then it follows from an easy super algebra that
this last statement is equivalent to the vanishing of the
associated even function $H^*$ on $N^*[1]Y$ at each point
$x \in Y$. This finishes the first part of the theorem.

For the proof of the converse, we first recall that any conic
Lagrangian submanifold in $T^*X$ has the form of $N^*Y$ for some
submanifold $Y \subset X$. Then the above argument shows that
vanishing of $H^*$ on $\L^* = N^*[1]Y$ is equivalent to $Y$ being
coisotropic. This finishes the proof.
\end{proof}

\begin{rem} Note that this lemma, as it is, applies to the
coisotropic submanifolds in Poisson manifolds. See [We2] for the
definition of coisotropic submanifolds in Poisson manifolds. In
our case of $U \subset E^*$, $\pi_X$ is nothing but $P$ above.
\end{rem}

Noting that $\L^* = N^*[1]Y$ is mapped to $\L = T\FF[1]$ under the
isomorphism $\widetilde \omega_X$, this lemma enables us to
restrict the odd vector field $Q$ to $T\FF[1]$. We need to
describe the Lagrangian embedding $T\FF[1] \subset T[1]X$ more
explicitly, and describe the induced directional derivative acting
on
$$
\Omega^\bullet(\FF)
$$
regarded as a subset of ``functions'' $T\FF[1]$. (Again we refer
to Appendix or [Gz] for the precise explanations of this). A more tensorial
way of saying this is as follows: Noting that the Poisson tensor
pairs with any section of $N^*\FF$ to give zero {\it due to the
coisotropic condition of $Y \subset U$}, the {\it Hamiltonian
operator} [GD]
$$
\delta_P:= [P, \cdot]
$$
restricts to $\Omega^\bullet(\FF) = \Gamma(\wedge^\bullet
(T^*\FF))$. We denote by $\delta'$ this restriction. More
precisely $\delta': \Omega^\bullet(\FF) \to \Omega^\bullet(\FF)$
is given by the formula
\begin{equation}\label{eq:delta'}
\delta'(\xi) = \{H, \widetilde \xi\}_\Omega\Big|_{\L}
\end{equation}
where $\widetilde \xi$ is the extension of $\xi$ in a neighborhood
of $\L \subset T[1]X$: the extension that we use is the lifting of
$\xi \in \Omega^\bullet(\FF)$ to an element of $\Omega^\bullet(U)$
obtained by the (local) Ehresman connection constructed in section
\ref{sec:master}. The condition $Q|_\L\equiv 0$ implies that this
formula is independent of the choice of  (local) Ehresman
connection. We will just denote $\delta'(\xi)= \{H, \xi\}_\Omega$
instead of (\ref{eq:delta'}) as long as there is no danger of
confusion.

Obviously, $\delta'$ satisfies $\delta'\delta' = 0$ because of
$\delta_P$ does. Now it remains to verify that this is translated
into the $L_\infty$ relation $\delta\delta = 0$ in the tensorial
language which is exactly what we wanted to prove. For this
purpose, we need to describe the map $\delta': \Omega^\bullet(\FF)
\to \Omega^\bullet(\FF)$ more explicitly.

Restricting ourselves to a Darboux neighborhood $\L = T\FF[1]
\subset T[1]U$, we identify the neighborhood with a neighborhood
of the zero section $T^*[1]\L$. Using the fact that
(\ref{eq:delta'}) depends only on $\xi$, not on the extension, we
will make a convenient choice of coordinates to write $H$ in the
Darboux neighborhood and describe how the derivation $Q = \{H,
\cdot \}_\Omega$ acts on $\Omega^*(\FF)$ in the canonical
coordinates of $T^*[1]\L$. In this way, we can apply the canonical
quantization which provides a canonical correspondence between
functions on ``the phase space'' $T^*[1]\L$ and the corresponding
operators acting on the functions on the ``configuration space''
$\L$, when we find out how $\delta'$ acts on
$\Omega^\bullet(\FF)$.

We denote by $(y^i,q^\alpha, p_\alpha, y^*_i, q^*_\alpha,
p^\alpha_*)$ the canonical coordinates $T^*\L$ associated with the
coordinates $(y^i,q^\alpha, p_\alpha)$ of $N^*\FF$. Note that
these coordinates are nothing but the canonical coordinates of
$N^*Y \subset T^*U$ pulled-back to $T\FF \subset TU$ and its
Darboux neighborhood, with the corresponding parity change: We
denote the (super) canonical coordinates of $T^*[1]\L$ associated
with $(y^i,q^\alpha\mid p_\alpha)$ by
$$
\Big(\begin{matrix}y^i, & q^\alpha  && \mid p_*^\alpha \\
y^*_i, & q^*_\alpha && \mid p_\alpha
\end{matrix}\Big)
$$
Here we note that the degree of $y^i,\, q^\alpha$ and $p_\alpha$
are $0$ while their anti-fields, i.e., those with $*$ in them have
degree $1$. And we want to emphasize that $\L$ is given by the
equation
\begin{equation}\label{eq:LL}
y_i^* = p_\alpha = p_*^\alpha = 0
\end{equation}\label{eq:LL}
and $(y^i, y^*_i)$, $(p_\alpha, q_\alpha^*)$ and $(p_*^\alpha,
q^\alpha)$ are conjugate variables.

The Poisson tensor $P$ (\ref{eq:Poisson}) becomes the even
function $H$ that has the form
\begin{eqnarray}
H =  \frac{1}{2} \widetilde\o^{ij}y_i^\# y_j^\# + p_*^\delta
q^*_\delta  \label{eq:H*}
\end{eqnarray}
in the canonical coordinates of $T^*[1]\L$. Here we define
$y_i^\#$ to be
$$
y_i^\#: = y_i^* + R_i^\delta p^\delta_* - p_\beta \frac{\del
R_i^\beta}{\del q^\delta}q^*_\delta.
$$
On the other hand, we have
$$
\widetilde \omega^{-1}_\alpha = \omega^{-1}_{\pi(\alpha)}
\sum_{\ell =0}^\infty (F^\# \rfloor \alpha)^\ell \quad \mbox{on }
\, TY/E
$$
(See (\ref{eq:omegaF}) in section \ref{sec:moduli} later.) which
is written as
$$
\widetilde\omega^{ij}_\alpha = \omega^{ij_0}_{\pi(\alpha)}
\sum_{\ell=0}^\infty
(p_{\beta_1}F^{\beta_1j_1}_{j_0})(p_{\beta_2}F^{\beta_2j_2}_{j_1})
\cdots (p_{\beta_\ell}F^{\beta_\ell j}_{j_{\ell-1}})
$$
in coordinates where $\alpha = p_\beta f_\beta^*$.
We make replacements
$$
y_i^* \mapsto \frac{\del}{\del y^i}, \, p_*^\alpha \mapsto
\frac{\del}{\del q^*_\alpha},  \, p^*_\alpha \mapsto
\frac{\del}{\del q^\alpha}
$$
following the canonical quantization process in the cotangent
bundle, and noting that the derivative $\frac{\del}{\del
q^*_\alpha}$ means the contraction by $\frac{\del}{\del q^\alpha}$
on $\Omega^\bullet(\FF)$.

Then by expanding the Poisson tensor $P$ or the even function $H$
above into the power series
$$
H = \sum_{\ell =1} H_\ell, \quad H_\ell \in \frak l^\ell,
$$
in terms of the degree (i.e., the number of
factors of odd variables $(y_i^*, p_*^\alpha, p_\alpha)$
or the `ghost number' in the physics language)
our definition of $\frak m$ exactly corresponds to the
$\ell$-linear operator
$$
(\xi_1, \xi_2, \cdots, \xi_\ell) \mapsto \{\cdots \{H_\ell,
\xi_1\}_\Omega, \cdots\}_\Omega, \xi_\ell\}_\Omega.
$$
Note that the above power series acting on
$(\xi_1, \cdots, \xi_\ell)$ always reduces to a finite sum and so is
well-defined as an operator. Then by definition, the coderivation
$$
\delta = \sum_{\ell =1}^\infty \widehat{\frak m}_\ell
$$
precisely corresponds to $\delta' = \{H, \cdot\}_\Omega$. The
$L_\infty$ relation $\delta \delta = 0$ then immediately follows
from $\delta' \delta' = 0$. This finishes the proof.
\end{proof}
For example, under the above translation,
the odd vector field
\begin{equation}\label{eq:q}
m_1 = \CQ \mid_\L
\end{equation}
acts on
$$
\frak l = \bigoplus_{\ell = 0}^{n-k}\frak l^\ell \cong
\bigoplus_{\ell = 0}^{n-k}\Omega^\ell(\FF)
$$
is translated into to the leafwise differential $d_{\FF}$.

\section{Gauge equivalence and formality question}
\label{sec:gauge}

In this section, based on the Lemma \ref{splittings}, we prove
that two strong homotopy Lie algebroids we have associated to two
different splittings are {\it gauge equivalent} or {\it
$L_\infty$-isomorphic}. In fact, there exists a {\it canonical}
$L_\infty$ isomorphism between the two which depends only on
$B_{\Pi_0\Pi}$ given in Lemma \ref{splittings}.

We first recall the definition of $L_\infty$ homomorphism
from [Definition 8.3.6, Fu].

\begin{defn}\label{Lmorphism}
Let $(C[1], \frak m)$, $(C'[1],\frak m')$ be $L_\infty$ algebras
and $\delta, \, \delta'$ be the associated coderivation. A
sequence $\varphi = \{\varphi_k\}_{k=1}^\infty$ with
$\varphi_k:E_kC[1] \to C'[1]$ is said to be an $L_\infty$
homomorphism if the corresponding coalgebra homomorphism $\widehat
\varphi: EC[1] \to EC'[1]$ satisfies
$$
\widehat \varphi \circ \delta
= \delta' \circ \widehat \varphi.
$$
We say that $\varphi$ is an {\it $L_\infty$ isomorphism}, if there
exists a sequence of  homomorphisms $\psi=
\{\psi_k\}_{k=1}^\infty$, $\psi: E_kC'[1] \to C'[1]$ such that its
associated coalgebra homomorphism $\widehat \psi: EC'[1] \to
EC[1]$ satisfies
$$
\widehat \psi \circ \widehat \varphi = id_{EC[1]},
\quad
\widehat \varphi \circ \widehat \psi = id_{EC'[1]}.
$$
In this case, we say that two $L_\infty$ algebras,
$(C[1], \frak m)$ and $(C'[1], \frak m')$ are $L_\infty$
isomorphic.
\end{defn}

We refer
to [Section 8.3, Fu] for more background materials on
the $L_\infty$ algebra and its homotopy theory, and also to [K1] in the
super language on a formal manifold.

\begin{thm} The two structures of strong homotopy Lie algebroid
on $T\FF \to Y$ induced by two  choices of splitting
$\Pi, \, \Pi^\prime$ are canonically $L_\infty$ isomorphic.
\end{thm}
\begin{proof}
We start with the expression of the symplectic form $\omega_U$
$$
\omega_U = \pi^*\omega - d\theta_G
$$
given in (\ref{eq:omega*}) that is canonically constructed on a
neighborhood $U$ of the zero section $E^* = T^*\FF$ when a
splitting $\Pi: TY = G\oplus T\FF$ is provided. To highlight
dependence on the splitting, we denote by $\theta_\Pi$ and
$\omega_\Pi$ the one form $\theta_G$ and the symplectic form
$\omega_U$. We will also denote by $\delta_\Pi$ the $\delta: EC[1]
\to EC[1]$ corresponding to the splitting $\Pi$.

Then for a  given splitting $\Pi_0$, we have
\begin{equation}\label{eq:difference}
\omega_\Pi - \omega_{\Pi_0} = d (\theta_{\Pi_0} - \theta_\Pi).
\end{equation}
In the super language, this is translated into
\begin{equation}\label{eq:super-difference}
H_{\Pi} - H_{\Pi_0}
= \{H_{\Pi_0}, \Gamma\}_\Omega
= - \{\Gamma, H_{\Pi_0} \}_\Omega
\end{equation}
where $\Gamma$ is the function associated to the one-form
$\theta_{\Pi _0}- \theta_\Pi$ which has $deg'(\Gamma) = 0$ (or
equivalently has $deg( \Gamma) = 1$). The last identity
comes from the super-commutativity of the bracket and
the fact that $deg(H_{\Pi_0}) = 2$ and $deg(\Gamma) =1$.
For the simplicity of
notations and also to make a connection with the more common
notation for the Gerstenhaber bracket as in Appendix, we will
simply write
$$
\{A, B\}_\Omega = [A,B]
$$
below.

Since any odd element commutes with itself under the bracket
$[\cdot, \cdot]$, we have $[\Gamma,\Gamma] = 0$ and so we have
$$
[[H_{\Pi_0}, \Gamma], \Gamma] = 0
$$
by the Jacobi identity. This then in turn implies
\begin{equation}\label{eq:nilpotent}
H_{\Pi_0} + [H_{\Pi_0}, \Gamma] = e^{ad_{(-\Gamma)}}(H_0)
\end{equation}
where $e^{ad_{(-\Gamma)}}$ is defined by
\begin{eqnarray*}
e^{ad_{(-\Gamma)}}(A) & := & \sum_{k=0} \frac{1}{k!}(ad_{(-\Gamma)})^k A \\
& = & A + [A, \Gamma] + \frac{1}{2!} [[A,\Gamma], \Gamma] +\cdots.
\end{eqnarray*}
Combining (\ref{eq:super-difference}) and (\ref{eq:nilpotent}), we have
obtained
$$
H_\Pi = e^{ad_{(-\Gamma)}}(H_{\Pi_0}).
$$
We recall the identity
$$
[e^{ad_{(-\Gamma)}}A, e^{ad_{(-\Gamma)}}B] = e^{ad_{(-\Gamma)}}[A, B]
$$
which can be rewritten as
\begin{equation}\label{eq:super-Linfty}
ad_{H_\Pi}(e^{ad_{(-\Gamma)}}A) = e^{ad_{(-\Gamma)}}(ad_{H_{\Pi_0}}A)
\end{equation}
when it is applied to $B = H_{\Pi_0}$.  Now we recall from
Corollary \ref{diffeo} that we have
$$
(\theta_{\Pi_0} - \theta_\Pi)|_{TE^*|_Y} \equiv 0.
$$
Noting that $\L = T\FF[1] \subset TE^*[1]|_Y$, this is translated
into
$$
\Gamma|_\L \equiv 0
$$
and so $ad_{\Gamma}$ naturally restricts to the functions on
$\L$. We denote by $\widetilde \varphi: EC[1] \to EC[1]$ the
coderivation associated to the restriction of $e^{ad_{(-\Gamma)}}$, then
(\ref{eq:super-Linfty}) is precisely translated into the identity
$$
\widetilde \varphi \circ \delta_{\Pi_0} = \delta_{\Pi}\circ
\widetilde\varphi.
$$
This proves that $\widetilde \varphi$ is an $L_\infty$
homomorphism. On the other hand, if $\Gamma'$ is the function
associated to the one-form $\theta_{\Pi} - \theta_{\Pi_0}$, then
we have $\Gamma' = -\Gamma$ and
$$
H_{\Pi_0} = H_\Pi + [H_\Pi, \Gamma']= e^{ad_{-\Gamma'}}(H_\Pi).
$$
Now we derive
\begin{eqnarray}\label{eq:Linverse}
e^{ad_{\Gamma'}}\circ e^{ad_{\Gamma}}(A) & = & e^{ad_{\Gamma}}(A)
+ [e^{ad_{\Gamma}}(A), \Gamma']
\nonumber \\
& = & (A + [A,\Gamma]) + [A + [A,\Gamma], \Gamma'] = A
\end{eqnarray}
for all $A$. Here for the last identity, we use the fact $\Gamma'
= -\Gamma$ and  the identity $[[A,\Gamma], \Gamma'] = 0$ . The
latter follows by the (super)-Jacobi identity using  the fact that
both $\Gamma$ and $\Gamma'$ are odd, and so $[\Gamma,\Gamma'] =
0$. If we denote by $\widetilde \varphi'$ the $L_\infty$
homomorphism associated to $\Gamma'$, then (\ref{eq:Linverse})
restricted to $\Omega^\bullet(\FF)$ as in (\ref{eq:delta'}) is
translated into the statement that $\widetilde \varphi'$ is the
inverse of $\widetilde \varphi$. This finishes the proof.
\end{proof}

This theorem then associates a canonical ($L_\infty$-)isomorphism
class of strong homotopy Lie algebras to each pre-symplectic
manifold and so to each coisotropic submanifold. It is obvious
from the construction that pre-Hamiltonian diffeomorphisms induce
canonical isomorphism by pull-backs in our strong homotopy Lie
algebroids. Pre-symplectic, in particular locally pre-Hamiltonian
diffeomorphisms also induce $L_\infty$ morphisms which however may
not be isomorphisms in general. For example, they do not induce
isomorphisms in $H^\bullet(\FF)$ in general, while global
pre-Hamiltonian diffeomorphisms do.

In the point of view of coisotropic embeddings this theorem implies
that our strong homotopy Lie algebroids for two
Hamiltonian isotopic coisotropic submanifolds are canonically
isomorphic and so the isomorphism class of the strong homotopy Lie
algebroids is an invariant of coisotropic submanifolds modulo the
Hamiltonian isotopy. We refer to the next section for the precise
explanation on the latter statement.

This enables us to study the moduli problem of deformations of
pre-symplectic structures on $Y$ in the similar way as done in
[K1], [FOOO], [Fu]. The followings are several interesting
questions to ask in this regard, which are analogs to Kontsevich's
formality theorem [K1] in our case.

\begin{ques}
\begin{enumerate}
\item Is the deformation problem formal in the sense of Kontsevich
[K1]? \item Does the $L_\infty$ structure on
$\Omega[1]^\bullet(Y,\omega)$ always canonically induce an
$L_\infty$ structure on its $\frak m_1$-cohomology
$H[1]^\bullet(Y,\omega)$
$$
\left(\bigoplus_\bullet H[1]^\bullet(Y,\omega), \{\overline{\frak
m}_\ell \}_{1 \leq \ell < \infty} \right)
$$
with $\overline{\frak m}_1 \equiv 0$. If not, what would be the
condition for this to be the case? \item If the answer is
affirmative in (2), are the two $L_\infty$ structures on
$\Omega[1]^\bullet(\FF)$ and its cohomology $H[1]^\bullet(\FF)$
quasi-isomorphic?
\end{enumerate}
\end{ques}

When $H^\bullet(\FF)$ is finite dimensional, the proof of [Theorem
8.3.5, Fu] can be imitated and so the answer is affirmative for
the questions (2) and (3) in that case. However in general
$H^\bullet(\FF)$ will be infinite dimensional. It would be very
interesting to see if the proof of [Theorem 8.3.5, Fu] can be
generalized to the case where $H^\bullet(\FF)$ is infinite
dimensional as in our case.

\section{Moduli problem and the Kuranishi map}
\label{sec:moduli}

In this section, we write down the defining equation
(\ref{eq:masterincoord}) for the graph $\operatorname{Graph}s
\subset TU \subset TE^*$ to be coisotropic in a formal
neighborhood, i.e., in terms of the power series of the section
$s$ with respect to the fiber coordinates in $U$.

In this section, we will study the moduli problem of the
Maurer-Cartan equation (\ref{eq:MC}) in the level of formal power
series. With respect to this strong homotopy Lie algebroid
constructed in section \ref{sec:algebroid}, the
formal power series version of (\ref{eq:masterincoord}) becomes
nothing but the Maurer-Cartan equation of $\frak
l^\infty_{(Y,\omega)}$. We refer to [Gz],
[GM], especially to [section 8.3, Fu] for more functorial formulation
of the formal moduli problem.

\begin{thm} The equation of the formal power series
solutions $\Gamma \in \frak l^1$
of $(\ref{eq:masterincoord})$
is given by
\begin{equation}\label{eq:MC}
\sum_{\ell = 1}^\infty \frac{1}{\ell !}
\frak m_\ell(\Gamma,\cdots, \Gamma) = 0 \quad \mbox{on } \,
\Omega^2(\FF)
\end{equation}
where
$$
\Gamma = \sum_{k = 1}^\infty \e^k \Gamma_k
$$
where $\Gamma_k$'s are sections of $T^*\FF$ and $\e$ is a formal
parameter.
\end{thm}
\begin{proof}
This immediately follows by substituting
$$
s = \Gamma = \sum_{k = 0}^\infty \e^k \Gamma_k
$$
into (\ref{eq:masterincoord}) and expanding the matrix
$(\widetilde\omega^{ij}) = (\omega_{ij} - s_\beta
F^\beta_{ij})^{-1}$ and comparing the result with the definition
of $\frak m_\ell$'s. Here we invoke the following matrix identity
$$
(A-B)^{-1} = A^{-1}(Id - BA^{-1})^{-1}
$$
for $A$ and $(A-B)$ invertible, and so we have
\begin{equation}\label{eq:omegaF}
(\omega - F \rfloor s)^{-1} = \omega^{-1}(Id - F^\# \rfloor s)^{-1}
= \omega^{-1} \sum_{\ell = 0}^\infty (F^\# \rfloor s)^\ell
\end{equation}
where we recall $F^\sharp: = F\omega^{-1}$ from (\ref{eq:Fsharp}).
Then the proof immediately follows from comparing (\ref{eq:MC}) and
the definition of $\frak m_k$ and
$\delta = \sum \widehat{\frak m}_\ell$ (\ref{eq:mell}) above.
\end{proof}

\begin{rem}\label{m0term} (\ref{eq:MC}) has the following
interpretation in terms of the deformation problem of
presymplectic structures on $(Y,\omega)$: Following the notation
from [Fu], [FOOO], we denote
$$
e^{\Gamma}:= \sum_{k = 0}\frac{1}{k!} \Gamma \otimes
\cdots \otimes \Gamma
$$
and write
$$
\frak m(e^{\Gamma}): =
\sum \frac{1}{\ell !}\frak m_\ell(\Gamma, \cdots, \Gamma).
$$
We define a new family of maps
$$
\frak m^\Gamma_k(\xi_1, \cdots, \xi_k) = \frak m(e^\Gamma, \xi_1,
\cdots, \xi_k)
= \sum_{\ell = 1}^\infty\frac{1}{\ell !}\frak m_{\ell+k} (\Gamma, \cdots,
\Gamma, \xi_1, \cdots, \xi_k)
$$
for $k \geq 1$, and
$$
\frak m^\Gamma_0 (1) : = \frak m(e^\Gamma) = \sum_{\ell = 1}^\infty
\frac{1}{\ell !}\frak m_\ell(\Gamma, \cdots, \Gamma).
$$
It was shown in [Fu] (or [FOOO] for the $A_\infty$
case) that the new coderivation
$$
\delta^\Gamma = \sum_{\ell = 0}^\infty \widehat{\frak m}^\Gamma_\ell
$$
satisfies the $L_\infty$-relation $\delta^\Gamma \delta^\Gamma =
0$ and so defines a weak $L_\infty$-algebra in general. By
definition, $\Gamma$ satisfying (\ref{eq:MC}) is equivalent to
$\frak m^\Gamma_0 = 0$. Therefore (\ref{eq:MC}) is precisely the
condition for this gauge changed (weak) $L_\infty$-structure to
define a strong $L_\infty$-structure, and is the Maurer-Cartan
equation for the deformation problem of the corresponding
presymplectic structure $(Y,\omega)$ as well.
\end{rem}

Now, we study (\ref{eq:MC}) inductively over the degrees of the
terms in the formal power series. Let and $\Gamma \in
\Omega^1(\FF)$ and
$$
\Gamma = \sum_{k=1}^\infty \e^k \Gamma_k.
$$
We fix a class $\alpha \in H^1(Y,\omega)$ and attempt to find
$\Gamma$ that satisfies (\ref{eq:MC}) in the formal power series and
that $d_\FF(\Gamma_1) =0$ and $[\Gamma_1] = \alpha \in H^1(Y,\omega)$.

For the sake of convenience, we will call $k$ the order of
the formal power series. Obviously the lowest order term of
(\ref{eq:MC}) is
$$
\frak m_1(\Gamma_1) = d_\FF(\Gamma_1) = 0,
$$
whose solution we assume is given in the class $\alpha \in H^1(Y,\omega)$.
Given $\Gamma_1$, the equation of the next order is
\begin{equation}\label{eq:weight2}
\frak m_1(\Gamma_2) + \frac{1}{2}\frak m_2(\Gamma_1, \Gamma_1) = 0.
\end{equation}
{}From the $L_\infty$-relation, we know that $\frak m_1$ is a derivation
with respect to $\frak m_2$. Therefore we have
$$
d_\FF(\frak m_2(\Gamma_1, \Gamma_1)) = 0 \quad \mbox{in } \,
\Omega^\bullet(\FF)
$$
since $d_\FF(\Gamma_1) = 0$. If we assume $[\frak
m_2(\Gamma_1,\Gamma_1)] = 0$ in $H^2(Y,\omega)$, then there exists
$\Gamma_2$ such that
$$
\frak m_1(\Gamma_2) + \frac{1}{2}\frak m_2(\Gamma_1,\Gamma_1)
= 0.
$$
We set $\Gamma^2 = \Gamma_1 + \Gamma_2$ which will then solve
(\ref{eq:MC}) up to the order of $2$. We can repeat this process
inductively over the degree $k$ to produce a solution
$$
\Gamma^k = \sum_{i=0}^k \e^i \Gamma_i
$$
up to the order of $k$ and then take the limit
$$
\Gamma = \lim_{k \to \infty} \Gamma^k = \sum_{\ell =1}^\infty
\e^\ell \Gamma_\ell
$$
provided the obstruction class vanishes in each step. We remark
that all the obstruction classes lie in the second cohomology
$H^2(Y,\omega)$. The limit exists as a formal power series
(or converges in the non-Archimedean topology
induced by the degrees of $\e^k$).  This proves the following
general theorem.

\begin{thm} Let $\FF$ be the null foliation of $(Y,\omega)$
and $\frak l = \oplus_{\ell =1}^{n-k} \frak l^{\ell}$ be the
associated complex. Suppose that $H^2(Y,\omega) = \{0\}$, i.e, any
$\FF$-closed two form is $\FF$-exact. Then for any given class
$\alpha \in H^1(Y,\omega)$, (\ref{eq:MC}) has a solution $\Gamma =
\sum_{k=1}^\infty \e^k \Gamma_k$ such that $d_\FF(\Gamma_1) =0$
and $[\Gamma_1] = \alpha \in H^1(Y,\omega)$. In other words, the
formal moduli problem is unobstructed.
\end{thm}

\begin{ques} Is it the case that whenever the formal moduli
problem is unobstructed, the corresponding $C^\infty$ moduli
problem is unobstructed? In other words, does the formal power
series obtained in the unobstructed case converge?
\end{ques}
We will investigate this question elsewhere.

\begin{defn} A pre-symplectic manifold $(Y,\omega)$ is called
{\it unobstructed} (resp. {\it formally unobstructed}) if the
corresponding moduli problem (resp. formal moduli problem) is
unobstructed. Otherwise it is called {\it obstructed} (resp. formally obstructed).
\end{defn}

\begin{cor}\label{hypersurface}
Let $(Y,\omega)$ have the nullity 1, i.e., $\operatorname{rank} E^* = 1$.
The moduli problem is unobstructed.
\end{cor}

In fact, it is easy to see that in this hypersurface case, the
genuine $C^\infty$ deformation problem is also unobstructed.
Therefore the first non-trivial case will be the one $(Y,\omega)$
with nullity 2.

We first provide a simple criterion for non-solvability of the
Maurer-Cartan equation. We denote by
$$
Z^\ell(Y,\omega):= Z^\ell(\FF) \subset \Omega^\ell(\FF)
$$
the set of $d_\FF$-closed $\ell$-forms (respectively by $B^\ell(Y,\omega)$
the set of $d_\FF$-exact one forms. Since $d_\FF = \frak m_1$ is
the derivation of $\frak m_2$, the bilinear map
$$
\frak l^1= \Omega^1(\FF) \to \frak l^2 = \Omega^2(\FF);\,
\Gamma_1 \to \frak m_2(\Gamma_1,\Gamma_1)
$$
canonically induces the map
$$
Kr: H^1(Y,\omega) \to H^2(Y,\omega); \, [\Gamma_1] \to
[\frak m_2(\Gamma_1,\Gamma_1)],
$$
where $[\cdot]$ is the cohomology class associated to the given
$d_\FF$-closed form. This is a version of the Kuranishi map for
this deformation problem which serves as the primary obstruction
to the deformation. The pairing
$$
H^1(Y,\omega) \otimes H^1(Y,\omega) \to H^2(Y,\omega); \quad
([\Gamma], [\Gamma^\prime]) \mapsto [\frak m_2(\Gamma,\Gamma')]
$$
is a special case of the so-called {\it Gerstenhaber bracket} [Ge].

\begin{thm}\label{kuranishi} Let $\alpha \in H^1(Y,\omega)$
such that $Kr(\alpha) \neq 0$ in $H^2(Y,\omega)$. Then
there is no solution
$$
\Gamma = \sum_{\ell = 1}^\infty \e^\ell \Gamma_\ell
$$
for (\ref{eq:MC}) with $[\Gamma_1] = \alpha$. In particular,
there is no smooth isotopy $Y_t \subset (X,\omega)$ of
coisotropic embeddings with
$$
Y_0 = Y, \quad \frac{d}{dt}\Big|_{t=0} Y_t = \Gamma_1.
$$
\end{thm}

One particular case is worth of mentioning

\begin{cor}\label{dgla}
Suppose that all $\frak m_k = 0$ for $k \geq 3$. Then $\alpha \in
H^1(Y,\omega)$ is formally unobstructed if and only if $Kr(\alpha)
= 0$.
\end{cor}

\begin{exm} We will analyze the example studied by M. Zambon
[Za] in the light of Theorem \ref{kuranishi}. This example
illustrates that the space of $C^1$-close coisotropic submanifolds
modulo the Hamiltonian isotopy is {\it not} smooth at the
Hamiltonian isotopy class of the torus. Let $(Y,\omega)$ be the
standard 4-torus $T^4 = \R^4/\Z^4$ with coordinates $(y^1, y^2,
q^1, q^2)$ with the closed two form
$$
\omega_Y = dy^1 \wedge dy^2.
$$
Note that the null foliation is provided by the 2-tori
$$
\{y^1 = const,\, y^2 = const \},
$$
and it also carries the transverse foliation given by
$$
\{q^1 = const, \, q^2 = const \}.
$$
The canonical symplectic thickening is given by
\begin{eqnarray*}
E^* & = & T^4 \times \R^2 = T^2 \times T^*(T^2), \\
\omega & = & dy^1\wedge dy^2 + (dq^1 \wedge dp^1 + dq^2 \wedge dp^2),
\end{eqnarray*}
where $p^1, \, p^2$ are the canonical conjugate coordinates of
$q^1, \, q^2$.
It follows that the transverse curvature $F \equiv 0$ and so
all $\frak m_\ell = 0$ for $\ell \geq 3$ and
the Maurer-Cartan equation (\ref{eq:MC}) becomes the quadratic equation
\begin{equation}\label{eq:tori}
d_\FF(\Gamma) + \frac{1}{2}\{\Gamma, \Gamma\} = 0
\end{equation}
where $\{\cdot,\cdot\} = \frak m_2$ given by the formula in
(\ref{eq:torim2}) below. In particular we have
\begin{equation}\label{eq:vanishing}
\int_{T^2} \{\Gamma,\Gamma\} = 0
\end{equation}
for any solution $\Gamma$ of (\ref{eq:tori}).
It is easy to compute
\begin{eqnarray*}
H^0(Y,\omega) & \cong & C^\infty(T^2), \\
H^1(Y, \omega) & \cong & C^\infty(T^2)\{\theta^1,\theta^2\}, \\
H^2(Y, \omega) & \cong & C^\infty(T^2)\{\theta^1 \wedge
\theta^2\},
\end{eqnarray*}
where $\theta^i = [dq^i] \in H^1(Y,\omega)$.
We consider a one form
$\Gamma = a_1(y,q) dq^1 + a_2(y,q) dq^2$ that is $d_\FF$-closed, i.e.,
satisfies
$$
\frac{\del a_2}{\del q^1} - \frac{\del a_1}{\del q^2} = 0.
$$
The map $Kr: H^1(Y,\omega) \to H^2(Y,\omega)$ is induced by the
bilinear map
$$
\Gamma  \to \frak m_2(\Gamma,\Gamma),
$$
which can be written in coordinates as
\begin{equation}\label{eq:torim2}
a_1(y,q) dq^1 + a_2(y,q) dq^2 \to \{a_1, a_2\}_y dq^1 \wedge dq^2,
\end{equation}
where the bracket is defined by
\begin{equation}
\{a_1, a_2\}_y = \Big((\frac{\del a_1}{\del y^2} \frac{\del
a_2}{\del y^1}) - (\frac{\del a_2}{\del y^2} \frac{\del a_1}{\del
y^1}) \Big).
\end{equation}
Therefore any infinitesimal deformation $\Gamma$ with the
non-vanishing integration over the fiber  $\Sigma_{(y^1,y^2)} =
T^2$
$$
\int_{\Sigma_{(y^1,y^2)}}
\{a_1, a_2\}_y dq^1 \wedge dq^2 \neq 0
$$
will be obstructed. For example, one can take
$$
\Gamma = \sin (2\pi y^1) dq^1 + \sin (2 \pi y^2) dq^2,
$$
which is the example Zambon looked at. In this case, we have
$$
\int \{a_1, a_2\}_y dq^1\wedge dq^2 = - 4\pi^2 \cos 2\pi y^2 \cos
2\pi y^1 \not \equiv 0.
$$
\end{exm}

In the next section, we study a more non-trivial example where the
structure of null foliations becomes more complicated.

\section{An example}
\label{sec:example}

In this section, we analyze one parameter family of examples
$(Y_\alpha,\omega_\alpha)$ from the mechanics of harmonic
oscillator. {}From the analysis of this example, it is manifest
that the deformation problem of coisotropic submanifolds is
closely tied to the geometry and dynamics of the null foliation. A
systematic study of geometry of coisotropic submanifolds in terms
of geometry of the foliation theory will be carried out elsewhere.

Consider the harmonic oscillator
$$
H = \frac{1}{2}(|Q|^2 + |P|^2)
$$
on the phase space $T^*\R^3 \cong \R^6$ where $Q = (Q^1,Q^2,Q^3)$
and $P=(P^1,P^2,P^3)$ are the position and momentum coordinates.
We can write
$$
H = H_1 + H_2 + H_3,
$$
where $H_i$ are the one dimensional harmonic oscillator Hamiltonians
$$
H_i = \frac{1}{2}((Q^i)^2 + (P^i)^2), \quad i = 1,\, 2, \, 3.
$$
It follows that  $\{H, H_i\} = \{H_i, H_j\} = 0$ where $\{\cdot,
\cdot\}$ is the canonical Poisson bracket. We fix two constants
$\alpha, \, \beta > 0$ such that
$$
\alpha > 1, \quad 0 < \beta < \frac{1}{2},
$$
and consider the submanifold $Y_{\alpha,\beta}
\subset S^5 \subset \R^6$ defined by
$$
Y_{\alpha,\beta} = \{ (q,p) \in \R^6 \mid
H(q,p) = \frac{1}{2}, \, (H_1 + \alpha H_2)(q,p) = \beta\}.
$$
It is easy to see that these provide two parameter family of
smooth coisotropic submanifolds of $\R^6$ whose images are all
contained in the unit sphere $H^{-1}(1/2)$.
To simplify the discussion, we fix $\beta = \frac{1}{4}$ and
denote
$$
Y_\alpha := Y_{\alpha,1/4}, \quad H_\alpha := H_1 + \alpha H_2.
$$
It is straightforward to check that the Hamiltonian vector fields
$X_H$ and $X_{H_\alpha}$ are linearly independent everywhere
on $Y_\alpha$ and so the characteristic distribution is
given by
$$
E = \operatorname{span}_\R\{X_H, X_{H_\alpha} \}
$$
which is a trivial bundle. In particular all leaves are
orientable. We recall
$$
Y_\alpha =\{(Q,P) \in \R^6 \mid |Q|^2 + |P|^2 = 1, \,
H_1 + \alpha H_2 = \frac{1}{4} \}.
$$
On $Y_\alpha$, we derive
\begin{eqnarray}\label{eq:L1L2L3}
H_1 & = & \frac{1}{\alpha -1} \left(\frac{2\alpha -1}{4} - \alpha H_3\right), \\
H_2 & = & \frac{1}{\alpha -1} \left(H_3 - \frac{1}{4}\right).
\end{eqnarray}
Since $H_1, \, H_2 \geq 0$, we have obtained the bound
\begin{equation}\label{eq:bound}
\frac{1}{4} \leq H_3 \leq \frac{2\alpha -1}{4\alpha}.
\end{equation}

We will denote by $\ell_i$ the value of the corresponding  $H_i$.
For $ \frac{1}{4} < \ell_3 < \frac{2\alpha -1}{4\alpha}$, $H_3$
are regular on $Y_\alpha$ and $H_3^{-1}(\ell_3)$ is a 3-torus
$$
H_3^{-1}(\ell_3) = S^1(\sqrt{2 \ell_1}) \times S^1(\sqrt{2\ell_2})
\times S^1(\sqrt{2\ell_3}),
$$
where $\ell_i, \, i =1, \,2$ are given by the formula
(\ref{eq:L1L2L3}) for a given $\ell_3$. It follows that leaves of
the null foliation on the open subset $H_3^{-1}(\frac{1}{4},
\frac{2\alpha -1} {4\alpha}) \subset Y_\alpha$ are generated by
$$
\{(1,1,1), (1, \alpha, 0) \}
$$
in the rectangularpid $[0,2\pi] \times [0,2\pi] \times [0,2\pi]$.
Note that this family of lines is uniquely determined by $\alpha$
and $\ell_3$.

When $\alpha$ is rational, then all the leaves are compact. When
$\alpha$ is irrational, all the leaves are non-compact. They are
all immersions of $\R \times S^1$ are dense in the three torus
$H_3^{-1}(\ell_3), \, \ell_3 \neq \frac{1}{4}, \, \frac{2\alpha-1}{4\alpha}$
which are invariant under the Hopf action
of $S^1$ on $Y_\alpha \subset S^5 \subset \C^3$.

When $\ell_3 = \frac{1}{4}$, we have $\ell_2 = 0$ and $\ell_1 = \frac{1}{4}$.
Hence $\{(Q,P) \in Y_\alpha \mid \ell_3 = \frac{1}{4}\}$ consists of
the unique leaf which is nothing but
$$
S^1(\frac{1}{\sqrt{2}}) \times \{(0,0)\} \times S^1(\frac{1}{\sqrt{2}}).
$$
Similarly when $\ell_3 = \frac{2\alpha -1}{4\alpha}$, we have
$\ell_1 =0 $ and $\ell_2 = \frac{1}{4\alpha}$ and hence
$\{(Q,P) \in Y_\alpha \mid \ell_3 = \frac{2\alpha -1}{4\alpha} \}$
is the torus
$$
\{(0,0)\} \times S^1(\frac{1}{\sqrt{4\alpha}}) \times S^1(\sqrt
{\frac{2\alpha-1}{4\alpha}}).
$$

Now we compute the transverse curvature of the leaves. For this
we consider the canonical splitting provided by the complex
structure on $\R^{6} \cong \C^3$. In other words, we choose
the splitting $TY = G \oplus E$ with $G$ given by
$$
G = (E \oplus JE)^\perp,
$$
where $\perp$ is the Euclidean orthogonal complement and $J$
denote the standard complex structure on $\C^3$. We choose
coordinates $(y^1, y^2, q^1, q^2)$ of $Y_\alpha$
by the functions
\begin{eqnarray}
y^1 & = & \alpha \theta_1 - \theta_2 - \frac{1+\alpha^2}{\alpha}
\theta_3,\nonumber\\
y^2 & = & H_3, \nonumber \\
q^1 & = &\theta_3, \nonumber \\
q^2 &= &\theta_1 + \alpha \theta_2 +
(1+\alpha^2)\theta_3.\nonumber
\end{eqnarray}
Here $(r_1,r_2,r_3,\theta_1,\theta_2,\theta_3)$ is the polar
coordinates of $\R^6$ and we have chosen $y^1, \, y^2$ so that
the leaves of the null foliation are given by
$$
\{y^1 =const, \, y^2 = const \}.
$$
A straightforward calculation leads to
\begin{equation}\label{eq:L's}
H_1  =  \frac{1}{\alpha -1}(\frac{2\alpha -1}{4} - \alpha y^2),
\qquad H_2 =
\frac{1}{\alpha -1}(y^2 - \frac{1}{4}),
\qquad H_3 = y^2.
\end{equation}
We derive
\begin{eqnarray}
\theta_1 & = & q^1 + \frac{q^2 + \alpha y^1}{1 + \alpha^2}, \nonumber \\
\theta_2 &= & q^1 + \frac{\alpha q^2 - y^1}{1 +\alpha^2}, \nonumber \\
\theta_3 & =& q^1.\nonumber
\end{eqnarray}
Therefore $Y_\alpha$ is parameterized by
\begin{equation}\label{eq:parameterization}
\Big(r_1 e^{i(q^1 + \frac{q^2 + \alpha y^1}{1+\alpha^2})}, r_2
e^{i(q^1 + \frac{\alpha q^2 - y^1}{1+\alpha^2})}, r_3 e^{i q^1}
\Big)
\end{equation}
in this coordinates,
where we note $r_i = \sqrt{2L_i}$.  From this, we derive
\begin{eqnarray}\label{eq:ddy's}
\frac{\del}{\del y^1} & = & \frac{\alpha}{1 + \alpha^2}
\frac{\del}{\del \theta_1} - \frac{1} {1 +
\alpha^2}\frac{\del}{\del\theta_2} = \frac{1}{1 +\alpha^2}X_{\alpha H_1-H_2},
\nonumber \\
\frac{\del}{\del y^2} & = & - \frac{\alpha}{2(\alpha
-1)}\frac{1}{r_1}\frac{\del}{\del r_1}+
\frac{1}{2(\alpha-1)}\frac{1}{r_2}\frac{\del}{\del r_2}  +
\frac{1}{2r_3} \frac{\del}{\del r_3}, \nonumber \\
\frac{\del}{\del q^1} & = &
\frac{\del}{\del \theta_1} +  \frac{\del}{\del \theta_2} +
\frac{\del}{\del \theta_3} = X_H,
\nonumber\\
\frac{\del}{\del q^2} & = & \frac{1}{1 + \alpha^2}
\frac{\del}{\del \theta_1} + \frac{\alpha} {1 +
\alpha^2}\frac{\del}{\del\theta_2} = \frac{1}{1 + \alpha^2}X_{H_\alpha}.
\end{eqnarray}
We consider the orthogonal splitting $TY = G \oplus E$ with
respect the induced metric on $Y_\alpha$ from the Euclidean metric
on $\R^6$. A calculation shows
\begin{eqnarray}
G & = & \operatorname{span}_\R \Big\{
-\frac{\alpha}{r_1^2}\frac{\del}{\del \theta_1} + \frac{1}{r_2^2}
\frac{\del}{\del \theta_2},
\frac{\del}{\del y^2} \Big\} \nonumber \\
& = & \operatorname{span}_\R \Big\{
-\frac{\alpha}{H_1}\frac{\del}{\del \theta_1} + \frac{1}{H_2}
\frac{\del}{\del \theta_2}, \frac{\del}{\del y^2} \Big\}.
\nonumber
\end{eqnarray}
A straightforward but lengthy calculation gives
$$
G = \operatorname{span}_\R \{e_1, e_2 \},
$$
where
\begin{eqnarray}
e_1 & = &\frac{\del}{\del y^1}
- \frac{\alpha(H_1-H_2)}{\alpha^2
H_2 +H_1} \frac{\del}{\del q^2},\nonumber \\
e_2 & = & \frac{\del}{\del y^2}.\nonumber
\end{eqnarray}
Hence we have the matrix $(R_i^\alpha)$ given by
\begin{equation}
R^1_2  =  0,\qquad
R^2_2 = - \frac{\alpha(H_1-H_2)}{\alpha^2 H_2 + H_1},\qquad
R^1_1  = R^2_1= 0.
\end{equation}
Noting that $H_i$'s are functions of $y^2$ alone in
(\ref{eq:L1L2L3}), it follows from \ref{eq:curvature}
that $F^\alpha_{ij} \equiv 0$, which
can be also seen from the observation that
the chosen splitting $G$ is also integrable. This proves
$$
\frak m_k = 0 \quad \mbox{for all } \, k\geq 3,
$$
and the corresponding homotopy algebra reduces a differential
graded Lie algebra. A straightforward calculation also shows
\begin{equation}\label{eq:formonY}
\omega|_{Y_\alpha} = \frac{1}{2(\alpha -1)}
dy^1 \wedge dy^2
\end{equation}
on $Y_\alpha \setminus H_3^{-1}(\frac{1}{4})
\cup H_3^{-1}(\frac{2\alpha-1}{4\alpha})$.

Finally we study the Kuranishi map
$$
Kr: H^1(Y_\alpha, \omega_\alpha) \to H^2(Y_\alpha,\omega_\alpha);\quad
[\Gamma] \to [\frak m_2(\Gamma, \Gamma)],
$$
where $\Gamma$ is $d_\FF$-closed one form.  We first need to characterize
the set of $d_\FF$-closed one forms.
Recalling that $E = \operatorname{span}_\R\{X_H, X_{H_\alpha}\}$,
we represent $E^* = \{f_1^*, f_2^* \}$ where $\{f_1^*, f_2^*\}$ is
the dual frame of $\{X_H,X_{H_\alpha}\}$. Let
$$
\Gamma = A f_1^* + B f_2^*
$$
be a one form in $\Omega^1(\FF)$ where $A, \, B$ are globally
defined function on $Y_\alpha$. From the definition of $d_\FF$, we
have
$$
d_\FF(\Gamma) = (X_H(B) - X_{H_\alpha}(A))f_1^* \wedge f_2^*.
$$
Therefore $\Gamma$ is $d_\FF$-closed if and only if
\begin{equation}\label{eq:closed}
X_H(B) = X_{H_\alpha}(A).
\end{equation}
On the other hand, for any smooth function $C \in
C^\infty(Y_\alpha)$, we have
$$
d_\FF(C) = X_H[C] f_1^* + X_{H_\alpha}[C] f_2^*.
$$
We will compute $H^1(Y,\omega)$ and $H^2(Y,\omega)$ and study the
Kuranishi map
$$
Kr: H^1(Y,\omega) \to H^2(Y,\omega).
$$

\subsection{Symplectic reduction and integration over fibers}

We first apply the symplectic reduction with respect to the Hopf
action of $S^1$ which induces the following commutative diagram
\begin{equation}
\begin{matrix}
Y_\alpha \,\,\,\,\, & \hookrightarrow  && S^5 \subset \C^3\setminus \{0\} \\
\downarrow \pi_1 & \quad && \downarrow \pi_1 \\
\widetilde{Y}_\alpha = Y_\alpha/S^1 & \hookrightarrow && \C P^2
\end{matrix}\nonumber
\end{equation}

We will denote by $\widetilde\omega$ the reduced symplectic form
on $H^{-1}(\frac{1}{2})/S^1 = S^5/S^1 \cong \C P^2$ which is
nothing but the standard Fubini-Study form, and by
$\widetilde\omega_\alpha$ the induced pre-symplectic form on
$\widetilde Y_\alpha \subset \C P^2$. Since $H_\alpha$ is
invariant under the Hopf action, it projects down to a function
$\widetilde H_\alpha$ on $\widetilde Y_\alpha$.

$\widetilde Y_\alpha \subset \C P^2$ is a hypersurface for which
the induced null foliation $\widetilde\FF$ on $\widetilde
Y_\alpha$ are given by $\Sigma/S^1$ for each $\Sigma$ is a leaf of
$\FF$. In fact, we have the obvious one-one correspondence
$$
Y_\alpha/\sim \, \cong \widetilde Y_\alpha/\sim.
$$
Furthermore $d_{\widetilde \FF}$ is given by
$$
d_{\widetilde \FF}(f) = X_{\widetilde H_\alpha}[f] \widetilde
f_2^*
$$
for $f \in C^\infty(\widetilde Y_\alpha)$ where $\widetilde f_2^*$ is the
push-forward of $f_2^*$ to $\widetilde Y_\alpha$.
We denote by
$$
(\pi_1)_*: \Omega^1(Y_\alpha,\omega_\alpha) \to
C^\infty(\widetilde Y_\alpha) = \Omega^0(\widetilde
Y_\alpha,\widetilde\omega_\alpha)
$$
the integration over fibers defined by
\begin{equation}\label{eq:pi*1}
(\pi_1)_*(A f_1^* + B f_2^*)(y) = \int_{\pi_1^{-1}(y)} A f_1^*.
\end{equation}
It is easy to check the identity
$$
d_{\widetilde\FF} \circ (\pi_1)_* = (\pi_1)_*\circ d_\FF
$$
and so $(\pi_1)_*$ induces a natural homomorphism
\begin{equation}
(\pi_1)_*: H^1(Y_\alpha,\omega_\alpha) \to C^\infty(\widetilde
Y_\alpha)\nonumber
\end{equation}
Similarly we define
$$
(\pi_1)_*: \Omega^2(Y_\alpha,\omega_\alpha) \to
\Omega^1(\widetilde Y_\alpha, \widetilde \omega_\alpha)
$$
by
\begin{equation}\label{eq:pi*2}
(\pi_1)_*(D\cdot f_1^* \wedge f_2^*) = \Big( \int_{\pi^{-1}(y)} D
f_1^*\Big) \widetilde f_2^*,
\end{equation}
which again intertwines $d_\FF$ and $d_{\widetilde \FF}$ and so
induces the homomorphism
\begin{equation}
(\pi_1)_*: H^2(Y_\alpha,\omega_\alpha) \to H^1(\widetilde
Y_\alpha,\widetilde\omega_\alpha).\nonumber
\end{equation}

Since the spaces of leaves $Y_\alpha/\sim$ and $\widetilde
Y_\alpha/\sim$ are isomorphic, it is easy to check that
$$
(\pi_1)_*: H^2(Y_\alpha,\omega_\alpha) \to H^1(\widetilde
Y_\alpha,\widetilde\omega_\alpha)
$$
is an isomorphism. We next quote the following general theorem by
Haefliger [Ha].
Here $\Omega_c^\bullet(Tr\FF)$ is the topological differential
graded vector space of forms on $T/H$, where $H$ is the holonomy
pseudogroup induced on a complete transversal submanifold $T$.
Since we will not use the theorem
except the presence of the isomorphism (\ref{eq:Omega0FF}),
we refer readers to [Ha] for more detailed explanation.

\begin{thm}\label{Haefliger} {\bf [Theorem 3.1, Ha]}
Let $\FF$ be a foliation on $X$ with leaves
of dimension $p$, and assume that the tangent bundle to the leaves
is oriented. Then there is a continuous open surjective linear map
$$
\int_\FF: \Omega^{p+k}_c (X) \to \Omega_c^k(\operatorname{Tr} \FF)
$$
which commutes with $d_\FF$. And it induces an isomorphism
\begin{equation}\label{eq:Omega0FF}
\int_\FF: H^p(\FF) \to \Omega^0(\operatorname{Tr}\FF).
\end{equation}
\end{thm}

Since we have the isomorphism
\begin{equation}\label{eq:0isom}
\Omega^0(\operatorname{Tr }\FF) \cong \Omega^0(\operatorname{Tr}
\widetilde \FF)
\end{equation}
which obviously follows from the fact that the spaces of leaves
$Y_\alpha/\sim$ and $\widetilde Y_\alpha/\sim$ are isomorphic,
$$
\pi_*: H^2(Y_\alpha,\omega_\alpha) \to H^1(\widetilde
Y_\alpha,\widetilde\omega_\alpha)
$$
is an isomorphism.  It is also known [Ha] that $\Omega^0(Tr \FF)$ is
naturally isomorphic to $C^\infty(Y/\sim)$ when the leaf space
$Y/\sim$ is an orbifold. One can also prove these statements directly without
referring to the above result from [Ha]. We did this to illustrate
a connection of the current deformation problem to the foliation
theory. A more systematic study in this regard will be carried out
elsewhere in the future.

From now on we treat two cases, $\alpha$ rational and irrational,
separately. We will postpone the study of irrational case to a
future work.

\subsection{The rational case}

We first note that the foliation $\FF$ is generated by the closed
one forms
\begin{eqnarray*}
dy^1 & = & \alpha d\theta_1 - d\theta_2 - \frac{1+\alpha_2}{\alpha}d\theta_3\\
dy^2 & = & dH_3.
\end{eqnarray*}
Both are globally well-defined and the period group is given by
$$
P := \operatorname{span}_\Z\{dy^1(H_1(Y_\alpha,\Z)),
dy^2(H_1(Y_\alpha,\Z))\} = \operatorname{span}_\Z \Big\{\alpha, -
\frac{1+\alpha^2}{\alpha}\Big\}.
$$
Suppose that $\alpha$ is rational and $\alpha = \frac{p}{q}$ with
relatively prime integers $p, \, q$. The period group $P$ of $\FF$
is generated by
$$
\Big\{ \frac{p}{q}, - \Big(\frac{p}{q} + \frac{q}{p}\Big) \Big\}
$$
and so
$$
P = \operatorname{span}_\Z\Big\{\frac{\operatorname{gcd}\{p^2,
-(p^2+q^2)\}}{qp}\Big\}.
$$
If we denote
$$
\frac{L}{M} = \frac{\operatorname{gcd}\{p^2, -pq, -q^2\}}{qp}
$$
with $L, \, M$ relatively prime to each other, then the local
holonomy group is trivial except the two leaves
$H_3^{-1}(\frac{1}{4})$ and $H_3^{-1}(\frac{2\alpha - 1}{4\alpha})$ at which the
holonomy groups are finite cyclic groups.

It turns out that the leaf space $Y_\alpha/\sim$ is a compact
Hausdorff symplectic orbifold induced by the canonical
transverse symplectic form. There are two orbifold points
mentioned above. We denote them by $z^-, \, z^+ \in Y_\alpha/\sim$
respectively. See [Sa] for the definition of orbifolds (or V-manifolds)
and [We3] for the definition of symplectic orbifolds (or symplectic V-manifolds).

Now we describe the leaf space more precisely. It remains to study
local structure near the two orbifold points $z_\pm$. Around  the
leaf $H_3^{-1}(\frac{1}{4})$,  we choose coordinates
$$
(H_2,\theta_2,q^1,q^2)
$$
where $q^1, \, q^2$ are the ones as before. It follows from
(\ref{eq:L1L2L3})
\begin{eqnarray}\label{eq:y1y2}
y^1 & = & \alpha q^2 + (1+\alpha^2)q^1 - (1+\alpha^2)\theta_2 \\
y^2 & = & (\alpha -1) H_2 + \frac{3}{4}
\end{eqnarray}
and so
$$
\frac{\del}{\del H_2} = (\alpha -1) \frac{\del}{\del y^2}, \,
\frac{\del}{\del \theta_2} = - (1 + \alpha^2) \frac{\del}{\del
y^1}.
$$
Substituting this into (\ref{eq:formonY}), we derive
$$
\omega|_{Y_\alpha} = (1 + \alpha^2)dH_2 \wedge d\theta_2 = (1 +
\alpha^2)r_2dr_2 \wedge d\theta_2.
$$
Recall that $(y^1, y^2)$ parameterizes the leaf space around $z_1$
and two leaves parameterized by $(y^1, y^2)$ and $(Y^1, Y^2)$
coincide if and only if
$$
y^2 = Y^2, \quad \mbox{and }\,\,  y^1\equiv Y^1\, \mod 2\pi.
$$
In the coordinates $(H_2,\theta_2)$, the latter condition is
equivalent to
$$
-(1+\alpha^2)\theta_2 \equiv -(1+\alpha^2)\Theta_2\,  \mod 2 \pi
\Longleftrightarrow \theta_2 \equiv \Theta_2\, \mod \frac{2\pi}{1
+ \alpha^2}.
$$
Therefore the local chart of $Y_\alpha/\sim$ is given by the
non-flat cone
$$
\R^2/\Z_{M_-}
$$
if we denote
$$
\frac{1}{1 + \alpha^2} = \frac{L_-}{M_-}
$$ with
$L_-, \, M_-$ relatively prime. Note that the symplectic form is
certainly invariant under the linear action by $\Z_{M_-}$. Similar
computation can be carried out at the leaf
$z^+=H_3^{-1}(\frac{2\alpha-1}{4\alpha})$ using the coordinates
$$
(H_1,\theta_1, q^1,q^2)
$$
and have the symplectic form
$$
\omega = \frac{1+\alpha^2}{2\alpha^2}dH_1 \wedge d\theta_1 =
\frac{1+\alpha^2}{\alpha^2}r_1dr_1 \wedge d\theta_1.
$$
This time the holonomy group is $Z_{M_+}$ where we write
$$
\frac{\alpha}{1+\alpha^2} = \frac{L_+}{M_+}
$$
for relatively prime integers $L_+, \, M_+$. It is easy to see
that the $Y_\alpha/\sim$ is topologically sphere. We denote the
corresponding symplectic orbifold by
$$
(S^2_{\alpha}, \omega_{\alpha}).
$$
In fact in the rational case, the flows of $H$ and $H_\alpha$
generate two circle actions which commute each other and so
defines a torus action. The leaves of the null foliation are then
just the orbits of this torus action and defines a fibration over
$S^2_\alpha$. See (\ref{eq:parameterization}). We denote by
$C^\infty(S^2_\alpha)$ the set of smooth functions on the orbifold
$S^2_\alpha$. Then we have

\begin{prop}\label{2ndcohom}{\bf [2nd Cohomology]} When $\alpha >1$ is rational,
we have
$$
H^2(Y_\alpha,\omega_\alpha) \cong C^\infty(S^2_{\alpha})\{ a^1
\wedge a^2\}.
$$
\end{prop}
\begin{proof} This follows from (\ref{eq:Omega0FF})-(\ref{eq:0isom})
and from the remark
right after Theorem \ref{Haefliger}.
\end{proof}

We next analyze $H^1(Y_\alpha,\omega_\alpha)$.
We consider an open covering of $S^2_\alpha = D^+ \cup D^-$
where $D^\pm$ is a neighborhood of the two orbifold
points $p^\pm$ such that $D^+_\alpha \cap D^-_\alpha$
is diffeomorphic to
$S^1 \times (-\e, \e)$. We denote
$$
Y^\pm = \pi^{-1}(D^\pm_\alpha)
$$
which provides a covering $Y_\alpha = Y^+ \cup Y^-$.
We emphasize that both $D^\pm_\alpha$ are orbifolds
diffeomorphic to $D^2/\Z_{M^\pm}$ respectively.
We then consider the Mayer-Vietoris sequence
\begin{eqnarray*}
&\to & H^0(Y^+ \cap Y^-,\omega_\alpha) \stackrel{\delta^0}\to
H^1(Y_\alpha,\omega_\alpha)
\stackrel{\alpha^1}\to \\
& \stackrel{\alpha^1} \to & H^1(Y^+,\omega_\alpha)\oplus
H^1(Y^-,\omega_\alpha)\stackrel{\beta^1} \to
H^1(Y^+ \cap Y^-, \omega_\alpha) \stackrel{\delta^1}
\to H^2(Y_\alpha,\omega_\alpha) \to
\end{eqnarray*}
From this exact sequence, we have derived
\begin{equation}\label{eq:H1}
H^1(Y_\alpha,\omega_\alpha) \cong \mbox{im } \delta^0
\oplus \ker \beta^1
\end{equation}
We now describe the two summands more explicitly.
We note that the forms
$$
f_1^*, \, f_2^*
$$
are leafwise closed, and
invariant under the torus action and also under the holonomy.
We denote the corresponding cohomology class thereof by
$a_j = [f_j^*] \in H^1(Y_\alpha, \omega_\alpha)$.
Then the above discussion gives rise to
\begin{eqnarray*}
H^1(Y^+\cap Y^-,\omega_\alpha) & \cong & C^\infty(S^1)\{a^1, a^2\} \\
H^1(Y^+,\omega_\alpha) & \cong & C^\infty(D^+_\alpha)\{a^1, a^2\} \\
H^1(Y^-,\omega_\alpha) & \cong & C^\infty(D^-_\alpha)\{a^1, a^2\}\\
\end{eqnarray*}
Furthermore again from the above exact sequence, we have
\begin{equation}\label{eq:kerbeta1}
\ker \beta^1 \cong H^1(Y^+,\omega_\alpha) \times_{H^1(Y^+\cap Y^-,
\omega_\alpha)} H^1(Y^-,\omega_\alpha)
\cong C^\infty(S^2_\alpha)\{a^1,a^2\}
\end{equation}
where the middle term is the obvious fiber product.
This finishes the description of the first cohomology.
We summarize our discussion into

\begin{prop}\label{1stcohom}{\bf [1st Cohomology]}
Let $\alpha >1$ be rational.
Then we have the isomorphism
\begin{eqnarray*}
H^1(Y_\alpha,\omega_\alpha) & \cong & \ker \beta^1 \oplus
\mbox{\rm im } \delta^0 \\
& \cong & C^\infty(S^2_\alpha)\{a^1,a^2\}
\oplus \mbox{\rm im } \delta^0.
\end{eqnarray*}
\end{prop}

Next we study the Gerstenhaber bracket $[\cdot, \cdot]$
and the Kuranishi map

\begin{prop}\label{Kura}{\bf [Gerstenhaber bracket and Kuranishi map]}
Under the isomorphisms in Theorem \ref{2ndcohom} and Theorem
\ref{1stcohom}, the Gerstenhaber bracket is given by
\begin{equation}\label{eq:gammaeta}
[\gamma, \eta] = 0
\end{equation}
for all $\gamma \in H^1(Y_\alpha, \omega_\alpha)$ and
$\eta \in \mbox{\rm im } \delta^0$
and
\begin{equation}\label{eq:gamma12}
[\gamma_1, \gamma_2] = (\{g_1, h_2\}_\alpha\circ \pi)
a^1 \wedge a^2
\end{equation}
where $\Gamma_j = g_j a^1 + h_j a^2, \quad j = 1, \, 2$
and $g_j, \, h_j \in C^\infty(S^2_\alpha)$, and
$\{\cdot, \cdot \}_\alpha$ is the Poisson bracket
on $S^2_\alpha$. In particular,
the Kuranishi map $Kr$ is given by the formula
$$
Kr(\Gamma) = (\{g, h\}_{\alpha} \circ \pi) a^1 \wedge a^2
$$
for $\Gamma = ga^1 + ha^2$.
\end{prop}
\begin{proof}
We first note that $\ker \alpha^1 = \mbox{\rm im }\delta^0$: any
element $a$ lying in $\mbox{\rm im }\delta^0 \subset
H^1(Y,\omega_\alpha)$ is represented by a leafwise one form
$\alpha$ determined by
$$
\alpha = \begin{cases}
dg_+ \quad \mbox{on }\, Y^+ \\
dg_- \quad \mbox{on }\, Y^-
\end{cases}
$$
such that $i^*_+g_+ - i^*_-g_- = f$ for some leafwise
constant function $f$ or equivalently for $f$ satisfying
$d_\FF(f) = 0$ on $Y^+ \cap Y^-$. Using this representation
of an element from $\mbox{im }\delta^0$ and Theorem 9.3,
which states that
$d_\FF$ is a derivation of $\{\cdot, \cdot\}_\Pi$ because
$F_\Pi\equiv 0$, it is easy to verify (\ref{eq:gammaeta}).

For the proof of (\ref{eq:gamma12}), we use the isomorphism
(\ref{eq:kerbeta1}) and the definition of
the bracket $\{\cdot, \cdot \}_\Pi$, the details of which
we omit. This finishes the proof.
\end{proof}

Using this proposition, if we choose any smooth function $a, \, b \in
C^\infty(S^2_\alpha)$ with $\{a, b\}_{S^2_\alpha} \neq 0$
and set $\Gamma =
(a\circ \pi)f_1^* + (b\circ \pi) f_2^*$, then we derive
$$
\int_\FF\{\Gamma,\Gamma\}_\Pi = \{a, b\}_{S^2_\alpha} \neq 0.
$$
Here $\pi: Y_\alpha \to Y_\alpha/\sim$ is the obvious projection.
Therefore the cohomology class $\alpha \in H^1(Y,\omega)$ is
obstructed. We summarize the above discussion into the following
theorem.

\begin{thm}\label{obstructed}
Let $\alpha > 1$ be rational. Then the coisotropic submanifold
$Y_\alpha \subset \R^6$ or equivalently the presymplectic manifold
$(Y_\alpha, \omega_\alpha)$ is obstructed.
\end{thm}

\begin{rem} In the irrational case of $\alpha$, the flow of $X_{H_\alpha}$
does not generate a circle action, and there is no simple analog of the
second projection
$$
\pi_2: Y_\alpha \to Y_\alpha/S^1_{H_\alpha}.
$$
Analysis of the first cohomology $H^1(Y_\alpha, \omega_\alpha)$ for the irrational
case requires
a fair amount of general foliation theory, whose study we will postpone to
a separate paper elsewhere.
\end{rem}

\section{Appendix: Description in  the super or graded language}
\label{sec:appendix}

In this appendix, we will provide a more physical description of
our deformation problem in the context of Batalin-Vilkovisky
formalism of supermanifolds. Because of this, we will not attempt
to make our discussion completely rigorous in the mathematical
sense in this appendix. However most of the discussions except few
explicitly stated conjectures can be made mathematically rigorous,
which we postpone to a future work.  We already gave one such
example in the proof of Theorem \ref{algebroid} in section 8.

We will also describe the formal deformation space
of the strong homotopy Lie algebra $\frak l^\infty_{(Y,\omega;
\Pi)}$. We refer readers to [OP] for more explanation of the BV
formalism and for an off-shell description of the $A$-model of
topological sigma models. We refer to [AKSZ] for a similar
approach to the closed $A$-model (the Gromov-Witten theory) on
K\"ahler manifold.

We shall interpret general open string $A$-model as a machine to
quantize the algebra of functions on $\L=\Pi E$ as an
$A_\infty$-algebra or, equivalently, as an $1$-algebra in the
sense of Kontsevich [K2]. First we review the basic set up for
the quantization of the $1$-algebra presented as in [OP], following
[P1], [P2]. Then we shall show how a coisotropic submanifold
naturally arises as the general boundary condition
for the open string $A$-model
and that deformations of the coisotropic submanifold correspond to
boundary deformations of the open string $A$-model.

Let $\mathbb{L}$ be a smooth $\Z$-graded space over $\C$ and let
$\mathfrak{l}$ be the super-commutative ring  of functions on
$\mathbb{L}$. As an abelian group we have a direct sum
decomposition $\mathfrak{l}=\oplus_n {\mathfrak{l}}_n$, where
$\mathfrak{l}_n$ is the maximal subspace of ${\mathfrak{l}}$
consisting of degree $n$ functions.  In this appendix, we will use
the complex {\it without shift of grading}. Then we have
super-commutative (and associative) product with degree $0$ such
that, for homogenous element $\g_1,\g_2 \in {\mathfrak{l}}$
$$
\g_1\g_2 -(-1)^{|\g_1||\g_2|}\g_2\g_1 =0,
$$
where $|\g|$ denotes the parity of $\g$ defined by
$degree(\g)\hbox{ mod } 2$. Let $T^*[1]\mathbb{L}$ be the total
space of the twisted by degree $1$ cotangent bundle to
$\mathbb{L}$ and let $\mathfrak{t}$ be the super-commutative (and
associative) ring of functions on it. We have a direct sum
decompositions $\mathfrak{t}=\oplus_\ell{\mathfrak{t}}_\ell$ of
${\CC}$ as an abelian group and super-commutative (and
associative) product with degree $0$. Note that there is an odd
symplectic structure $\Omega$ on $T^*[1]\mathbb{L}$ induced from
the canonical symplectic structure of the total space of tangent
bundle $T^*\mathbb{L}$ to $\mathbb{L}$

\begin{defn}
The pair $(\mathfrak{t}, [\bullet,\bullet])$,
where $[\bullet,\bullet]:{\mathfrak{t}}_{\ell_1}\otimes
{\mathfrak{t}}_{\ell_2}\rightarrow {\mathfrak{t}}_{\ell_1+\ell_2 -1}$
is the graded Poisson bracket of degree $-1$ induced from the
canonical symplectic form $\Omega$ of degree $1$ on
$T^*[1]\mathbb{L}$, is called the structure of
symplectic  $2$-algebra on $T^*[1]\mathbb{L}$ or on
$\mathfrak{t}$.
For any homogeneous elements $A,B,C \in \mathfrak{t}$
\begin{enumerate}
\item super-commutativity
$$
[A,B]=-(-1)^{(|{A}|+1)(|{B}|+1)}[B,A]
$$
\item super-Jacobi
$$
[A,[B,C]]=[[A,B],C]+(-1)^{(|{A}|+1)(|{B}|+1)}[B,[A,C]]
$$
\item  super-Leibnitz
$$
[A,B\cdot C]=[A,B] C +(-1)^{(|A|+1)|B|}B[A,C]
$$
\item linearity
$$
[A, B+C]= [A,B] +[A,C]
$$
\end{enumerate}
\end{defn}

\begin{cor}
The  symplectic $2$-algebra $(\mathfrak{t}, [\bullet,\bullet])$ is
a Gerstenhaber algebra
\end{cor}

\begin{cor}
The bracket $[\bullet,\bullet]$, after forgetting the product, induces
a structure of Lie algebra on $\mathfrak{t}_1$;
$$
[\bullet,\bullet]:\mathfrak{t}_1\otimes\mathfrak{t}_1\longrightarrow \mathfrak{t}_1.
$$
\end{cor}

{}From the corollary above we have degree preserving adjoint action
$ad_{(-B)}(A)=[A, B]$
by an element $B \in\mathfrak{t}_1$ on any homogeneous
element $A \in\mathfrak{t}$ and the associated transformation
$$
e^{ad_{(-B)}}(A) := A + [A,B] +\frac{1}{2!}[[A, B],B] + \ldots,
$$
The above does not make sense unless the Lie algebra $\mathfrak{t}^1$
is nilpotent, otherwise we may tensor it with suitable Artinian ring [Ge]. Then
$$
\left[e^{ad_{(-B)}}(A_1),e^{ad_{(-B)}}(A_2)\right]
=e^{ad_{(-B)}}(\left[A_1,A_2\right]).
$$
We remark that the adjoint action $e^{ad_{B}}$, for $B \in \mathfrak{t}_1$
is equivalent to a degree
preserving canonical transformation connected to the identity.

Note that any odd element in $\mathfrak{t}$ automatically commutes
with itself on the bracket, while the condition that an even element
commutes with itself on the bracket is non-trivial.

We say there is a structure of {\it weak homotopy Lie
$1$-algebroid} on $\mathfrak{l}$ or on $\mathbb{L}$ if there exist
a non-vanishing element $H \in\mathfrak{t}_2$ satisfying $[H, H]=0$.
Two structures $H$ and $H^\pr$ of weak homotopy Lie $1$-algebroid
on $\mathbb{L}$ are defined to be equivalent if there exists some
$B \in \mathfrak{t}^1$ such that $H^\pr = e^{ad_B}(H)$.

Consider a natural $\C^*$ action of weight $1$ on the fiber, over
the zero section $\mathbb{L}$, of $T^*[1]\mathbb{L}$ such that the
degree $1$ symplectic form $\O$ on  $T^*[1]\mathbb{L}$ has weight
$1$. For given $H$ above, we may expand it as
$$
H = \sum_{\ell=0}^\infty H_\ell
$$
in the neighborhood of $\mathbb{L}$ according to the integral
weight $n$ of the $\C^*$ action. We may identify $H_0$ above as
the restriction $H|_\mathbb{L}$ of $H$ to $\mathbb{L}$.

A structure of {\it strong homotopy Lie $1$-algebroid} on
$\mathbb{L}$ or on $\mathfrak{l}$ is defined by an element $H \in
\mathfrak{t}_2$ satisfying
$$
[H, H]=0,\qquad H|_\mathbb{L}=0.
$$
Consider a structure $H$ of strong homotopy Lie $1$-algebroid on
$\mathbb{L}$ which has a decomposition
$$
H = \sum_{\ell=1}^\infty H_\ell
$$
according to the weights of the representation of $\C^*$ mentioned
above. The equation $[H,H]=0$ has the corresponding decompositions
\begin{eqnarray}\label{eq:xxx}
[H_1,H_1]&=&0,\cr [H_1,H_2]&=&0,\cr \frac{1}{2}[H_2,H_2]+ [H_1,
H_3]&=&0,\cr \vdots&&
\end{eqnarray}
etc. For a given $H_n$ in the sequence $(H_1,H_2,\ldots)$ of above
we can associate $n$-multi-linear map  $\mm_n$ of degree $2-n$;
$$
\mm_n :\mathfrak{l}^{\otimes n} \longrightarrow \mathfrak{l}
$$
by, for any set $\g_1,\ldots,\g_n$ of  homogenous elements of
$\mathfrak{l}$
$$
\mm_n(\g_1,\ldots,\g_n):=[[\cdots[H_n,\g_1],\cdots],\g_n]
$$
Then the relation (\ref{eq:xxx}) together with the super-Jacobi
identity of the bracket $[\bullet,\bullet]$ implies that
$(\mm_1,\mm_2,\mm_3,\ldots)$ satisfies the relation equivalent to
that of strong homotopy Lie algebra ($L_\infty$-algebra in short).

\begin{rem}
Our definition includes the standard notion
of the Lie algebroid as a special case, as shown
below.
Consider a smooth manifold $M$ and a vector bundle $E$ over $M$.
Let $\mathbb{L}=\Pi E$ be the total space of $E$ after applying
twisting functor $\Pi$ by degree $1$ to the fiber. Let $\{x^I\}$ be
local coordinates on $X$, $\{e_a\}$ be a local frame on $E$
and $\{e^\a\}$ be the dual frame.  Then we may identify
$\{x^I,c^\a :=\Pi e^\a\}$ as a coordinate system on $\Pi E$.
Consider $T^*[1]\mathbb{L}$ with Darboux coordinates
$\{x^I,c^\a |\chi_I, e_\a\}$ with degree $\{0,1|1,0\}$
and the canonical degree $1$-symplectic structure
$$
\O = d\chi_I dx^I + de_\a d c^\a.
$$
Now consider $H \in \mathfrak{t}_2$ given by
$$
H= c^\a \Gamma(x)_{\a}{}^I\chi_I + \Fr{1}{2}C(x)_{\a\b}{}^\g c^\a
c^\b e_\g
$$
Then the condition $[H,H]=0$ means that $(\Gamma_{\a}{}^I,
C_{\a\b}{}^\g)$ are the anchor and structure function of Lie
algebroid. Note that $H|_\L=0$. Let $\s =s^\a e_\a$ and $\t =\t^\a
e_\a$ by any sections of $E$ and let $g$ be a smooth function on
$M$ then
$$
[[H,\s], g\t] = \s^\a \Gamma_{\a}{}^I \Fr{\rd g}{\rd x^I} + g\cdot
C_{\a\b}{}^\g \s^\a\t^\b e_\g.
$$
The graded supercommutative algebra $\mathfrak{l}$ functions on
$\mathbb{L}$ is isomorphic the exterior algebra of E differential
forms. The corresponding E differential operator is given by $\CQ
=[H,\ldots]$ restricted to $\L$;
$$
\mm_1\equiv\CQ_1|_\L = c^\a \Gamma_{\a}{}^I\Fr{\rd}{\rd x^I} +
\Fr{1}{2}C_{\a\b}{}^\g c^\a c^\b\Fr{\rd}{\rd c^\g}.
$$
Note that $H = H_1$ and, thus, $\mm_2=\mm_3=\ldots = 0$.
\end{rem}

Two structures $H$ and $H^\prime$ of strong homotopy Lie
$1$-algebroid on $\mathbb{L}$ are equivalent if they are related
by the adjoint action of $\b \in \mathfrak{t}_1$ satisfying
$\b|_\mathbb{L}=0$. Note that such an adjoint action preserves the
conditions $H|_\mathbb{L}=H^\prime |_\mathbb{L}=0$, and equivalent
to a change of Lagrangian complimentary in $T^*[1]\L$.

Let $\G \in \mathfrak{l}_1\subset\mathfrak{t}_1$ and let
$H^\G\equiv e^{ad_\G}(H)$ denotes the resulting canonical
transformations of $H$. Then it is obvious that $[H^\G,H^\G]=0$,
while $H^\G |_\mathbb{L} \neq 0$ in general. For each $\G$ leading
to $H^\G|_\mathbb{L} =0$ we have another structure of strong
homotopy Lie $1$-algebroid on $\L$. The condition
$H^\G|_\mathbb{L}=0$ is equivalent to Maurer-Cartan equation \be
\label{eq:bmas} \sum_{\ell=1}^\infty
\frac{1}{\ell!}
\mathfrak m_\ell(\G,\ldots,\G) =0. \ee

Now we explain how the above general story is relevant to
our subject.

We consider a symplectic manifold $(X,\o_X)$ and the graded space
$T[1]X$, which is the total space of twisted by the degree $1$
tangent bundle to $X$ where $U$ is the degree or the {\it ghost
number}. Let $\mathfrak{t}=\oplus_{k=0}^{2n}\mathfrak{t}_k$ be the
$\Z$-graded supercommutative algebra of smooth functions on
$T[1]X$. The algebra $\mathfrak{t}$ is isomorphic to the exterior
algebra  of differential forms on $X$, such that
$\Omega^k(X)\simeq \mathfrak{t}_k$ and the wedge product is
replaced with the supercommutative product. The exterior
derivative induces an odd degree $1$ vector field $\CQ$ on $T[1]X$
such that $\CQ:\mt_k \rightarrow \mt_{k+1}$ and $\CQ^2=0$. Thus we
have a structure of differential graded algebra $(\mathfrak{t},
\CQ, \cdot)$. The cohomology of $(\mathfrak{t}, \CQ, \cdot)$ is
isomorphic to de Rham cohomology of $X$.  The symplectic structure
$\o_X$ on $X$ induces a degree $1$ (odd) symplectic form $\O$ on
$T[1]X$ via the standard isomorphism $TX \rightarrow T^*X$
together with the twisting. Thus we have a structure of symplectic
$2$-algebra on $T[1]X$ by the pair $(\mathfrak{t},
[\bullet,\bullet])$, we the degree $-1$ (odd) Poisson bracket is
defined by $\O$. The above bracket is equivalent to Koszul
bracket, which is the covariant version Schouten-Nijenhuis bracket
[Ko].

There exists an element $H \in \mathfrak{t}_2$, which is
isomorphic to the symplectic form $\o_X$.
Using the closedness of $\omega$, it is not difficult to check the
following identities
$$
[H, H] =0,\qquad\CQ = [H,\phantom{\bullet}].
$$
Thus  $\CQ$, isomorphic to
the exterior derivative, is realized as the Hamiltonian vector
field of the function $H$.

\begin{rem}

It is instructive to give a coordinate representation of above. We
introduce a local coordinates $\{x^I\}$, $I=1,\ldots, 2n$, on $X$.
We denote the corresponding fiber coordinates on $T[1]X$ by
$\{\p^I\}$ carrying the degree $1$. Then $\CQ = \p^I\Fr{\rd}{\rd
x^I}$. In the sense of ordinary geometry, $\psi^I$ is nothing but
$dx^I$ considered as a fiberwise linear function on $T[1]X$. Now
the symplectic structure $\o= \Fr{1}{2}\o_{IJ}dx^I\wedge dx^J$ on
$X$ induces a non-degenerate function
$$
H=\Fr{1}{2}\o_{IJ}\p^I\p^J \in \mt_2
$$
of degree $2$ on $T[1]X$. Thus the condition $d\o_X=0$ is equivalent to
$\CQ H=0$. The
symplectic structure $\o_X$ on $X$ induces a degree $U=1$ symplectic
form
$$
\O = d(\o_{IJ}\p^I) dx^J
$$
on $T[1]X$ and such that the corresponding graded Poisson bi-vector
on $T[1]X$ is given by
$$
\O^* = [(\o^{-1})^{IJ}\Fr{\rd}{\rd \p^I}\wedge \Fr{\rd}{\rd x^J}
+\Fr{1}{2}\left(\Fr{\rd (\o^{-1})^{IJ}}{\rd x^K}\right)
\p^K\Fr{\rd}{\rd \p^I}\wedge\Fr{\rd}{\rd \p^J}
$$
which define the  graded Poisson bracket $[\bullet,\bullet]$.
One may check $[H,H]=0$ and $\CQ=[H,\bullet]$ by an explicit computation.

\end{rem}

Let $\mathbb{L}$ be a Lagrangian subspace of $(T[1]X,\O)$. Then
the following is easy to see from the form of the Hamiltonian $H$
in coordinates.
\begin{lem}
Any (conic) Lagrangian subspace $\mathbb{L}$ satisfying
$H|_\mathbb{L}=0$ is equivalent to $\Pi E =\Pi(TY)^\o$, where $Y$
is a coisotropic submanifold of $X$ and $\Pi E$ is the total space
of $E$ after twisting the fiber by $1$.
\end{lem}
Consequently each coisotropic submanifold $Y\subset X$ inherits a
structure of strong homotopy Lie $1$-algebroid induced from the
symplectic structure of $X$, and the (formal) deformation problem
of coisotropic submanifold is equivalent to that of strong
homotopy Lie $1$-algebroid. We remark that $\mathfrak{l}=\oplus
\mathfrak{l}_\bullet$ is isomorphic to $\oplus\G( \wedge^\bullet
E^*)$ and in particular $\G \in \mathfrak{l}_1$ is isomorphic to
$\G( E^*)$. The condition $H^\G|_\mathbb{L} =0$ can be identified
with $H|_{\mathbb{L}^\G}=0$ where $\mathbb{L}^\G$ is the graph of
Lagrangian subspace generated by $\G$. Thus the condition
$H^\G|_\mathbb{L} =0$ means $\mathbb{L}^\G$ is given by $\Pi
E^\prime \equiv\Pi(TY^\prime)^\o$, where $Y^\prime$ is another
coisotropic submanifold.

It is also natural to consider extended deformations of
coisotropic submanifold. Let us consider a graded Artin ring with
maximal ideal $\mathfrak{a}$;
$$
\mathfrak{a} = \bigoplus_{-(n-k)\leq j \leq 1}\mathfrak{a}_j
$$
where $(n-k)$ is the rank of $E=(TY)^\o$. Let $\Upsilon \in
(\mathfrak{l}\otimes \mathfrak{a})_1$. Then the condition
$$
H^\Upsilon|_\mathbb{L}\equiv ad_\Upsilon(H)
\equiv\sum_{n=1}^\infty\frac{1}{\ell!}
m_\ell(\Upsilon,\ldots,\Upsilon)=0
$$
may be regarded as the condition for the extended deformations of
coisotropic submanifold $Y$.

We note that the deformation problem of coisotropic submanifold
is obstructed in general.
On the other hand our setting allows us to consider an
extended deformation problem  allowing both $\Gamma$ and $H$
to vary. It may turn out that this extended deformation problem could
be unobstructed. This is a subject of future study.

Now we like to motivate the problem of quantization of
coisotropic submanifold.
We begin with recalling a lemma of Kontsevich [K2] stating
\begin{quote}
The cohomology of the Hochschild complex of the algebra of
functions $A_1(\L)$ on $\L$, regarded as an $1$-algebra, is
isomorphic to the space $\mathfrak{t}$ of functions on
$T^*[1]\mathbb{L}$.
\end{quote}
In the spirit of the above lemma, we may regard a structure $H$ of
strong homotopy Lie $1$-algebroid is an element of the $1$st
cohomology of Hochschild complex of $A_1(\L)$, which is the $1$st
order deformations of $A_1(\mathbb{L})$ as a strong homotopy
associative ($A_\infty$ in short) algebra (or simply as an
$1$-algebra). In the spirit of Kontsevich's formality theorem
[K1], we may define the quantization of the $1$-algebra
$A_1(\mathbb{L})$ as a quasi-isomorphism between the Hochschild
complex of $A_1(\mathbb{L})$ and its cohomology.

Open string or quantum field theory in a two-dimensional manifold
with boundary may be regarded as a universal machine to quantize
the $1$-algebra. The following can be shown
\begin{quote}
For each structure $H$ of the strong homotopy Lie $1$-algebroid on
any $\L$ there exists a pre-quantum field theory in two-dimension with
boundary, whose boundary condition is defined in terms of
$\mathbb{L}$, depending on the topology of two-manifolds, so that
it satisfies the classical BV master equation.
\end{quote}
In particular the set of solutions of (\ref{eq:bmas}), modulo
equivalence, is isomorphic to the moduli space of boundary
interactions.

The second named author [P2] called such a QFT as an {\it open
$1$-braneoid}. Assume that the resulting QFT actually satisfies
the quantum BV master equation. Then a conjecture is that the path
integral generates a quasi-isomorphism between the Hochschild
complex of $A_1(\mathbb{L})$ and its cohomology.

Now we may interpret the open string $A$-model with the
coisotropic boundary condition as a machine to quantize the
algebra of functions on $\mathbb{L}=\Pi E$ as an
$A_\infty$-algebra. We may simply call the problem as {\it
quantization of a coisotropic submanifold} $Y$ on a symplectic
manifold $X$. In [OP] we shall see that the genus zero open string
$A$-model governs the maps $\Phi: T[1]D\rightarrow T[1]X$
satisfying $\Phi(T[1](\rd D)) \subset \mathbb{L}=\Pi E$. Then the
path integral is formally defined as an integral over an (infinite
dimensional) Lagrangian subspace of the space all maps $\Phi$
determined by conformal structure on $D$ and an almost complex
structure on $X$.
\medskip

\begin{exm} Let $Y=X$, then $\mathbb{L} = X$ and $m_\ell\neq
0$ only for $\ell=2$ such that $m_2$ is the usual Poisson bracket
on $X$, and consider the zero-instanton sector governed by
constant maps of the $A$-model. Then the quantization problem is
nothing but that of $X$ as a Poisson manifold. This is the
original context of Kontsevich's formality theorem [K1], as
interpreted in the path integral approach of [CF1]. What is the
result after including instanton corrections even in this case?
\end{exm}

We conjecture, which can be justified at the physical level of
rigor, that the Hochschild complex of $A_1(\mathbb{L} =\Pi E)$ is
quasi-isomorphic to its cohomology. The above conjecture is based
on the following fact.
\begin{quote}
There exist a degree $-1$ odd differential operator $\Delta$;
$\Delta :\mathfrak{t}_\ell \rightarrow \mathfrak{t}_{\ell-1}$
satisfying $\Delta^2=0$ and generates the bracket
$[\bullet,\bullet]$ such that
$$
\Delta H =\Delta H^\Upsilon =0.
$$
\end{quote}
The above properties can be used to formally show that the open
string $A$-model with coisotropic boundary condition actually
satisfies the quantum BV master equation, which implies that the
path integral would give the formality map. The rigorous proof of
this assertion will be investigated in the future.
\medskip

Now we go back to our study of coisotropic submanifolds using this
BV formalism and give a proof of Theorem \ref{algebroid} in this
context. We
introduce a local coordinates system $(y^i, q^\a| p_\a)$ of $X$ in
a neighborhood of $Y$ as in section 4. It was shown in section $4$
that the symplectic form $\o$ in $X$ can be written as
\begin{eqnarray}\label{eq:omega}
\omega & = & \frac{1}{2}\Big(\omega_{ij} - p_\beta F^\beta_{ij}
\Big)
dy^i \wedge dy^j \nonumber \\
& \quad &  - (dp_\beta + p_\beta \frac{\del R_i^\beta}{\del
q^\gamma} dy^i) \wedge (dq^\delta - R_j^\delta dy^j)
\end{eqnarray}
It is straightforward to derive
\begin{eqnarray*}
\widetilde\omega\Big(\frac{\del}{\del q^\delta}\Big) & = &
dp_\delta +
p_\beta \frac{\del R_i^\beta}{\del q^\delta}dy^i \\
\widetilde\omega\Big(\frac{\del}{\del p_\delta}\Big) & = &
-(dq^\delta -
R_j^\delta dy^j) \\
\widetilde\omega\Big(\frac{\del}{\del y^i}\Big) & = & (\omega_{ij}
- p_\beta F^\beta_{ij}) dy^j - p_\beta \frac{\del R_i^\beta}{\del
q^\delta}(dq^\delta - R_j^\delta dy^j) \nonumber \\
& \quad & - R_i^\delta(dp_\delta + p_\beta \frac{\del
R_i^\beta}{\del q^\delta}dy^i)
\end{eqnarray*}
The last identity can be rewritten as
\begin{equation}
\widetilde\omega(e_i) = \widetilde\o_{ij}dy^j
\end{equation}
where we recall to have defined
$$
e_i =\frac{\del}{\del y^i} + R_i^\delta\frac{\del}{\del q^\delta}
- p_\beta\frac{\del R_i^\beta}{\del q^\delta}\frac{\del}{\del
p_\delta}.
$$

 Combining these we have derived the following formula for
the inverse $(\widetilde\o)^{-1}: T^*X \to TX$:
\begin{eqnarray}\label{eq:inverse}
\begin{cases}\widetilde\o^{-1}(dy^i)& =
\widetilde\o^{ij}e_j \\
\widetilde\o^{-1}(dp_\delta) & = \frac{\del}{\del q^\delta} -
p_\beta\frac{\del R_i^\beta}{\del q^\delta} \widetilde\o^{ij}
e_j \\
\widetilde\o^{-1}(dq^\delta) & =-\frac{\del}{\del p_\delta} -
R_j^\delta\widetilde\o^{ij} e_j
\end{cases}
\end{eqnarray}
In the canonical coordinates of $T[1]X$ associated $(y^i,
q^\alpha\mid p_\alpha)$, the even Hamiltonian $H$ becomes
\begin{equation}\label{eq:H}
H = \frac{1}{2}(\o_{ij} - p_\beta F^\beta_{ij})\psi^i\psi^j +
(\eta^\delta - R_j^\delta \psi^j)(\chi^\delta + p_\beta\frac{\del
R_i^\beta}{\del q^\delta}\psi^i)
\end{equation}
It would be more convenient to write down the Hamiltonian on
$T^*\L \cong T[1]X$ near $\L = (TY)^\omega[1] \subset T[1]X$ using
the super version of Darboux-Weinstein theorem. We denote the
(super) canonical coordinates of $T^*[1]\L$ associated with
$(y^i,q^\alpha\mid p_\alpha)$ by
$$
\Big(\begin{matrix}y^i, & q^\alpha  && \mid p_*^\alpha \\
y^*_i, & q^*_\alpha && \mid p_\alpha
\end{matrix}\Big)
$$
Here we note that the degree of $y^i,\, q^\alpha$ and $p_\alpha$
are $0$ while their anti-fields, i.e., those with $*$ in them have
degree $1$. And we want to emphasize that $\L$ is given by the
equation
\begin{equation}\label{eq:LL}
y_i^* = p_\alpha = p_*^\alpha = 0
\end{equation}\label{eq:LL}
and $(y^i, y^*_i)$, $(p_\alpha, q_\alpha^*)$ and $(p_*^\alpha,
q^\alpha)$ are conjugate variables.
Then we have the canonical odd symplectic form of degree $1$
and associated canonical odd Poisson bracket $[\bullet,\bullet]^*$
of degree $-1$.

It follows from (\ref{eq:inverse}) that we have
\begin{eqnarray*}
\widetilde\omega^{-1}(dq^\delta - R_j^\delta dy^j) & = &-
\frac{\del}{\del p_\delta} \\
\widetilde\omega^{-1}(dp_\delta + p_\beta \frac{\del
R_j^\beta}{\del q^\delta} dy^j) & = & \frac{\del}{\del q^\delta}
\end{eqnarray*}
and so that $H$ has the form
\begin{eqnarray}
H  =  \frac{1}{2} \widetilde\o^{ij}y_i^\# y_j^\# + p_*^\delta
q^*_\delta  \label{eq:H*}
\end{eqnarray}
in the canonical coordinates of $T^*[1]\L$. Here we define
$y_i^\#$ to be
$$
y_i^\#: = y_i^* + R_i^\delta p^\delta_* - p_\beta \frac{\del
R_i^\beta}{\del q^\delta}q^*_\delta.
$$
We then derive
$$
\widetilde \omega^{-1}_\alpha = \omega^{-1}_{\pi(\alpha)}
\sum_{\ell =0}^\infty (F^\# \rfloor \alpha)^\ell \quad \mbox{on }
\, TY/E
$$
(See (\ref{eq:omegaF}) in section \ref{sec:moduli} later.) which
is written as
$$
\widetilde\omega^{ij}_\alpha = \omega^{ij_0}_{\pi(\alpha)}
\sum_{\ell=0}^\infty
(p_{\beta_1}F^{\beta_1j_1}_{j_0})(p_{\beta_2}F^{\beta_2j_2}_{j_1})
\cdots (p_{\beta_\ell}F^{\beta_\ell j_\ell}_{j_{\ell-1}})
$$
in coordinates where $\alpha = p_\beta f_\beta^*$. Since $H$
vanishes on $\L$, its Hamiltonian vector field $\CQ$ is tangent to
$\L$ and so can be canonically restricted to $\L$ as an odd vector
field on $\L$. The odd vector field
\begin{equation}\label{eq:q}
m_1 = \CQ \mid_\L
\end{equation}
acts on
$$
\frak l = \bigoplus_{\ell = 0}^{n-k}\frak l^\ell \cong
\bigoplus_{\ell = 0}^{n-k}\Omega^\ell(\FF)
$$
as the corresponding directional derivative and equivalent to
$d_{\FF}$.

The general discussions mentioned above then implies that the
strong homotopy Lie algebroid structure is obtained by expanding
$H$ into the series
$$
H = \sum_{\ell = 1} H_\ell, \quad H_\ell \in \frak l^\ell
$$
in the normal direction of $\L$, i.e., in terms of $(y_i^*,
p_*^\alpha, p_\alpha)$ and use the odd Poisson bracket
$[\bullet,\bullet]$.

Now Theorem \ref{algebroid} immediately follows
from our general discussions.

\end{document}